\documentclass[a4paper,12pt,amsart,frenchb]{article}

\usepackage{amsmath,amsbsy,amsfonts,amssymb}
\usepackage[french]{babel}
\newcommand{\ass}[2]{\vskip0.3cm\noindent
{\bf {#1}}. { \sl {#2}}\vskip0.3cm\noindent
}
  
\begin{document}

    \title{ Stabilisation de la formule des traces tordue VIII:  l'application $\epsilon_{\tilde{M}}$ sur un corps de base local non-archim\'edien}
\author{ J.-L. Waldspurger}
\date{3 octobre  2014}
\maketitle

\bigskip

{\bf Introduction.} 

\bigskip

Dans une s\'erie d'articles, nous avons \'enonc\'e les th\'eor\`emes qui conduisent \`a la stabilisation de la partie g\'eom\'etrique de la formule des traces tordue. Le th\'eor\`eme cl\'e est local. Rappelons-le tr\`es sommairement. Le corps de base $F$ est local de caract\'eristique nulle. On consid\`ere un groupe r\'eductif connexe $G$ d\'efini sur $F$, un espace tordu $\tilde{G}$ sous $G$ et une classe de cocycle ${\bf a}\in H^1(W_{F};Z(\hat{G}))$, auquel est associ\'e un caract\`ere $\omega$ de $G(F)$. On suppose $\omega$ unitaire. Soit $\tilde{M}$ un espace de Levi de $\tilde{G}$. Pour un \'el\'ement $\gamma\in \tilde{M}(F)$ qui est fortement r\'egulier dans $\tilde{G}(F)$ et pour une fonction $f\in C_{c}^{\infty}(\tilde{G}(F))$, on d\'efinit l'int\'egrale orbitale pond\'er\'ee $\omega$-\'equivariante $I_{\tilde{M}}^{\tilde{G}}(\gamma,\omega,f)$ (la d\'efinition exacte n\'ecessite d'introduire des mesures dont nous ne tenons pas compte dans cette introduction). Dans l'article [II], nous avons d\'efini un avatar endoscopique de ce terme, que l'on note $I_{\tilde{M}}^{\tilde{G},{\cal E}}(\gamma,\omega,f)$. Le th\'eor\`eme principal affirme que, pour tous $\gamma$ et $f$ comme ci-dessus, on a l'\'egalit\'e
$$I_{\tilde{M}}^{\tilde{G},{\cal E}}(\gamma,\omega,f)=I_{\tilde{M}}^{\tilde{G}}(\gamma,\omega,f).$$

Nous supposons ici que $F$ est non-archim\'edien. Sous cette hypoth\`ese, nous effectuons le premier pas dans la d\'emonstration du th\'eor\`eme. Il consiste \`a prouver que pour $f\in C_{c}^{\infty}(\tilde{G}(F))$, il existe une fonction $\epsilon_{\tilde{M}}(f)$ sur $\tilde{M}(F)$ v\'erifiant la condition suivante. Pour tout  \'el\'ement $\gamma\in \tilde{M}(F)$ qui est fortement r\'egulier dans $\tilde{G}(F)$, la diff\'erence entre les deux termes dont nous voulons prouver l'\'egalit\'e est \'egale \`a l'int\'egrale orbitale de $\epsilon_{\tilde{M}}(f)$ au point $\gamma$. Cette derni\`ere est une int\'egrale orbitale    sur $\tilde{M}$,  ordinaire c'est-\`a-dire non pond\'er\'ee, mais tenant \'evidemment compte du caract\`ere $\omega$. Autrement dit
$$(1) \qquad I^{\tilde{M}}(\gamma,\omega,\epsilon_{\tilde{M}}(f))=I_{\tilde{M}}^{\tilde{G},{\cal E}}(\gamma,\omega,f)-I_{\tilde{M}}^{\tilde{G}}(\gamma,\omega,f).$$
La suite de la d\'emonstration consistera \`a appliquer la formule des traces dans $\tilde{M}$ \`a cette fonction $\epsilon_{\tilde{M}}(f)$  et \`a en d\'eduire que  celle-ci est nulle. Cela d\'emontrera le th\'eor\`eme principal. Au point o\`u nous en sommes, nous pouvons seulement prouver l'existence d'une telle fonction $\epsilon_{\tilde{M}}(f)$. Nous ne pouvons m\^eme pas prouver que l'on peut la choisir localement constante et \`a support compact. Nous prouvons toutefois qu'on peut lui imposer les deux conditions

- il existe un sous-groupe ouvert compact de $M(F)$ tel que $\epsilon_{\tilde{M}}(f)$ soit biinvariante par ce sous-groupe;

- la restriction de $\epsilon_{\tilde{M}}(f)$ \`a toute fibre de l'application usuelle $\tilde{H}_{\tilde{M}}:\tilde{M}(F)\to \tilde{{\cal A}}_{\tilde{M}}$ (cf. 1.1) est \`a support compact.  

Ces conditions suffisent pour d\'efinir les int\'egrales orbitales $I^{\tilde{M}}(\gamma,\omega,\epsilon_{\tilde{M}}(f))$. En fait, nous d\'emontrons plus que la seconde propri\'et\'e ci-dessus: la fonction $\epsilon_{\tilde{M}}(f)$ est de Schwartz. Nous d\'efinissons en 1.7 ce que nous entendons par l\`a (on trouve dans la litt\'erature d'autres d\'efinitions des fonctions de Schwartz, qui ne co\"{\i}ncident sans doute pas avec la n\^otre). 

Il y a une restriction \`a notre r\'esultat. Pour certains triplets $(G,\tilde{G},{\bf a})$, on impose que les int\'egrales orbitales ordinaires de $f$ ($\omega$-\'equivariantes) sont nulles en tout point de la r\'eunion d'un certain ensemble fini de classes de conjugaison  (cf. 4.4 pour un \'enonc\'e pr\'ecis). 

Notre r\'esultat est l'exact analogue dans le cas tordu de la proposition 3.1 de [A1], dont nous reprenons la preuve. En utilisant les r\'esultats sur les germes de Shalika prouv\'es en [II] et [III], il est assez bref de d\'emontrer que, pour tout \'el\'ement semi-simple $\eta\in \tilde{M}(F)$, il existe une fonction $\epsilon_{\tilde{M}}(f)$ lisse et \`a support compact de sorte que (1) soit v\'erifi\'e pour $\gamma$ au voisinage de $\eta$. On a envie ensuite de recoller les fonctions ainsi construites en utilisant une partition de l'unit\'e. Mais on ne peut contr\^oler ni l'uniforme lissit\'e de la fonction ainsi construite, ni sa croissance \`a l'infini. On a besoin d'une deuxi\`eme construction qui, elle, nous fournit une fonction $\epsilon_{\tilde{M}}(f)$ qui a les bonnes propri\'etes de lissit\'e et de croissance et qui v\'erifie l'\'egalit\'e (1), cette fois pour $\gamma$ hors d'un certain compact. On arrive \`a bon port en utilisant \`a la fois les deux constructions. Cette deuxi\`eme construction utilise des variantes des int\'egrales orbitales pond\'er\'ees $\omega$-\'equivariantes, que l'on note $^cI_{\tilde{M}}^{\tilde{G}}(\gamma,\omega,f)$ en suivant comme toujours Arthur. La magnifique propri\'et\'e de ces termes est que, pour $f$ fix\'ee, ils sont \`a support compact en $\gamma$, modulo conjugaison. On doit d'abord d\'efinir et \'etudier ces termes, ainsi que  certaines applications n\'ecessaires \`a leur d\'efinition. C'est l'objet de la section 1. On doit ensuite les stabiliser (section 2) et en d\'efinir des avatars endoscopiques (section 3). La d\'efinition de l'application $\epsilon_{\tilde{M}}$ et la preuve de ses propri\'et\'es est donn\'ee dans la derni\`ere section. 

\bigskip

\section{L'application $^c\theta_{\tilde{M}}$}

\bigskip

\subsection{D\'efinition de fonctions combinatoires}
Dans tout l'article, le corps de base $F$ est local non archim\'edien et de caract\'eristique nulle. On consid\`ere un triplet $(G,\tilde{G},{\bf a})$ d\'efini sur $F$. On suppose que le caract\`ere $\omega$ associ\'e \`a ${\bf a}$ est unitaire. 

Soit $\tilde{M}$ un espace de Levi de $\tilde{G}$. On note $\Sigma(A_{\tilde{M}})$ l'ensemble des racines de $A_{\tilde{M}}$ dans $G$.  Un \'el\'ement $\alpha\in \Sigma(A_{\tilde{M}})$ peut \^etre consid\'er\'e comme un \'el\'ement de ${\cal A}_{\tilde{M}}^*$. Il lui est associ\'e une coracine $\check{\alpha}\in {\cal A}_{\tilde{M}}$. Sa d\'efinition pr\'ecise est un peu arbitraire, l'ensemble $\Sigma(A_{\tilde{M}})$ n'\'etant pas en g\'en\'eral un syst\`eme de racines au sens de Bourbaki. Toutefois, la demi-droite port\'ee par $\check{\alpha}$ est d\'efinie sans ambigu\"{\i}t\'e et c'est la seule chose qui nous importera. Tout sous-espace parabolique $\tilde{P}\in {\cal P}(\tilde{M})$ d\'etermine un sous-ensemble positif dans  $\Sigma(A_{\tilde{M}})$ que l'on note $\Sigma^{\tilde{P}}(A_{\tilde{M}})$. On en d\'eduit des chambres positives 
$${\cal A}_{\tilde{P}}^+=\{X\in {\cal A}_{\tilde{M}}; <\alpha,X>>0\,\, \forall \alpha\in \Sigma^{\tilde{P}}(A_{\tilde{M}})\},$$
$${\cal A}_{\tilde{P}}^{*,+}=\{\mu\in {\cal A}_{\tilde{M}}^*; <\mu,\check{\alpha}>>0\,\, \forall \alpha\in \Sigma^{\tilde{P}}(A_{\tilde{M}})\}.$$
On rappelle que, quand $\tilde{P}$ d\'ecrit ${\cal P}(\tilde{M})$, les ensembles ${\cal A}_{\tilde{P}}^+$ d\'ecrivent les composantes connexes de ${\cal A}_{\tilde{M}}$ priv\'e des hyperplans annul\'es par les racines $\alpha\in \Sigma(A_{\tilde{M}})$. 

Dans notre situation tordue, il y  a un espace affine $\tilde{{\cal A}}_{\tilde{M}}$ sur ${\cal A}_{\tilde{M}}$. Rappelons  sa d\'efinition. On note $M(F)^1$ le noyau de l'application $H_{\tilde{M}}:M(F)\to {\cal A}_{\tilde{M}}$ et  ${\cal A}_{\tilde{M},F}$ son image. Notons $\tilde{{\cal A}}_{\tilde{M},F}$ le quotient $M(F)^1\backslash \tilde{M}(F)$. Le groupe ${\cal A}_{\tilde{M},F}$ agit par translations sur ce quotient et celui-ci est  un  espace principal homog\`ene sous  cette action.  On pose $\tilde{{\cal A}}_{\tilde{M}}=\tilde{{\cal A}}_{\tilde{M},F}\times_{{\cal A}_{\tilde{M},F}}{\cal A}_{\tilde{M}}$.  On note $\tilde{H}_{\tilde{M}}:\tilde{M}(F)\to \tilde{{\cal A}}_{\tilde{M}}$ l'application naturelle.

Pour tout $\tilde{P}\in {\cal P}(\tilde{M})$, on fixe une fonction $\omega_{\tilde{P}}:\tilde{{\cal A}}_{\tilde{M}}\to [0,1]$, soumise aux conditions suivantes:

(1) pour tout $X_{0}\in \tilde{{\cal A}}_{\tilde{M}}$,  il existe $c\in {\mathbb R}$ tel que, pour $X\in \tilde{{\cal A}}_{\tilde{M}}$ et $\alpha\in \Sigma^{\tilde{P}}(A_{\tilde{M}})$, la condition $<\alpha,X-X_{0}>< c$ entra\^{\i}ne $\omega_{\tilde{P}}(X)=0$;

(2) $\sum_{\tilde{P}\in {\cal P}(\tilde{M})}\omega_{\tilde{P}}(X)=1$ pour tout $X\in \tilde{{\cal A}}_{\tilde{M}}$.

De telles fonctions existent. Remarquons que ces deux conditions impliquent

(3)  pour tout $X_{0}\in \tilde{{\cal A}}_{\tilde{M}}$,  il existe $c\in {\mathbb R}$ tel que, pour $X\in \tilde{{\cal A}}_{\tilde{M}}$, les conditions $<\alpha,X-X_{0}>>c $ pour tout $\alpha\in \Sigma^{\tilde{P}}({\cal A}_{\tilde{M}})$ entra\^{\i}nent $\omega_{\tilde{P}}(X)=1$.

Les fonctions $\omega_{\tilde{P}}$ sont fix\'ees pour tout espace de Levi $\tilde{M}$. On leur  impose  la condition (4) suivante. Soit $x\in \tilde{G}(F)$. L'automorphisme $ad_{x}$ induit un isomorphisme $\tilde{{\cal A}}_{\tilde{M}}\to \tilde{{\cal A}}_{ad_{x}(\tilde{M})}$. Alors

(4) pour tout $\tilde{P}\in {\cal P}(\tilde{M})$ et tout $X\in \tilde{{\cal A}}_{\tilde{M}}$, $\omega_{ad_{x}(\tilde{P})}(ad_{x}(X))=\omega_{\tilde{P}}(X)$.

C'est possible. En effet, pour tout espace de Levi $\tilde{M}$, notons $W(\tilde{M})$ le quotient $Norm_{G(F)}(\tilde{M})/M(F)$. Ce groupe agit sur ${\cal P}(\tilde{M})$, sur ${\cal A}_{\tilde{M}}$ en permutant les chambres ${\cal A}_{\tilde{P}}^+$ et sur $\tilde{{\cal A}}_{\tilde{M},F}$. Fixons un ensemble $\underline{{\cal L}}$ de repr\'esentants des classes de conjugaison par $G(F)$ d'espaces de Levi. Pour $\tilde{M}\in \underline{{\cal L}}$ et $\tilde{P}\in {\cal P}(\tilde{M})$, 
  on peut remplacer la  fonction $\omega_{\tilde{P}}$ par la fonction
$$X\mapsto \vert W(\tilde{M})\vert ^{-1}\sum_{w\in W(\tilde{M})}\omega_{w\tilde{P}}(wX).$$
Les nouvelles fonctions v\'erifient encore (1) et (2) et de plus $\omega_{w\tilde{P}}(wX)=\omega_{\tilde{P}}(X)$ pour tout $\tilde{P}$, tout $X$ et tout $w\in W(\tilde{M})$. Pour $\tilde{M}'$ quelconque, soit $\tilde{M}$ l'unique \'el\'ement de $\underline{{\cal L}}$ qui est conjugu\'e \`a $\tilde{M}'$ et fixons $x\in G(F)$ tel que $ad_{x}(\tilde{M})=\tilde{M}'$. Pour $\tilde{P}'\in {\cal P}(\tilde{M}')$, on d\'efinit $\omega_{\tilde{P}'}$ par $\omega_{\tilde{P}'}(X')=\omega_{ad_{x}^{-1}(\tilde{P}')}(ad_{x^{-1}}(X'))$. Cela ne d\'epend pas du choix de $x$ et le syst\`eme de fonctions obtenu v\'erifie (4). 

Enfin, on peut faire varier le groupe ambiant $\tilde{G}$. On suppose des fonctions $\omega_{\tilde{P}}$ fix\'ees comme ci-dessus. Soient $\tilde{L}$ un espace de Levi de $\tilde{G}$ et $\tilde{M}$ un espace de Levi de $\tilde{L}$. Il y a une application naturelle
$$\begin{array}{ccc}{\cal P}^{\tilde{G}}(\tilde{M})&\to&{\cal P}^{\tilde{L}}(\tilde{M})\\ \tilde{P}&\mapsto &\tilde{P}\cap \tilde{L}.\\ \end{array}$$
Pour $\tilde{P}'\in{\cal P}^{\tilde{L}}(\tilde{M})$, on pose

(5)  $\omega_{\tilde{P}'}=\sum_{\tilde{P}}\omega_{\tilde{P}}$, 

\noindent o\`u l'on somme sur les $\tilde{P}\in {\cal P}(\tilde{M})$ tels que $\tilde{P}\cap \tilde{L}=\tilde{P}'$.  Ces fonctions v\'erifient encore les conditions (1), (2) et (4).

Des fonctions $\omega_{\tilde{P}}$ v\'erifiant les conditions ci-dessus sont d\'esormais fix\'ees pour tout l'article.

\bigskip

\subsection{Fonctions rationnelles}
Soit $\tilde{M}$ un espace de Levi de $\tilde{G}$. On note ${\cal A}_{\tilde{M},F}^{\vee}$ le sous-groupe des $\lambda\in {\cal A}_{\tilde{M}}^*$ tels que $<\lambda,X>\in 2\pi {\mathbb Z}$ pour tout $X\in {\cal A}_{\tilde{M},F}$. On pose ${\cal A}^*_{\tilde{M},F}={\cal A}_{\tilde{M}}^*/{\cal A}_{\tilde{M},F}^{\vee}$. C'est un groupe compact que l'on munit de la mesure de masse totale $1$.
Le quotient ${\cal A}_{\tilde{M},{\mathbb C}}^*/i{\cal A}_{\tilde{M},F}^{\vee}$ s'identifie \`a un tore complexe, ce qui permet de parler de fonctions polynomiales ou rationnelles sur ce quotient.  Consid\'erons une telle
  fonction rationnelle $\varphi$. Supposons qu'il existe une famille finie $(\check{\alpha}_{i},c_{i})_{i=1,...,n}$ v\'erifiant les conditions suivantes:
  
  - pour tout $i=1,...,n$, $c_{i}$ est un nombre complexe et $\check{\alpha}_{i}$ est une coracine associ\'ee \`a une racine $\alpha_{i}\in \Sigma(A_{\tilde{M}})$, normalis\'ee de sorte que $\check{\alpha}_{i}\in {\cal A}_{\tilde{M},F}$;
  
  - la fonction $\lambda\mapsto \varphi(\lambda)\prod_{i=1,...,n}(e^{<\lambda,\check{\alpha}_{i}>}-c_{i})$ est polynomiale.
  
  A cette condition, nous dirons que 
    n'a qu'un nombre fini d'hyperplans polaires d'\'equations de la forme $e^{<\lambda,\check{\alpha}>}=c$, pour $\alpha\in \Sigma(A_{\tilde{M}})$. 
    
    {\bf Remarque.} La condition  $\check{\alpha}_{i}\in {\cal A}_{\tilde{M},F}$ est n\'ecessaire pour que la fonction $\lambda\mapsto e^{<\lambda,\check{\alpha}_{i}>}$ soit invariante par $i{\cal A}_{\tilde{M},F}^{\vee}$. Evidemment, cette fonction d\'epend de la normalisation choisie.
    
    \bigskip
   Rappelons que l'on note $\tilde{{\cal A}}_{\tilde{M}}^*$ l'espace des formes lin\'eaires affines sur l'espace r\'eel $\tilde{{\cal A}}_{\tilde{M}}$. On a donc une suite exacte
$$  0\to {\mathbb R}\to \tilde{{\cal A}}_{\tilde{M}}^*\to {\cal A}_{\tilde{M}}^*\to 0.$$
On ajoute un indice ${\mathbb C}$ pour d\'esigner les complexifi\'es.   On a une suite exacte analogue
$$(1) \qquad 0\to {\mathbb C}\to \tilde{{\cal A}}_{\tilde{M},{\mathbb C}}^*\to {\cal A}_{\tilde{M},{\mathbb C}}^*\to 0.$$
On note $\tilde{{\cal A}}_{\tilde{M},F}^{\vee}$ le sous-groupe des $\tilde{\lambda}\in \tilde{{\cal A}}_{\tilde{M}}^*$ tels que $<\tilde{\lambda},X>\in 2\pi {\mathbb Z}$ pour tout $X\in \tilde{{\cal A}}_{\tilde{M},F}$. On a une suite exacte
$$0\to 2\pi{\mathbb Z}\to \tilde{{\cal A}}_{\tilde{M},F}^{\vee}\to {\cal A}_{\tilde{M},F}^{\vee}\to 0.$$
On note usuellement $\tilde{\lambda}$ un \'el\'ement de $\tilde{{\cal A}}_{\tilde{M},{\mathbb C}}^*$ et $\lambda$ sa projection dans ${\cal A}_{\tilde{M},{\mathbb C}}^*$.    Consid\'erons une fonction $\varphi:\tilde{{\cal A}}_{\tilde{M},{\mathbb C}}^*/i\tilde{{\cal A}}_{\tilde{M},F}^{\vee}\to {\mathbb C}$. Pour $X\in \tilde{{\cal A}}_{\tilde{M},F}$, consid\'erons les conditions

(2)  le produit $\varphi(\tilde{\lambda})e^{-<\tilde{\lambda},X>}$ ne d\'epend que de la projection $\lambda$.

\noindent Ce produit ne d\'epend alors que de la projection de $\lambda$ dans  ${\cal A}_{\tilde{M},{\mathbb C}}^*/i{\cal A}_{\tilde{M},F}^{\vee}$. On note $\lambda\mapsto \varphi(\tilde{\lambda})e^{-<\tilde{\lambda},X>}$ la fonction obtenue sur cet ensemble.

 (3) Sous la condition (2),  la fonction $\lambda\mapsto \varphi(\tilde{\lambda})e^{-<\tilde{\lambda},X>}$       est polynomiale; 
 
  (4) sous la condition (2),  la fonction $\lambda\mapsto \varphi(\tilde{\lambda})e^{-<\tilde{\lambda},X>}$       est  rationnelle;  

 (5) sous la condition (2),  la fonction $\lambda\mapsto \varphi(\tilde{\lambda})e^{-<\tilde{\lambda},X>}$       est rationnelle  et n'a qu'un nombre fini d'hyperplans polaires d'\'equations de la forme $e^{<\lambda,\check{\alpha}>}=c$, pour $\alpha\in \Sigma(A_{\tilde{M}})$

Pour un autre point $X'\in \tilde{{\cal A}}_{\tilde{M},F}$, la fonction $\varphi(\tilde{\lambda})e^{-<\tilde{\lambda},X'>}$ est le produit de $\varphi(\tilde{\lambda})e^{-<\tilde{\lambda},X>}$ et du polyn\^ome $e^{<\lambda,X-X'>}$. Il en r\'esulte que les conditions (2), (3), (4) et (5) ne d\'ependent pas du choix de $X$. Plus pr\'ecisement, sous la condition (5), les hyperplans polaires ne d\'ependent pas de ce choix.   Si $\varphi$ v\'erifie (2) et (3), resp. et (4), resp. et (5), nous dirons simplement que $\varphi$ est  polynomiale sur $\tilde{{\cal A}}_{\tilde{M},{\mathbb C}}^*/i\tilde{{\cal A}}_{\tilde{M},F}^{\vee}$, resp.  est rationnelle, resp. est rationnelle et n'a qu'un nombre fini d'hyperplans polaires d'\'equations de la forme $e^{<\lambda,\check{\alpha}>}=c$, pour $\alpha\in \Sigma(A_{\tilde{M}})$). 

Une variante de la propri\'et\'e (2) vaut aussi pour des fonctions $\varphi$ d\'efinies sur le sous-ensemble  $i\tilde{{\cal A}}_{\tilde{M}}^*/i\tilde{{\cal A}}_{\tilde{M},F}^{\vee}$. 

On peut scinder la suite exacte  (1)  en fixant un point $X_{0}\in \tilde{{\cal A}}_{\tilde{M},F}$ et en relevant tout $\lambda\in {\cal A}_{\tilde{M},{\mathbb C}}^*$ en l'unique rel\`evement $\tilde{\lambda}$ tel que $\tilde{\lambda}(X_{0})=0$.  Une telle section envoie ${\cal A}_{\tilde{M},F}^{\vee}$ dans $\tilde{{\cal A}}_{\tilde{M},F}^{\vee}$.  Fixons une telle section. Sous la condition (2),  on a   l'\'egalit\'e $\varphi(\tilde{\lambda})e^{-<\tilde{\lambda},X>}=\varphi(\lambda)e^{-<\lambda,X>}$. 
\bigskip

\subsection{L'application $^c\phi_{\tilde{M}}$}
Fixons un espace de Levi minimal $\tilde{M}_{0}$ et un sous-groupe compact sp\'ecial $K$ de $G(F)$ en bonne position relativement \`a $M_{0}$. Soit $\tilde{M}\in {\cal L}(\tilde{M}_{0})$.

Pour une $\omega$-representation $\tilde{\pi}$ de $\tilde{M}(F)$ et pour $\tilde{\lambda}\in \tilde{{\cal A}}_{\tilde{M},{\mathbb C}}^*/i\tilde{{\cal A}}_{\tilde{M},F}^{\vee}$, on sait d\'efinir la repr\'esentation $\tilde{\pi}_{\tilde{\lambda}}$: on a $\tilde{\pi}_{\tilde{\lambda}}(x)=e^{<\tilde{\lambda},\tilde{H}_{\tilde{M}}(x)>}\tilde{\pi}(x)$ pour tout $x\in \tilde{M}(F)$. Pour ${\bf f}\in C_{c}^{\infty}(\tilde{M}(F))\otimes Mes(M(F))$, la fonction $\tilde{\lambda}\mapsto trace(\tilde{\pi}_{\tilde{\lambda}}({\bf f})) $ est polynomiale sur  $\tilde{{\cal A}}_{\tilde{M},{\mathbb C}}^*/i\tilde{\cal A}_{\tilde{M},F}^{\vee}$.

 On a d\'efini en [W] 2.9 l'espace $D_{temp}(\tilde{M}(F),\omega)$ engendr\'e par les caract\`eres de $\omega$-repr\'esentations temp\'er\'ees de $\tilde{M}(F)$.
 
{\bf Remarque.} On consid\`ere ici les caract\`eres comme des distributions, c'est-\`a-dire des formes lin\'eaires sur $C_{c}^{\infty}(\tilde{M}(F))$. Ils d\'ependent donc des mesures. Vus comme fonctions localement int\'egrables, les caract\`eres vivent dans $D_{temp}(\tilde{M}(F),\omega)\otimes Mes(M(F))^*$. On note $I^{\tilde{M}}(\tilde{\pi},.)$ l'\'el\'ement de $D_{temp}(\tilde{M}(F),\omega)\otimes Mes(M(F))^*$ associ\'e \`a une $\omega$-repr\'esentation $\tilde{\pi}$. 

\bigskip

On a d\'efini  en [W] 2.12 le sous-espace $D_{ell}(\tilde{M}(F),\omega)$ engendr\'e par les caract\`eres de repr\'esentations elliptiques. 
L'induction fournit un isomorphisme
$$D_{temp}(\tilde{M}(F),\omega)\otimes Mes(M(F))^*=$$
$$\oplus_{\tilde{R}\in {\cal L}^{\tilde{M}}(\tilde{M}_{0})/W^M(\tilde{M}_{0})}Ind_{\tilde{R}}^{\tilde{M}}(D_{ell}(\tilde{R}(F),\omega)\otimes Mes(R(F))^*)^{W^M(\tilde{R})}.$$

Soient $\tilde{R}\in {\cal L}^{\tilde{M}}(\tilde{M}_{0})$, $\tilde{\pi}$ une $\omega$-repr\'esentation elliptique de $\tilde{R}(F)$ et ${\bf f}\in C_{c}^{\infty}(\tilde{G}(F))\otimes Mes(G(F))$. Pour $X\in \tilde{{\cal A}}_{\tilde{R},F}$, la fonction
$$\lambda\mapsto J_{\tilde{M}}^{\tilde{G}}(Ind_{\tilde{R}}^{\tilde{M}}(\tilde{\pi}_{\tilde{\lambda}}),{\bf f})e^{-<\tilde{\lambda},X>}$$
d\'efinie en [W] 2.7 est rationnelle sur ${\cal A}_{\tilde{R},{\mathbb C}}^*/i{\cal A}_{\tilde{R},F}^{\vee}$. Elle a un nombre fini d'hyperplans polaires d'\'equations de la forme $e^{<\lambda,\check{\alpha}>}=c$, pour $\alpha\in \Sigma(A_{\tilde{R}})$.

 Ces hyperplans ne d\'ependent pas de $X$ comme on l'a  dit en 1.2. Ils ne d\'ependent pas non plus de ${\bf f}$, au sens qu'il existe une famille finie d'hyperplans de la forme indiqu\'ee de sorte que, pour tout ${\bf f}$, les hyperplans polaires de la fonction ci-dessus appartiennent \`a cette famille. Soit $\tilde{S}\in {\cal P}^{\tilde{G}}(\tilde{R})$. Fixons un point $\nu_{\tilde{S}}\in {\cal A}_{\tilde{R}}^*$ tel que $<\nu_{\tilde{S}},\check{\alpha}>$ soit assez grand pour tout $\alpha\in \Sigma^{\tilde{S}}(\tilde{R})$. Alors l'int\'egrale
$$\int_{\nu_{\tilde{S}}+i{\cal A}_{\tilde{R},F}^*}J_{\tilde{M}}^{\tilde{G}}(Ind_{\tilde{R}}^{\tilde{M}}(\tilde{\pi}_{\tilde{\lambda}}),{\bf f})e^{-<\tilde{\lambda},X>}\,d\lambda$$
ne d\'epend pas du choix de $\nu_{\tilde{S}}$.

Posons
$$(1) \qquad ^cJ_{\tilde{M}}^{\tilde{G}}(Ind_{\tilde{R}}^{\tilde{M}}(\tilde{\pi}),{\bf f})=\sum_{X\in \tilde{{\cal A}}_{\tilde{R},F}}\sum_{\tilde{S}\in {\cal P}^{\tilde{G}}(\tilde{R})}\omega_{\tilde{S}}(X)\int_{\nu_{\tilde{S}}+i{\cal A}_{\tilde{R},F}^*}J_{\tilde{M}}^{\tilde{G}}(Ind_{\tilde{R}}^{\tilde{M}}(\tilde{\pi}_{\tilde{\lambda}}),{\bf f})e^{-<\tilde{\lambda},X>}\,d\lambda,$$
o\`u les points $\nu_{\tilde{S}}$ sont choisis comme ci-dessus. Il n'est pas clair que cette somme converge. Il n'est pas non plus clair qu'elle ne d\'epende que de $Ind_{\tilde{R}}^{\tilde{M}}(\tilde{\pi})$. Ces deux propri\'et\'es sont assur\'ees par la proposition suivante.

\ass{Proposition}{(i) L'expression (1) est une somme finie et ne d\'epend que de $Ind_{\tilde{R}}^{\tilde{M}}(\tilde{\pi})$.

(ii) Il existe une unique application lin\'eaire 
$$\begin{array}{ccc}C_{c}^{\infty}(\tilde{G}(F))\otimes Mes(G(F))&\to& I(\tilde{M}(F),\omega)\otimes Mes(M(F))\\ {\bf f}&\mapsto& {^c\phi}_{\tilde{M}}({\bf f})\\ \end{array}$$
de sorte que, pour tout $\tilde{R}\in {\cal L}^{\tilde{M}}(\tilde{M}_{0})$, toute $\omega$-repr\'esentation elliptique $\tilde{\pi}$ de $\tilde{R}(F)$ et tout ${\bf f}\in C_{c}^{\infty}(\tilde{G}(F))\otimes Mes(G(F))$, on ait l'\'egalit\'e
$$I^{\tilde{M}}(Ind_{\tilde{R}}^{\tilde{M}}(\tilde{\pi}),{^c\phi}_{\tilde{M}}({\bf f}))={^cJ}_{\tilde{M}}^{\tilde{G}}(Ind_{\tilde{R}}^{\tilde{M}}(\tilde{\pi}),{\bf f}).$$}

Preuve. Fixons des mesures de Haar sur tous les groupes pour nous d\'ebarrasser des espaces de mesures. Soit $f\in C_{c}^{\infty}(\tilde{G}(F))$. On va commencer par prouver une formule de descente, sous une forme assez g\'en\'erale car cela nous servira ult\'erieurement.  Fixons un espace de Levi $\tilde{M}'$ tel que $\tilde{R}\subset \tilde{M}'\subset \tilde{M}$. On a
$$Ind_{\tilde{R}}^{\tilde{M}}(\tilde{\pi}_{\tilde{\lambda}})=Ind_{\tilde{M}'}^{\tilde{M}}(Ind_{\tilde{R}}^{\tilde{M}'}(\tilde{\pi}_{\tilde{\lambda}})).$$
Utilisons la formule de descente du lemme 5.4(iv) de [W]. On obtient
$$J_{\tilde{M}}^{\tilde{G}}(Ind_{\tilde{R}}^{\tilde{M}}(\tilde{\pi}_{\tilde{\lambda}}),f)=\sum_{\tilde{L}\in {\cal L}(\tilde{M}')}d_{\tilde{M}'}^{\tilde{G}}(\tilde{M},\tilde{L})J_{\tilde{M}'}^{\tilde{L}}(Ind_{\tilde{R}}^{\tilde{M}'}(\tilde{\pi}_{\tilde{\lambda}}),f_{\tilde{Q},\omega}),$$
o\`u $\tilde{Q}\in {\cal P}(\tilde{L})$ est d\'etermin\'e par le choix d'un param\`etre auxiliaire. Cela transforme (1) en
$$^cJ_{\tilde{M}}^{\tilde{G}}(Ind_{\tilde{R}}^{\tilde{M}}(\tilde{\pi}),f)=\sum_{\tilde{L}\in {\cal L}(\tilde{M}')}d_{\tilde{M}'}^{\tilde{G}}(\tilde{M},\tilde{L}){\cal J}(\tilde{L}),$$
o\`u
$${\cal J}(\tilde{L})=\sum_{X\in \tilde{{\cal A}}_{\tilde{R},F}}\sum_{\tilde{S}\in {\cal P}^{\tilde{G}}(\tilde{R})}\omega_{\tilde{S}}(X)\int_{\nu_{\tilde{S}}+i{\cal A}_{\tilde{R},F}^*}J_{\tilde{M}'}^{\tilde{L}}(Ind_{\tilde{R}}^{\tilde{M}'}(\tilde{\pi}_{\tilde{\lambda}}),f_{\tilde{Q},\omega})e^{-<\tilde{\lambda},X>}\,d\lambda.$$
Fixons $\tilde{L}$. Pour $\tilde{S}'\in {\cal P}^{\tilde{L}}(\tilde{R})$, fixons un point $\nu_{\tilde{S}'}\in {\cal A}_{\tilde{R}}^*$ tel que $<\nu_{\tilde{S}'},\check{\alpha}>$ soit assez grand pour tout $\alpha\in \Sigma^{\tilde{S}'}(A_{\tilde{R}})\subset \Sigma^{\tilde{L}}(A_{\tilde{R}})$. Soit $\tilde{S}\in {\cal P}(\tilde{R})$, posons $\tilde{S}'=\tilde{S}\cap \tilde{L}$. Alors le segment joignant $\nu_{\tilde{S}}$ \`a $\nu_{\tilde{S}'}$ est form\'e de points $\nu$ tels que $<\nu,\check{\alpha}>$ est grand pour toute racine $\alpha\in \Sigma^{\tilde{S}'}(A_{\tilde{R}})$. Il ne coupe aucun hyperplan polaire de la fonction $J_{\tilde{M}'}^{\tilde{L}}(Ind_{\tilde{R}}^{\tilde{M}'}(\tilde{\pi}_{\tilde{\lambda}}),f_{\tilde{Q},\omega})e^{-<\tilde{\lambda},X>}$. On peut donc d\'eplacer le contour d'int\'egration et remplacer l'int\'egrale sur $\nu_{\tilde{S}}+i{\cal A}_{\tilde{R},F}^*$ par celle sur $\nu_{\tilde{S}'}+i{\cal A}_{\tilde{R},F}^*$. La somme sur les $\tilde{S}$ devient une somme sur les $\tilde{S}'$ d'int\'egrales ne d\'ependant que de $\tilde{S}'$ et de la somme
$$\sum_{\tilde{S}; \tilde{S}\cap \tilde{L}=\tilde{S}'}\omega_{\tilde{S}}(X).$$
D'apr\`es 1.1(5), ceci n'est autre que $\omega_{\tilde{S}'}(X)$. On obtient alors
$${\cal J}(\tilde{L})={^cJ}_{\tilde{M}'}^{\tilde{L}}(Ind_{\tilde{R}}^{\tilde{M}'}(\tilde{\pi}),f_{\tilde{Q},\omega}).$$
Nous retenons la formule obtenue
$$(2) \qquad ^cJ_{\tilde{M}}^{\tilde{G}}(Ind_{\tilde{R}}^{\tilde{M}}(\tilde{\pi}),f)=\sum_{\tilde{L}\in {\cal L}(\tilde{M}')}d_{\tilde{M}'}^{\tilde{G}}(\tilde{M},\tilde{L}){^cJ}_{\tilde{M}'}^{\tilde{L}}(Ind_{\tilde{R}}^{\tilde{M}'}(\tilde{\pi}),f_{\tilde{Q},\omega}).$$
Le calcul ci-dessus est formel puisqu'on n'a pas encore prouv\'e que les sommes \'etaient convergentes. Mais il entra\^{\i}ne que, pour d\'emontrer cette convergence, plus pr\'ecis\'ement pour d\'emontrer que les sommes sont finies, il suffit de le faire pour chaque  terme  de la somme ci-dessus. En appliquant ceci au cas $\tilde{M}'=\tilde{R}$, on est  ramen\'e au cas o\`u $\tilde{R}=\tilde{M}$. Dans ce cas, fixons $\tilde{S}\in {\cal P}^{\tilde{G}}(\tilde{M})$. On voit que l'int\'egrale figurant dans (1) n'est non nulle que si la projection $X_{\tilde{G}}$ de $X$ dans $\tilde{{\cal A}}_{\tilde{G},F}$ appartient \`a la projection du support de $f$. Cette projection est finie puisque $f$ est \`a support compact. On peut donc d\'ecomposer la somme en $X$ selon cette projection et prouver que la somme en les $X$ ayant une projection $X_{\tilde{G}}$ fix\'ee est finie. Fixons  comme en 1.2 une section de la projection $\tilde{{\cal A}}_{\tilde{M},{\mathbb C}}^*\to {\cal A}_{\tilde{M},{\mathbb C}}^*$. 
 La fonction 
$$\lambda\mapsto J_{\tilde{M}}^{\tilde{G}}( \tilde{\pi}_{\lambda}, f)$$
est rationnelle. Pour $\vert Re(\lambda)\vert $ assez grand, elle v\'erifie une majoration
$$\vert J_{\tilde{M}}^{\tilde{G}}( \tilde{\pi}_{\lambda}, f)\vert \leq c_{1}e^{C_{1}\vert Re(\lambda)\vert }$$
pour des constantes positives convenables $c_{1}$ et $C_{1}$. 
Faisons tendre le point $\nu_{\tilde{S}}$ vers l'infini de sorte que $\vert \nu_{\tilde{S}}\vert $ reste \'equivalent \`a $<\nu_{\tilde{S}},\check{\alpha}>$ pour tout $\alpha\in \Sigma^{\tilde{S}}(A_{\tilde{M}})$. Puisque les $<\alpha,X>$ sont minor\'es pour de tels $\alpha$ quand $\omega_{\tilde{S}}(X)\not=0$, il existe des nombres r\'eels strictement positifs $c_{1}$, $C_{2}$ et $C_{3}$ de sorte que
$$\vert e^{<\lambda,X>}\vert\leq c_{2} e^{C_{2}\vert  \lambda\vert -C_{3}\vert X\vert \vert \lambda\vert }$$
pour $X$ tel que $\omega_{\tilde{S}}(X)\not=0$ et que la projection $X_{\tilde{G}}$ soit fix\'ee. 
En faisant tendre $\nu_{\tilde{S}}$ vers l'infini, l'int\'egrale devient nulle pourvu que $C_{3}\vert X\vert > C_{1}+C_{2}$. Or il n'y a qu'un nombre fini de $X$ v\'erifiant les conditions pr\'ec\'edentes et ne v\'erifiant pas cette in\'egalit\'e. Cela d\'emontre que la somme en $X$ est finie. 

On doit prouver que la somme (1) ne d\'epend que de $Ind_{\tilde{R}}^{\tilde{M}}(\tilde{\pi})$. C'est-\`a-dire, fixons $x\in M(F)$ tel que $ad_{x}(\tilde{R})\in {\cal L}(\tilde{M}_{0})$. On doit voir que le membre de droite de (1) ne change pas si l'on remplace $\tilde{R}$ par $ad_{x}(\tilde{R})$ et $\tilde{\pi}$ par $\tilde{\pi}\circ ad_{x}^{-1}$. On v\'erifie imm\'ediatement qu'un terme index\'e par $X$ et $\tilde{S}$ de la formule (1) est \'egal au terme  de la nouvelle formule index\'e par $ad_{x}(X)$ et $ad_{x}(\tilde{S})$. 

Pour d\'emontrer (ii), on utilise le th\'eor\`eme de Paley-Wiener ([HL] th\'eor\`eme 3.3, repris en [W] th\'eor\`eme 6.1). On doit d'abord prouver une propri\'et\'e de finitude. A savoir que, pour tout $\tilde{R}$, il existe un ensemble  fini $\Xi$ de  $\omega$-repr\'esentations elliptiques de $\tilde{R}(F)$ tel que, si $ \tilde{\pi}$ v\'erifie $^cJ_{\tilde{M}}^{\tilde{G}}(Ind_{\tilde{R}}^{\tilde{M}}(\tilde{\pi}),f)\not=0$, alors il existe $\tilde{\lambda}\in i\tilde{{\cal A}}_{\tilde{R}}^*$ de sorte    que $\tilde{\pi}_{\tilde{\lambda}}\in \Xi$.   Puisque $f$ est biinvariante par un sous-groupe ouvert compact de $G(F)$, il existe aussi un sous-groupe ouvert compact de $R(F)$ tel que la non-nullit\'e de $^cJ_{\tilde{M}}^{\tilde{G}}(Ind_{\tilde{R}}^{\tilde{M}}(\tilde{\pi}),f)$ entra\^{\i}ne que $\pi$ admet des invariants non nuls par ce compact. Or, \`a torsion pr\`es par  $i\tilde{{\cal A}}_{\tilde{R}}^* $, il n'existe qu'un  nombre fini de $\omega$-repr\'esentations elliptiques de $\tilde{R}(F)$ v\'erifiant cette propri\'et\'e. D'o\`u la finitude requise. 
Fixons $\tilde{R}$ et $\tilde{\pi}$. Fixons aussi un scindage ${\cal A}_{\tilde{R}}^*\to \tilde{{\cal A}}_{\tilde{R}}^*$. On doit montrer que la fonction
$$\mu\mapsto {^cJ}_{\tilde{M}}^{\tilde{G}}(Ind_{\tilde{R}}^{\tilde{M}}(\tilde{\pi}_{\mu}),f)$$
est de Paley-Wiener sur $i{\cal A}_{\tilde{R},F}^*$. Par un changement de variables, on a
 $$^cJ_{\tilde{M}}^{\tilde{G}}(Ind_{\tilde{R}}^{\tilde{M}}(\tilde{\pi}_{\mu}),{\bf f})=\sum_{X\in \tilde{{\cal A}}_{\tilde{R},F}}\sum_{\tilde{S}\in {\cal P}^{\tilde{G}}(\tilde{R})}\omega_{\tilde{S}}(X)e^{<\mu,X>}\int_{\nu_{\tilde{S}}+i{\cal A}_{\tilde{R},F}^*}J_{\tilde{M}}^{\tilde{G}}(Ind_{\tilde{R}}^{\tilde{M}}(\tilde{\pi}_{\tilde{\lambda}}),{\bf f})e^{-<\tilde{\lambda},X>}\,d\lambda.$$
 Le raisonnement fait ci-dessus montre que la somme en $X$ est finie ind\'ependamment de $\mu$. Donc, comme fonction de $\mu$, l'expression ci-dessus est une somme finie de termes $e^{<\mu,X>}$. C'est donc une fonction de Paley-Wiener. Cela ach\`eve la d\'emonstration. $\square$
 
 \bigskip
 
 \subsection{Propri\'et\'es de l'application $^c\phi_{\tilde{M}}$}
 On fixe $\tilde{M}\in {\cal L}(\tilde{M}_{0})$. Conform\'ement \`a nos habitudes, on ajoute un exposant $\tilde{G}$ si besoin est pour indiquer l'espace ambiant: $^c\phi_{\tilde{M}}^{\tilde{G}}$ au lieu de $^c\phi_{\tilde{M}}$.
 
 Soit $\tilde{M}'\in {\cal L}(\tilde{M}_{0})$, supposons $\tilde{M}'\subset \tilde{M}$. Pour ${\bf f}\in C_{c}^{\infty}(\tilde{G}(F))\otimes Mes(G(F))$, on a
 
$$ (1) \qquad (^c\phi^{\tilde{G}}_{\tilde{M}}({\bf f}))_{\tilde{M}',\omega}=\sum_{\tilde{L}\in {\cal L}(\tilde{M}')}d_{\tilde{M}'}^{\tilde{G}}(\tilde{M},\tilde{L}){^c\phi}^{\tilde{L}}_{\tilde{M}'}({\bf f}_{\tilde{Q},\omega}),$$
o\`u $\tilde{Q}\in {\cal P}(\tilde{L})$ est d\'etermin\'e par le choix d'un param\`etre auxiliaire. 

Preuve. Soit $\tilde{R}\in {\cal P}(\tilde{M}_{0})$ avec $\tilde{R}\subset \tilde{M}'$ et soit $\tilde{\pi}$ une $\omega$-repr\'esentation elliptique de $\tilde{R}(F)$. On doit prouver que la distribution $I^{\tilde{M}'}(Ind_{\tilde{R}}^{\tilde{M}'}(\tilde{\pi}),.)$ prend la m\^eme valeur sur les deux membres de (1).  Par d\'efinition de l'induction, on a
$$I^{\tilde{M}'}(Ind_{\tilde{R}}^{\tilde{M}'}(\tilde{\pi}),(^c\phi^{\tilde{G}}_{\tilde{M}}({\bf f}))_{\tilde{M}',\omega})=I^{\tilde{M}}(Ind_{\tilde{R}}^{\tilde{M}}(\tilde{\pi}),{^c\phi}^{\tilde{G}}_{\tilde{M}}({\bf f}))$$
$$= {^cJ}_{\tilde{M}}^{\tilde{G}}(Ind_{\tilde{R}}^{\tilde{M}}(\tilde{\pi}),{\bf f}).$$
C'est le membre de gauche de 1.3(2). On voit aussi que la valeur de $I^{\tilde{M}'}(Ind_{\tilde{R}}^{\tilde{M}'}(\tilde{\pi}),.)$ sur le membre de droite de (1) est le membre de droite de  1.3(2). Cette \'egalit\'e 1.3(2) conclut. $\square$

Pour ${\bf f}\in C_{c}^{\infty}(\tilde{G}(F))\otimes Mes(G(F))$, $y\in G(F)$ et $\tilde{Q}=\tilde{L}U_{Q}\in {\cal F}(\tilde{M}_{0})$, Arthur d\'efinit en [A2] paragraphe 3 une fonction ${\bf f}_{\tilde{Q},y}\in C_{c}^{\infty}(\tilde{L}(F))\otimes Mes(L(F))$ (il convient de glisser dans la d\'efinition notre caract\`ere $\omega$).   Pour $\tilde{Q}=\tilde{G}$, on a simplement ${\bf f}_{\tilde{Q},y}={\bf f}$. On a

$$(2) \qquad \omega(y)^{-1}{^c\phi}^{\tilde{G}}_{\tilde{M}}({\bf f}\circ ad_{y})=\sum_{\tilde{Q}=\tilde{L}U_{Q}\in {\cal F}(\tilde{M})}{^c\phi}_{\tilde{M}}^{\tilde{L}}({\bf f}_{\tilde{Q},y}).$$

Preuve.  Soit $\tilde{R}\in {\cal P}(\tilde{M}_{0})$ avec $\tilde{R}\subset \tilde{M}$ et soit $\tilde{\pi}$ une $\omega$-repr\'esentation elliptique de $\tilde{R}(F)$. On doit prouver que la distribution $I^{\tilde{M}}(Ind_{\tilde{R}}^{\tilde{M}}(\tilde{\pi}),.)$ prend la m\^eme valeur sur les deux membres de (1). Il revient au m\^eme de prouver que
$$(3) \qquad \omega(y)^{-1}{^cJ}_{\tilde{M}}^{\tilde{G}}(Ind_{\tilde{R}}^{\tilde{M}}(\tilde{\pi}),{\bf f}\circ ad_{y})=\sum_{\tilde{Q}=\tilde{L}U_{Q}\in {\cal F}(\tilde{M})}{^cJ}_{\tilde{M}}^{\tilde{L}}(Ind_{\tilde{R}}^{\tilde{M}}(\tilde{\pi}),{\bf f}_{\tilde{Q},y}).$$
Le membre de gauche est d\'efini par 1.3(1) o\`u l'on remplace ${\bf f}$ par ${\bf f}\circ ad_{y}$. On utilise la formule
$$\omega(y)^{-1}J_{\tilde{M}}^{\tilde{G}}(Ind_{\tilde{R}}^{\tilde{M}}(\tilde{\pi}_{\tilde{\lambda}}),{\bf f}\circ ad_{y})=
\sum_{\tilde{Q}=\tilde{L}U_{Q}\in {\cal F}(\tilde{M})}J_{\tilde{M}}^{\tilde{L}}(Ind_{\tilde{R}}^{\tilde{M}}(\tilde{\pi}_{\tilde{\lambda}}),{\bf f}_{\tilde{Q},y}),$$
cf. [A3] lemme 6.2.  On obtient
$$(4) \qquad \omega(y)^{-1}{^cJ}_{\tilde{M}}^{\tilde{G}}(Ind_{\tilde{R}}^{\tilde{M}}(\tilde{\pi}),{\bf f}\circ ad_{y})=\sum_{\tilde{Q}=\tilde{L}U_{Q}\in {\cal F}(\tilde{M})}\sum_{X\in \tilde{{\cal A}}_{\tilde{R},F}}\sum_{\tilde{S}\in {\cal P}^{\tilde{G}}(\tilde{R})}\omega_{\tilde{S}}(X)$$
$$\int_{\nu_{\tilde{S}}+i{\cal A}_{\tilde{R},F}^*}J_{\tilde{M}}^{\tilde{L}}(Ind_{\tilde{R}}^{\tilde{M}}(\tilde{\pi}_{\tilde{\lambda}}),{\bf f}_{\tilde{Q},y})e^{-<\tilde{\lambda},X>}\,d\lambda.$$
Fixons $\tilde{L}$. Pour la m\^eme raison que dans la preuve de 1.3(2), on peut remplacer un point $\nu_{\tilde{S}}$ par $\nu_{\tilde{S}'}$, o\`u $\tilde{S}'=\tilde{S}\cap \tilde{L}$. La somme en $\tilde{S}$ devient une somme en $\tilde{S}'\in {\cal P}^{\tilde{L}}(\tilde{R})$ d'une int\'egrale ne d\'ependant que de $\tilde{S}'$ et de la somme 
$$\sum_{\tilde{S};\tilde{S}\cap \tilde{L}=\tilde{S}' }\omega_{\tilde{S}}(X).$$
D'apr\`es 1.1(5), ceci n'est autre que  $\omega_{\tilde{S}'}(X)$. Alors le membre de droite de (4) devient celui de (3). Cela ach\`eve la preuve. $\square$

Soit $b$ une fonction sur $\tilde{{\cal A}}_{\tilde{G},F}$, \`a valeurs complexes. Pour tout ${\bf f}\in C_{c}^{\infty}(\tilde{G}(F))\otimes Mes(G(F))$, on a

(5) $^c\phi_{\tilde{M}}({\bf f}(b\circ\tilde{H}_{\tilde{G}}))=(^c\phi_{\tilde{M}}({\bf f}))(b\circ\tilde{H}_{\tilde{G}})$.

Preuve. Il revient au m\^eme de dire que, pour $Y\in \tilde{{\cal A}}_{\tilde{G},F}$, si ${\bf f}$ est support\'ee par l'ensemble des $\gamma\in \tilde{G}(F)$ tels que $\tilde{H}_{\tilde{G}}(\gamma)=Y$, alors $^c\phi_{\tilde{M}}({\bf f})$ est support\'ee par l'ensemble des $\gamma\in \tilde{M}(F)$ tels que $\tilde{H}_{\tilde{G}}(\gamma)=Y$. Fixons donc $Y$ et une fonction ${\bf f}$ v\'erifiant la condition de support ci-dessus.  Par transformation de Fourier, la conclusion esp\'er\'ee \'equivaut \`a ce que, pour toute $\omega$-repr\'esentation temp\'er\'ee $\tilde{\sigma}$ de $\tilde{M}(F)$ et tout $\tilde{\mu}\in i\tilde{{\cal A}}_{\tilde{G}}^*$, on a l'\'egalit\'e
$$I^{\tilde{M}}(\tilde{\sigma}_{\tilde{\mu}},{^c\phi}_{\tilde{M}}({\bf f}))=e^{<\tilde{\mu},Y>}I^{\tilde{M}}(\tilde{\sigma},{^c\phi}_{\tilde{M}}({\bf f})).$$
On peut supposer que $\tilde{\sigma}=Ind_{\tilde{R}}^{\tilde{M}}$, o\`u $\tilde{R}\in {\cal L}^{\tilde{M}}(\tilde{M}_{0})$ et $\tilde{\pi}$ est une $\omega$-repr\'esentation elliptique de $\tilde{R}(F)$. On doit alors prouver que
$$^cJ_{\tilde{M}}^{\tilde{G}}(Ind_{\tilde{R}}^{\tilde{M}}(\tilde{\pi}_{\tilde{\mu}}),{\bf f})=e^{<\tilde{\mu},Y>}{^cJ}_{\tilde{M}}^{\tilde{G}}(Ind_{\tilde{R}}^{\tilde{M}}(\tilde{\pi}),{\bf f}).$$
Il suffit pour cela d'utiliser la d\'efinition 1.3(1) et la relation facile
$$J_{\tilde{M}}^{\tilde{G}}(Ind_{\tilde{R}}^{\tilde{M}}(\tilde{\pi}_{\tilde{\lambda}+\tilde{\mu}}),{\bf f})=e^{<\tilde{\mu},Y>}J_{\tilde{M}}^{\tilde{G}}(Ind_{\tilde{R}}^{\tilde{M}}(\tilde{\pi}_{\tilde{\lambda}}),{\bf f}),$$
qui r\'esulte de l'hypoth\`ese sur le support de ${\bf f}$. $\square$

On a d\'efini en [II] 1.6 l'espace $C_{ac}^{\infty}(\tilde{G}(F))$ des fonctions $f:\tilde{G}(F)\to {\mathbb C}$ qui sont biinvariantes par un sous-groupe ouvert compact de $G(F)$ et qui v\'erifient $f(b\circ\tilde{H}_{\tilde{G}})\in C_{c}^{\infty}(\tilde{G}(F))$ pour tout $b\in C_{c}^{\infty}(\tilde{{\cal A}}_{\tilde{G},F})$ (ce dernier espace n'est autre que l'espace des fonctions sur $\tilde{{\cal A}}_{\tilde{G},F}$ \`a support fini). On  a aussi d\'efini l'espace $I_{ac}(\tilde{G}(F),\omega)$. C'est le quotient de  $C_{ac}^{\infty}(\tilde{G}(F))$ par le sous-espace des \'el\'ements dont toutes les int\'egrales orbitales sont nulles. On a

(6) l'application $^c\phi_{\tilde{M}}$ s'\'etend de fa\c{c}on unique en une application lin\'eaire 
$$C_{ac}^{\infty}(\tilde{G}(F))\otimes Mes(G(F))\to I_{ac}(\tilde{M}(F),\omega)\otimes Mes(M(F))$$
qui v\'erifie encore (5).

Preuve. Un \'el\'ement ${\bf f}\in C_{ac}^{\infty}(\tilde{G}(F))\otimes Mes(G(F))$ s'\'ecrit de fa\c{c}on unique ${\bf f}=\sum_{Y\in \tilde{{\cal A}}_{\tilde{G},F}}{\bf f}_{Y}$,  o\`u ${\bf f}_{Y}$ est support\'ee par l'ensemble des $\gamma\in \tilde{G}(F)$ tels que $\tilde{H}_{\tilde{G}}(\gamma)=Y$. La seule d\'efinition  de $^c\phi_{\tilde{M}}({\bf f})$ qui soit compatible avec (5) est
$$^c\phi_{\tilde{M}}({\bf f})=\sum_{Y\in \tilde{{\cal A}}_{\tilde{G},F}} {^c\phi}_{\tilde{M}}({\bf f}_{Y}).$$
Justement, (5) assure la convergence de cette expression. Pour que la somme appartienne \`a $I_{ac}(\tilde{M}(F),\omega)\otimes Mes(M(F))$, il faut prouver qu'il existe un sous-groupe ouvert compact de $M(F)$ tel que, pour tout $Y$, ${^c\phi}_{\tilde{M}}({\bf f}_{Y})$ soit l'image d'un \'el\'ement de $C_{c}^{\infty}(\tilde{M}(F))\otimes Mes(M(F))$ qui soit biinvariant par le sous-groupe en question. On a vu dans la preuve de la proposition 1.3 que,  pour tout $Y$, il y avait un sous-groupe ouvert compact de $M(F)$ tel que ${^c\phi}_{\tilde{M}}({\bf f}_{Y})$ soit l'image d'un \'el\'ement de $C_{c}^{\infty}(\tilde{M}(F))\otimes Mes(M(F))$ qui soit biinvariant par ce sous-groupe. Cette preuve fournit plus: ce sous-groupe ne d\'epend que d'un sous-groupe ouvert compact de $G(F)$ tel que ${\bf f}_{Y}$ soit elle-m\^eme invariante par ce sous-groupe. Par d\'efinition de $C_{ac}^{\infty}(\tilde{G}(F))$, on peut choisir ce dernier sous-groupe ind\'ependant de $Y$. Le sous-groupe de $M(F)$ qu'on en d\'eduit est donc lui-aussi ind\'ependant de $Y$. $\square$

 Pour $f\in C_{c}^{\infty}(\tilde{G}(F))$ et $z\in Z(G;F)^{\theta}$, on note $f^z$ la fonction $\gamma\mapsto f(z\gamma)$. Pour ${\bf f}=f\otimes dg \in C_{c}^{\infty}(\tilde{G}(F))\otimes Mes(G(F))$, on pose ${\bf f}^z=f^z\otimes dg$. Soit $Z$ un sous-groupe de $Z(G)^{\theta}$. 
Notons ${\cal Z}$  l'image dans ${\cal A}_{\tilde{G},F}$  de $Z(F)$ par l'application $H_{\tilde{G}}$. Supposons

(7)  pour tout espace parabolique $\tilde{S}=\tilde{L}U_{S}$ de $\tilde{G}$, tout $X\in \tilde{{\cal A}}_{\tilde{L}}$ et tout $H\in {\cal Z}$, on a l'\'egalit\'e $\omega_{\tilde{S}}(X+H)=\omega_{\tilde{S}}(X)$. 

Montrons que 

(8) sous l'hypoth\`ese (7), on a l'\'egalit\'e $^c\phi^{\tilde{G}}_{\tilde{M}}({\bf f}^z)=(^c\phi^{\tilde{G}}_{\tilde{M}}({\bf f}))^z$ pour tout $z\in Z(F)$ et tout ${\bf f}\in  C_{c}^{\infty}(\tilde{G}(F))\otimes Mes(G(F))$.

Preuve. Soient $\tilde{R}$ un espace de Levi contenu dans $\tilde{M}$ et $\tilde{\pi}$ une $\omega$-repr\'esentation elliptique de $\tilde{R}(F)$. On doit prouver que la distribution $I^{\tilde{M}}(Ind_{\tilde{R}}^{\tilde{M}}(\tilde{\pi}),.)$ prend la m\^eme valeur sur les deux membres de l'\'egalit\'e \`a d\'emontrer. Notons $\mu$ le caract\`ere central de $\pi$ (rappelons que $\pi$ n'est pas irr\'eductible en g\'en\'eral, mais toutes ses composantes irr\'eductibles ont m\^eme caract\`ere central). On a imm\'ediatement
$$I^{\tilde{M}}(Ind_{\tilde{R}}^{\tilde{M}}(\tilde{\pi}),(^c\phi^{\tilde{G}}_{\tilde{M}}({\bf f}))^z)=\mu(z)^{-1}I^{\tilde{M}}(Ind_{\tilde{R}}^{\tilde{M}}(\tilde{\pi}),{^c\phi}^{\tilde{G}}_{\tilde{M}}({\bf f}))=\mu(z)^{-1}{^cJ}_{\tilde{M}}^{\tilde{G}}(Ind_{\tilde{R}}^{\tilde{M}}(\tilde{\pi}),{\bf f})).$$
On a aussi
$$I^{\tilde{M}}(Ind_{\tilde{R}}^{\tilde{M}}(\tilde{\pi}),{^c\phi}^{\tilde{G}}_{\tilde{M}}({\bf f}^z))={^cJ}_{\tilde{M}}^{\tilde{G}}(Ind_{\tilde{R}}^{\tilde{M}}(\tilde{\pi}),{\bf f}^z)$$
$$=\sum_{X\in \tilde{{\cal A}}_{\tilde{R},F}}\sum_{\tilde{S}\in {\cal P}^{\tilde{G}}(\tilde{R})}\omega_{\tilde{S}}(X)\int_{\nu_{\tilde{S}}+i{\cal A}_{\tilde{R},F}^*}J_{\tilde{M}}^{\tilde{G}}(Ind_{\tilde{R}}^{\tilde{M}}(\tilde{\pi}_{\tilde{\lambda}}),{\bf f}^z)e^{-<\tilde{\lambda},X>}\,d\lambda.$$
On a
$$J_{\tilde{M}}^{\tilde{G}}(Ind_{\tilde{R}}^{\tilde{M}}(\tilde{\pi}_{\tilde{\lambda}}),{\bf f}^z)=\mu(z)^{-1}e^{-<\lambda,H >}J_{\tilde{M}}^{\tilde{G}}(Ind_{\tilde{R}}^{\tilde{M}}(\tilde{\pi}_{\tilde{\lambda}}),{\bf f}),$$
o\`u $H=H_{\tilde{G}}(z)$. D'o\`u
$$\int_{\nu_{\tilde{S}}+i{\cal A}_{\tilde{R},F}^*}J_{\tilde{M}}^{\tilde{G}}(Ind_{\tilde{R}}^{\tilde{M}}(\tilde{\pi}_{\tilde{\lambda}}),{\bf f}^z)e^{-<\tilde{\lambda},X>}\,d\lambda=\mu(z)^{-1}\int_{\nu_{\tilde{S}}+i{\cal A}_{\tilde{R},F}^*}J_{\tilde{M}}^{\tilde{G}}(Ind_{\tilde{R}}^{\tilde{M}}(\tilde{\pi}_{\tilde{\lambda}}),{\bf f})e^{-<\tilde{\lambda},X+H>}\,d\lambda.$$
En changeant $X$ en $X-H$, on obtient
$$I^{\tilde{M}}(Ind_{\tilde{R}}^{\tilde{M}}(\tilde{\pi}),{^c\phi}^{\tilde{G}}_{\tilde{M}}({\bf f}^z))=\mu(z)^{-1}\sum_{X\in \tilde{{\cal A}}_{\tilde{R},F}}\sum_{\tilde{S}\in {\cal P}^{\tilde{G}}(\tilde{R})}\omega_{\tilde{S}}(X-H)$$
$$\int_{\nu_{\tilde{S}}+i{\cal A}_{\tilde{R},F}^*}J_{\tilde{M}}^{\tilde{G}}(Ind_{\tilde{R}}^{\tilde{M}}(\tilde{\pi}_{\tilde{\lambda}}),{\bf f})e^{-<\tilde{\lambda},X>}\,d\lambda.$$
L'hypoth\`ese (7) nous permet de faire dispara\^{\i}tre $H$ de cette formule et on obtient
$$I^{\tilde{M}}(Ind_{\tilde{R}}^{\tilde{M}}(\tilde{\pi}),{^c\phi}^{\tilde{G}}_{\tilde{M}}({\bf f}^z))=\mu(z)^{-1}{^cJ}_{\tilde{M}}^{\tilde{G}}(Ind_{\tilde{R}}^{\tilde{M}}(\tilde{\pi}),{\bf f}).$$
Cela prouve l'\'egalit\'e souhait\'ee. $\square$

Soit $\mu$ un caract\`ere de $G(F)$ invariant par $\theta$ (c'est-\`a-dire $\mu\circ ad_{\gamma}=\mu$ pour tout $\gamma\in \tilde{G}(F)$). Soit $\tilde{\mu}$ une fonction non nulle sur $\tilde{G}(F)$ telle que $\tilde{\mu}(x\gamma)=\mu(x)\tilde{\mu}(\gamma)$ pour tous $x\in G(F)$ et $\gamma\in \tilde{G}(F)$.  On a

(9) $^c\phi_{\tilde{M}}({\bf f}\tilde{\mu})=\tilde{\mu}{^c\phi}_{\tilde{M}}({\bf f})$.

Preuve. Soient $\tilde{R}$ et $\tilde{\pi}$ comme dans la preuve pr\'ec\'edente.  Le produit $\tilde{\pi}\tilde{\mu}$ est encore une $\omega$-repr\'esentation elliptique de $\tilde{R}(F)$. On v\'erifie ais\'ement que
$$I^{\tilde{M}}(Ind_{\tilde{R}}^{\tilde{M}}(\tilde{\pi}),{^c\phi}_{\tilde{M}}({\bf f}\tilde{\mu}))=I^{\tilde{M}}(
Ind_{\tilde{R}}^{\tilde{M}}(\tilde{\pi}\tilde{\mu}),{^c\phi}_{\tilde{M}}({\bf f}))=I^{\tilde{M}}(Ind_{\tilde{R}}^{\tilde{M}}(\tilde{\pi}),\tilde{\mu}{^c\phi}_{\tilde{M}}({\bf f})).$$
D'o\`u (9). $\square$

\bigskip

\subsection{D\'efinition de l'application $^c\theta_{\tilde{M}}$}
Soit $\tilde{M}\in {\cal L}(\tilde{M}_{0})$. Nous allons d\'efinir une application lin\'eaire
$$^c\theta^{\tilde{G}}_{\tilde{M}}:C_{c}^{\infty}(\tilde{G}(F))\otimes Mes(G(F))\to I_{ac}(\tilde{M}(F),\omega)\otimes Mes(M(F)).$$
Comme toujours, on a besoin d'en conna\^{\i}tre  une propri\'et\'e par r\'ecurrence. A savoir qu'elle   est $\omega$-\'equivariante, c'est-\`a-dire qu'elle se quotiente en une application lin\'eaire d\'efinie sur $I(\tilde{G}(F),\omega)\otimes Mes(G(F))$. On peut alors poser par r\'ecurrence la d\'efinition 
$$(1) \qquad ^c\theta^{\tilde{G}}_{\tilde{M}}({\bf f})=\phi^{\tilde{G}}_{\tilde{M}}({\bf f})-\sum_{\tilde{L}\in {\cal L}(\tilde{M}), \tilde{L}\not=\tilde{G}}{^c\theta}_{\tilde{M}}^{\tilde{L}}({^c\phi}^{\tilde{G}}_{\tilde{L}}({\bf f}))$$
 pour tout ${\bf f}\in C_{c}^{\infty}(\tilde{G}(F))\otimes Mes(G(F))$. Le terme $\phi^{\tilde{G}}_{\tilde{M}}({\bf f})$ est celui d\'efini en [W] 6.4.

 On montre facilement que l'application $^c\theta_{\tilde{M}}^{\tilde{G}}$ se prolonge en une application d\'efinie sur  $C_{ac}^{\infty}(\tilde{G}(F))\otimes Mes(G(F))$, qui se quotiente en une application d\'efinie sur $I_{ac}(\tilde{G}(F),\omega)\otimes Mes(G(F))$.  Nous laissons ce point au lecteur.
 
 Prouvons la propri\'et\'e d'\'equivariance.
 
 \ass{Proposition}{L'application $^c\theta^{\tilde{G}}_{\tilde{M}}$ se quotiente en une application lin\'eaire d\'efinie sur $I(\tilde{G}(F),\omega)\otimes Mes(G(F))$.}
 
Preuve. On doit prouver que, pour tout ${\bf f}\in C_{c}^{\infty}(\tilde{G}(F))\otimes Mes(G(F))$ et tout $y\in G(F)$, on a l'\'egalit\'e
$$\omega(y)^{-1}{^c\theta}^{\tilde{G}}_{\tilde{M}}({\bf f}\circ ad_{y})={^c\theta}^{\tilde{G}}_{\tilde{M}}({\bf f}).$$
On utilise la relation 1.4(2). On sait que l'application $\phi_{\tilde{M}}^{\tilde{G}}$ v\'erifie la m\^eme relation. D'apr\`es la d\'efinition (1), on obtient 
$$\omega(y)^{-1}{^c\theta}^{\tilde{G}}_{\tilde{M}}({\bf f}\circ ad_{y})=\sum_{\tilde{Q}'=\tilde{L}'U_{Q'}\in {\cal F}(\tilde{M})}\phi_{\tilde{M}}^{\tilde{L}'}({\bf f}_{\tilde{Q}',y})$$
$$-\sum_{\tilde{L}\in {\cal L}(\tilde{M}), \tilde{L}\not=\tilde{G}}\sum_{\tilde{Q}'=\tilde{L}'U_{Q'}\in {\cal F}(\tilde{L})}{^c\theta}_{\tilde{M}}^{\tilde{L}}({^c\phi}_{\tilde{L}}^{\tilde{L}'}({\bf f}_{\tilde{Q}',y}))$$
$$=\sum_{\tilde{Q}'=\tilde{L}'U_{Q'}\in {\cal F}(\tilde{M}),\tilde{Q}'\not=\tilde{G}}\left(\phi_{\tilde{M}}^{\tilde{L}'}({\bf f}_{\tilde{Q}',y})-\sum_{\tilde{L}\in {\cal L}^{\tilde{L}'}(\tilde{M})}{^c\theta}_{\tilde{M}}^{\tilde{L}}({^c\phi}_{\tilde{L}}^{\tilde{L}'}({\bf f}_{\tilde{Q}',y}))\right)$$
$$+\phi_{\tilde{M}}^{\tilde{G}}({\bf f}_{\tilde{G},y})-\sum_{\tilde{L}\in {\cal L}^{\tilde{G}}(\tilde{M}),\tilde{L}\not=\tilde{G}}{^c\theta}_{\tilde{M}}^{\tilde{L}}({^c\phi}_{\tilde{L}}^{\tilde{G}}({\bf f}_{\tilde{G},y})).$$
La premi\`ere somme est nulle d'apr\`es la d\'efinition (1) appliqu\'ee avec $\tilde{G}$ remplac\'e par $\tilde{L}'$. Comme on l'a dit, on a ${\bf f}_{\tilde{G},y}={\bf f}$ et la deuxi\`eme somme est \'egale \`a $^c\theta_{\tilde{M}}^{\tilde{G}}({\bf f})$. Cela conclut. $\square$

\bigskip

\subsection{Propri\'et\'es de l'application $^c\theta^{\tilde{G}}_{\tilde{M}}$}
 On a fix\'e un espace de Levi minimal $\tilde{M}_{0}$ et un sous-groupe compact sp\'ecial $K$ de $G(F)$ en bonne position relativement \`a $M_{0}$. On a alors d\'efini l'application $^c\theta^{\tilde{G}}_{\tilde{M}}$ pour un espace de Levi $\tilde{M}$ contenant $\tilde{M}_{0}$. L'espace $\tilde{M}_{0}$ \'etant fix\'e, on montre qu'elle ne d\'epend pas de $K$. L'argument est que, quand on remplace $K$ par un autre groupe compact sp\'ecial $K'$ en bonne position relativement \`a $M_{0}$, les ingr\'edients basiques $J_{\tilde{M}}^{\tilde{G}}(\tilde{\sigma},{\bf f})$ de toutes nos distributions sont chang\'es en d'autres qui se d\'eduisent des premiers par une formule de la m\^eme forme que 1.4(2). Un raisonnement similaire \`a celui de la preuve de la proposition pr\'ec\'edente permet alors de conclure. Nous renvoyons pour plus de d\'etails \`a [A2] proposition 13.2.
Cela \'etant, ajoutons pour un instant l'espace $\tilde{M}_{0}$ dans la notation: $^c\theta^{\tilde{G}}_{\tilde{M},\tilde{M}_{0}}$ au lieu de $^c\theta^{\tilde{G}}_{\tilde{M}}$. Par simple transport de structure, on a une \'egalit\'e
$$(1) \qquad ^c\theta^{\tilde{G}}_{\tilde{M},\tilde{M}_{0}}({\bf f}\circ ad_{y})=({^c\theta}^{\tilde{G}}_{ad_{y}(\tilde{M}),ad_{y}(\tilde{M}_{0})}({\bf f}))\circ ad_{y}$$
pour tout ${\bf f}\in C_{c}^{\infty}(\tilde{G}(F))\otimes Mes(G(F))$ et tout $y\in G(F)$. Pour $y\in M(F)$, gr\^ace \`a la proposition pr\'ec\'edente, cette \'egalit\'e devient
$$\omega(y){^c\theta}^{\tilde{G}}_{\tilde{M},\tilde{M}_{0}}({\bf f})=\omega(y){^c\theta}^{\tilde{G}}_{\tilde{M},ad_{y}(\tilde{M}_{0})}({\bf f}).$$
Il en r\'esulte que notre application $^c\theta^{\tilde{G}}_{\tilde{M},\tilde{M}_{0}}$ ne change pas quand on remplace $\tilde{M}_{0}$ par $ad_{y}(\tilde{M}_{0})$. Puisque tous les espaces de Levi minimaux contenus dans $\tilde{M}$ sont conjugu\'es par un \'el\'ement de $M(F)$, cette application ne d\'epend  pas de $\tilde{M}_{0}$. Cet espace peut de nouveau dispara\^{\i}tre de la notation. D'autre part, l'application $^c\theta^{\tilde{G}}_{\tilde{M}}$ peut \^etre d\'efinie pour tout $\tilde{M}$ puisque, un tel espace de Levi \'etant fix\'e, on peut toujours choisir un $\tilde{M}_{0}$ minimal contenu dans $\tilde{M}$. 

On fixe un espace de Levi $\tilde{M}$ quelconque. Le groupe $W(\tilde{M})$ agit naturellement sur $I_{ac}(\tilde{M}(F),\omega) $. Rappelons la d\'efinition de cette action. Notons $Norm_{G(F)}(\tilde{M})$ le normalisateur de $\tilde{M}$ dans $G(F)$. On d\'efinit une action de ce groupe sur l'espace des fonctions sur $\tilde{M}(F)$: \`a $n\in Norm_{G(F)}(\tilde{M})$ et \`a une fonction   $\varphi$ sur $\tilde{M}(F)$, on associe la fonction $n\varphi$ sur $\tilde{M}(F)$ d\'efinie par $(n\varphi)(\gamma)=\omega(n)\varphi(n^{-1}\gamma n)$.  De cette action se d\'eduit une action de $Norm_{G(F)}(\tilde{M})$ sur $I(\tilde{M}(F),\omega) $ et $I_{ac}(\tilde{M}(F),\omega) $. Sur ces espaces, l'action du sous-groupe $M(F)\subset Norm_{G(F)}(\tilde{M})$ est triviale et l'action se quotiente en une action du quotient  $W(\tilde{M})$. Par tensorisation avec l'action triviale de ce groupe sur $Mes(M(F))$, on obtient une action de $W(\tilde{M})$ sur  $I_{ac}(\tilde{M}(F),\omega)\otimes Mes(M(F))$.
On peut 
appliquer (1) \`a un \'el\'ement $y$ qui normalise $\tilde{M}$. On obtient le r\'esultat suivant 

(2) $^c\theta^{\tilde{G}}_{\tilde{M}}$ prend ses valeurs dans le sous-espace des \'el\'ements de $I_{ac}(\tilde{M}(F),\omega)\otimes Mes(M(F))$ qui sont invariants par $W(\tilde{M})$.

Soit $b$ une fonction sur $\tilde{{\cal A}}_{\tilde{G},F}$, \`a valeurs complexes. Pour tout ${\bf f}\in C_{c}^{\infty}(\tilde{G}(F))\otimes Mes(G(F))$, on a l'\'egalit\'e
$$(3)\qquad ^c\theta^{\tilde{G}}_{\tilde{M}}({\bf f}(b\circ \tilde{H}_{\tilde{G}}))=(^c\theta^{\tilde{G}}_{\tilde{M}}({\bf f}))(b\circ \tilde{H}_{\tilde{G}}).$$
Cela r\'esulte par r\'ecurrence de 1.4(5) et de la relation similaire v\'erifi\'ee par l'application $\phi^{\tilde{G}}_{\tilde{M}}$.

Soit $Z$ un sous-groupe de $Z(G)^{\theta}$. On a

(4)  sous l'hypoth\`ese  1.4(7), on a l'\'egalit\'e $^c\theta^{\tilde{G}}_{\tilde{M}}({\bf f}^z)=(^c\theta^{\tilde{G}}_{\tilde{M}}({\bf f}))^z$ pour tout $z\in Z(F)$ et tout ${\bf f}\in  C_{c}^{\infty}(\tilde{G}(F))\otimes Mes(G)$.

Cela r\'esulte par r\'ecurrence de 1.4(8) et de la relation similaire v\'erifi\'ee par l'application $\phi^{\tilde{G}}_{\tilde{M}}$ ( dans le cas de $\phi^{\tilde{G}}_{\tilde{M}}$, cette propri\'et\'e est ind\'ependante   de l'hypoth\`ese 1.4(7)).

Soit $\tilde{\mu}$ un "caract\`ere affine" comme en 1.4(9). On a de m\^eme

(5) $^c\theta^{\tilde{G}}_{\tilde{M}}({\bf f}\tilde{\mu})=\tilde{\mu}{^c\theta}^{\tilde{G}}_{\tilde{M}}({\bf f})$ pour  tout ${\bf f}\in  C_{c}^{\infty}(\tilde{G}(F))\otimes Mes(G)$.

\ass{Lemme}{Soit $\tilde{M}'$ un espace de Levi contenu dans $\tilde{M}$. Pour tout ${\bf f }\in I(\tilde{G}(F),\omega)\otimes Mes(G(F))$, on a l'\'egalit\'e
$$(^c\theta^{\tilde{G}}_{\tilde{M}}({\bf f}))_{\tilde{M}',\omega}=\sum_{\tilde{L}\in {\cal L}(\tilde{M}')}d_{\tilde{M}'}(\tilde{M},\tilde{L}){^c\theta}_{\tilde{M}'}^{\tilde{L}}({\bf f}_{\tilde{L},\omega}).$$}

Preuve.  On applique la formule de d\'efinition 1.5(1)  et on lui applique l'application "terme constant" relative \`a $\tilde{M}'$. On  obtient  l'\'egalit\'e
$$ (\phi_{\tilde{M}}^{\tilde{G}}({\bf f}))_{\tilde{M}',\omega}=\sum_{\tilde{L}_{1}\in {\cal L}(\tilde{M})}(^c\theta_{\tilde{M}}^{\tilde{L}_{1}}(^c\phi_{\tilde{L}_{1}}({\bf f})))_{\tilde{M}',\omega}.$$
Pour un  terme de la somme index\'e par $\tilde{L}_{1}\not=\tilde{G}$, on peut 
utiliser par r\'ecurrence la formule du pr\'esent \'enonc\'e. Pour le terme index\'e par $\tilde{L}_{1}=\tilde{G}$, on ne peut pas. Mais on peut quand m\^eme, \`a condition d'ajouter la diff\'erence entre les deux membres de cet \'enonc\'e. Pr\'ecis\'ement, posons
$$X=(^c\theta^{\tilde{G}}_{\tilde{M}}({\bf f}))_{\tilde{M}',\omega}-\sum_{\tilde{L}\in {\cal L}(\tilde{M}')}d_{\tilde{M}'}(\tilde{M},\tilde{L}){^c\theta}_{\tilde{M}'}^{\tilde{L}}({\bf f}_{\tilde{L},\omega}).$$
On obtient alors
$$ (\phi_{\tilde{M}}^{\tilde{G}}({\bf f}))_{\tilde{M}',\omega}=X+\sum_{\tilde{L}_{1}\in {\cal L}(\tilde{M})}\sum_{\tilde{L}_{2}\in {\cal L}(\tilde{M}'), \tilde{L}_{2}\subset \tilde{L}_{1}}d_{\tilde{M}'}^{\tilde{L}_{1}}(\tilde{M},\tilde{L}_{2}){^c\theta}_{\tilde{M}'}^{\tilde{L}_{2}}((^c\phi_{\tilde{L}_{1}}^{\tilde{G}}({\bf f}))_{\tilde{L}_{2},\omega}).$$
Utilisons la relation 1.4(1) et la relation analogue concernant l'application $\phi_{\tilde{M}}^{\tilde{G}}$. On obtient
$$(6)\qquad  \sum_{\tilde{L}\in {\cal L}(\tilde{M}')}d_{\tilde{M}'}(\tilde{M},\tilde{L})\phi_{\tilde{M}'}^{\tilde{L}}({\bf f}_{\tilde{Q},\omega})=X+Y,$$
o\`u
 $$  Y=\sum_{\tilde{L}_{1}\in {\cal L}(\tilde{M})}\sum_{\tilde{L}_{2}\in {\cal L}(\tilde{M}'), \tilde{L}_{2}\subset \tilde{L}_{1}}d_{\tilde{M}'}^{\tilde{L}_{1}}(\tilde{M},\tilde{L}_{2})\sum_{\tilde{L}_{3}\in {\cal L}(\tilde{L}_{2})}d_{\tilde{L}_{2}}^{\tilde{G}}(\tilde{L}_{1},\tilde{L}_{3}){^c\theta}_{\tilde{M}'}^{\tilde{L}_{2}}(^c\phi_{\tilde{L}_{2}}^{\tilde{L}_{3}}({\bf f}_{\tilde{Q}_{3},\omega})).$$
 
 Revenons sur la d\'efinition des espaces $\tilde{Q}$ et $\tilde{Q}_{3}$ qui interviennent ici.  On fixe un point $\xi\in {\cal A}_{\tilde{M}'}^{\tilde{G}}$ en position g\'en\'erale. Soit $\tilde{L}\in {\cal L}(\tilde{M}')$ tel que $d_{\tilde{M}'}(\tilde{M},\tilde{L})\not=0$. On a alors la d\'ecomposition
 $${\cal A}_{\tilde{M}'}^{\tilde{G}}={\cal A}^{\tilde{G}}_{\tilde{M}}\oplus {\cal A}_{\tilde{L}}^{\tilde{G}}.$$
  Conform\'ement \`a celle-ci, on projette $\xi$ en un point $\xi[\tilde{L}]\in {\cal A}_{\tilde{L}}^{\tilde{G}}$. L'espace parabolique $\tilde{Q}$ est l'unique \'el\'ement de ${\cal P}(\tilde{L})$ tel que $\xi[\tilde{L}]$ appartienne \`a la chambre positive associ\'ee \`a $\tilde{Q}$.  De m\^eme, fixons $\tilde{L}_{1}\in {\cal L}(\tilde{M})$, $\tilde{L}_{2}\in {\cal L}(\tilde{M}')$ avec $\tilde{L}_{2}\subset \tilde{L}_{1}$ et $d_{\tilde{M}'}^{\tilde{L}_{1}}(\tilde{M},\tilde{L}_{2})\not=0$. Fixons un point $\xi(\tilde{L}_{1},\tilde{L}_{2})\in {\cal A}_{\tilde{L}_{2}}^{\tilde{G}}$ en position g\'en\'erale.  Soit $\tilde{L}_{3}\in {\cal L}(\tilde{L}_{2})$ tel que $d_{\tilde{L}_{2}}^{\tilde{G}}(\tilde{L}_{1},\tilde{L}_{3})\not=0$. On a alors la d\'ecomposition
  $${\cal A}_{\tilde{L}_{2}}^{\tilde{G}}={\cal A}_{\tilde{L}_{1}}^{\tilde{G}}\oplus {\cal A}_{\tilde{L}_{3}}^{\tilde{G}}.$$
  Conform\'ement \`a celle-ci, on projette $\xi(\tilde{L}_{1},\tilde{L}_{2})$ en un point $\xi(\tilde{L}_{1},\tilde{L}_{2})[\tilde{L}_{3}]\in {\cal A}_{\tilde{L}_{3}}^{\tilde{G}}$.  L'espace parabolique $\tilde{Q}_{3}$ est l'unique \'el\'ement de ${\cal P}(\tilde{L}_{3})$ tel que $\xi(\tilde{L}_{1},\tilde{L}_{2})[\tilde{L}_{3}]$ appartienne \`a la chambre positive associ\'ee \`a $\tilde{Q}$.  Le point $\xi$ \'etant fix\'e, nous allons choisir le point $\xi(\tilde{L}_{1},\tilde{L}_{2})$ de la fa\c{c}on suivante. On note comme toujours $\xi_{\tilde{L}_{1}}$ et $\xi^{\tilde{L}_{1}}$ les projections orthogonales sur ${\cal A}^{\tilde{G}}_{\tilde{L}_{1}}$ et ${\cal A}_{\tilde{M}'}^{\tilde{L}_{1}}$. On a suppos\'e  
$d_{\tilde{M}'}^{\tilde{L}_{1}}(\tilde{M},\tilde{L}_{2})\not=0$. On a donc la d\'ecomposition
$${\cal A}_{\tilde{M}'}^{\tilde{L}_{1}}={\cal A}_{\tilde{M}}^{\tilde{L}_{1}}\oplus {\cal A}_{\tilde{L}_{2}}^{\tilde{L}_{1}}.$$
On peut d\'ecomposer conform\'ement $\xi^{\tilde{L}_{1}}$ en $\xi^{\tilde{L}_{1}}[\tilde{M}]+\xi^{\tilde{L}_{1}}[\tilde{L}_{2}]$. On choisit $\xi(\tilde{L}_{1},\tilde{L}_{2})=\xi_{\tilde{L}_{1}}+\xi^{\tilde{L}_{1}}[\tilde{L}_{2}]=\xi-\xi^{\tilde{L}_{1}}[\tilde{M}]$. 

On r\'ecrit la d\'efinition
 $$Y=\sum_{\tilde{L}_{3}\in {\cal L}(\tilde{M}')}\sum_{\tilde{L}_{2}\in {\cal L}^{\tilde{L}_{3}}(\tilde{M}')}Y(\tilde{L}_{2},\tilde{L}_{3}),$$
 o\`u
 $$(7) \qquad  Y(\tilde{L}_{2},\tilde{L}_{3})=\sum_{\tilde{L}_{1}\in {\cal L}(\tilde{M}), \tilde{L}_{2}\subset \tilde{L}_{1}}d_{\tilde{M}'}^{\tilde{L}_{1}}(\tilde{M},\tilde{L}_{2})d_{\tilde{L}_{2}}^{\tilde{G}}(\tilde{L}_{1},\tilde{L}_{3}){^c\theta}_{\tilde{M}'}^{\tilde{L}_{2}}(^c\phi_{\tilde{L}_{2}}^{\tilde{L}_{3}}({\bf f}_{\tilde{Q}_{3},\omega})).$$
Fixons $\tilde{L}_{2}$ et $\tilde{L}_{3}$. On a vu en [II] 1.7(5) que la somme ci-dessus en $\tilde{L}_{1}$ \'etait vide si $d_{\tilde{M}'}^{\tilde{G}}(\tilde{M},\tilde{L}_{3})=0$. Si $d_{\tilde{M}'}^{\tilde{G}}(\tilde{M},\tilde{L}_{3})\not=0$, la somme est r\'eduite \`a un seul \'el\'ement pour lequel
$$d_{\tilde{M}'}^{\tilde{L}_{1}}(\tilde{M},\tilde{L}_{2})d_{\tilde{L}_{2}}^{\tilde{G}}(\tilde{L}_{1},\tilde{L}_{3})=d_{\tilde{M}'}^{\tilde{G}}(\tilde{M},\tilde{L}_{3}).$$
L'unique \'el\'ement $\tilde{L}_{1}$ de la somme est caract\'eris\'e par l'\'egalit\'e
$$(8) \qquad {\cal A}_{\tilde{L}_{1}}^{\tilde{G}}={\cal A}_{\tilde{M}}^{\tilde{G}}\cap {\cal A}_{\tilde{L}_{2}}^{\tilde{G}}.$$
Supposons $d_{\tilde{M}'}^{\tilde{G}}(\tilde{M},\tilde{L}_{3})\not=0$.  Dans (7) appara\^{\i}t un espace parabolique $\tilde{Q}_{3}$. En posant $\tilde{L}=\tilde{L}_{3}$, cet espace de Levi intervient dans le membre de gauche de la formule (6) et il lui est associ\'e un espace parabolique $\tilde{Q}$. Montrons que

(9) $\tilde{Q}_{3}=\tilde{Q}$.

D'apr\`es la construction rappel\'ee ci-dessus, il suffit de prouver que $\xi[\tilde{L}]=\xi(\tilde{L}_{1},\tilde{L}_{2})[\tilde{L}_{3}]$. On a 
$$\xi[\tilde{L}]\in \xi+{\cal A}_{\tilde{M}}^{\tilde{G}}=\xi(\tilde{L}_{1},\tilde{L}_{2})+\xi^{\tilde{L}_{1}}[\tilde{M}]+{\cal A}_{\tilde{M}}^{\tilde{G}}=\xi(\tilde{L}_{1},\tilde{L}_{2})+{\cal A}_{\tilde{M}}^{\tilde{G}}$$
puisque $\xi^{\tilde{L}_{1}}[\tilde{M}]\in {\cal A}_{\tilde{M}}^{\tilde{G}}$. Les \'el\'ements $\xi[\tilde{L}]$ et $\xi(\tilde{L}_{1},\tilde{L}_{2})$ appartenant tous deux \`a ${\cal A}_{\tilde{L}_{2}}^{\tilde{G}}$, la relation pr\'ec\'edente se renforce en
$$\xi[\tilde{L}]\in \xi(\tilde{L}_{1},\tilde{L}_{2})+{\cal A}_{\tilde{M}}^{\tilde{G}}\cap {\cal A}_{\tilde{L}_{2}}^{\tilde{G}},$$
ou encore, d'apr\`es (8), en
$$\xi[\tilde{L}]\in \xi(\tilde{L}_{1},\tilde{L}_{2})+{\cal A}_{\tilde{L}_{1}}^{\tilde{G}}.$$
 Alors $\xi[\tilde{L}]$  appartient \`a l'intersection de cet espace affine avec ${\cal A}_{\tilde{L}_{3}}^{\tilde{G}}$. Or cette intersection est r\'eduite au point $\xi(\tilde{L}_{1},\tilde{L}_{2})[\tilde{L}_{3}]$. Cela prouve (9).
 
 Modifions les notations en rempla\c{c}ant $\tilde{L}_{3}$ par $\tilde{L}$ et $\tilde{L}_{2}$ par $\tilde{L}'$. On obtient 
 $$Y(\tilde{L}',\tilde{L})=d_{\tilde{M}'}^{\tilde{G}}(\tilde{M},\tilde{L}){^c\theta}_{\tilde{M}'}^{\tilde{L}'}(^c\phi_{\tilde{L}'}^{\tilde{L}}({\bf f}_{\tilde{Q},\omega}))$$
 puis
 $$Y=\sum_{\tilde{L}\in {\cal L}(\tilde{M})}d_{\tilde{M}'}^{\tilde{G}}(\tilde{M},\tilde{L})\sum_{\tilde{L}'\in {\cal L}^{\tilde{L}}(\tilde{M}')}{^c\theta}_{\tilde{M}'}^{\tilde{L}'}(^c\phi_{\tilde{L}'}^{\tilde{L}}({\bf f}_{\tilde{Q},\omega})).$$
 Reportons cette relation dans la formule (6). Les termes index\'es par $\tilde{L}$ de chaque c\^ot\'e de la formule sont \'egaux par d\'efinition des applications $^c\theta_{\tilde{M}'}^{\tilde{L}}$. On obtient simplement $X=0$. C'est ce qu'affirme l'\'enonc\'e. $\square$

\bigskip

\subsection{Fonctions de Schwartz}
 Soit $\tilde{\pi}$ une $\omega$-repr\'esentation temp\'er\'ee de $\tilde{G}(F)$.  
Pour une fonction $\boldsymbol{\varphi}\in I(\tilde{G}(F),\omega)\otimes Mes(G(F))$, la fonction $\tilde{\lambda}\mapsto I^{\tilde{G}}(\tilde{\pi}_{\tilde{\lambda}},\boldsymbol{\varphi})$ est bien d\'efinie pour tout $\tilde{\lambda}\in \tilde{{\cal A}}_{\tilde{G},{\mathbb C}}^*/i\tilde{{\cal A}}_{\tilde{G},F}^{\vee}$. C'est une fonction polynomiale sur cet espace. Pour $X\in \tilde{{\cal A}}_{\tilde{G},F}$, on d\'efinit le coefficient de Fourier
$$I^{\tilde{G}}(\tilde{\pi},X,\boldsymbol{\varphi})=\int_{i{\cal A}^*_{\tilde{G},F}}I^{\tilde{G}}(\tilde{\pi}_{\tilde{\lambda}},\boldsymbol{\varphi})e^{-<\tilde{\lambda},X>}\,d\lambda.$$
Comme en 1.2, le terme $\tilde{\lambda}$ apparaissant dans l'int\'egrale est un rel\`evement quelconque de $\lambda$ dans $i\tilde{{\cal A}}_{\tilde{G}}^*$. 
On a les \'egalit\'es
$$(1) \qquad I^{\tilde{G}}(\tilde{\pi},X,\boldsymbol{\varphi})=I^{\tilde{G}}(\tilde{\pi},({\bf 1}_{X}\circ \tilde{H}_{\tilde{G}})\boldsymbol{\varphi}),$$
o\`u ${\bf 1}_{X}$ est la fonction caract\'eristique de $X$ dans $\tilde{{\cal A}}_{\tilde{G},F}$ et
$$(2) \qquad I^{\tilde{G}}(\tilde{\pi}_{\tilde{\lambda}},\boldsymbol{\varphi})=\sum_{X\in \tilde{{\cal A}}_{\tilde{G},F}}e^{<\tilde{\lambda},X>} I^{\tilde{G}}(\tilde{\pi},X,\boldsymbol{\varphi})$$
pour tout $\tilde{\lambda}\in \tilde{{\cal A}}_{\tilde{G},{\mathbb C}}^*/i\tilde{{\cal A}}_{\tilde{G},F}^{\vee}$.
Cette somme ne contient qu'un nombre fini de  termes non nuls. 

Pour une fonction $\boldsymbol{\varphi}\in I_{ac}(\tilde{G}(F),\omega)\otimes Mes(G(F))$,  on ne peut plus en g\'en\'eral d\'efinir la fonction $\tilde{\lambda}\mapsto I^{\tilde{G}}(\tilde{\pi}_{\tilde{\lambda}},\boldsymbol{\varphi})$.  Par contre, pour  $X\in \tilde{{\cal A}}_{\tilde{G},F}$, on peut d\'efinir le coefficient de Fourier $ I^{\tilde{G}}(\tilde{\pi},X,\boldsymbol{\varphi})$, cf. [W] 6.4. Dans notre cas o\`u le corps de base est non-archim\'edien, il est d\'efini par la formule (1): le membre de droite de cette formule a un sens puisque la fonction $ ({\bf 1}_{X}\circ \tilde{H}_{\tilde{G}})\boldsymbol{\varphi}$ est \`a support compact. Nous dirons que $\boldsymbol{\varphi}$ est de Schwartz si, pour toute $\omega$-repr\'esentation temp\'er\'ee $\tilde{\pi}$ de $\tilde{G}(F)$, la fonction
$$(3) \qquad X\mapsto  I^{\tilde{G}}(\tilde{\pi},X,\boldsymbol{\varphi})$$
est \`a d\'ecroissance rapide sur $\tilde{{\cal A}}_{\tilde{G},F}$. Dans ce cas, on peut d\'efinir une fonction $\tilde{\lambda}\mapsto I^{\tilde{G}}(\tilde{\pi},\tilde{\lambda},\boldsymbol{\varphi})$ sur $i\tilde{{\cal A}}_{\tilde{G}}^*/i\tilde{{\cal A}}_{\tilde{G},F}^{\vee}$ par la formule d'inversion de Fourier
$$(4) \qquad  I^{\tilde{G}}(\tilde{\pi},\tilde{\lambda},\boldsymbol{\varphi})=\sum_{X\in \tilde{{\cal A}}_{\tilde{G},F}}e^{<\tilde{\lambda},X>} I^{\tilde{G}}(\tilde{\pi},X,\boldsymbol{\varphi}).$$
C'est une fonction $C^{\infty}$ de $\tilde{\lambda}$. Inversement, supposons que, pour tout $\tilde{\pi}$, il existe une fonction $C^{\infty}$ sur $i\tilde{{\cal A}}_{\tilde{G}}^*/i\tilde{{\cal A}}_{\tilde{G},F}^{\vee}$ v\'erifiant la condition 1.2(2) et  dont la fonction (3) soit la transform\'ee de Fourier. Alors cette fonction (3) est \`a d\'ecroissance rapide, donc $\boldsymbol{\varphi}$ est de Schwartz.

Pour une fonction de Schwartz $\boldsymbol{\varphi}$, il arrive que la fonction (4) d\'efinie sur $i\tilde{{\cal A}}_{\tilde{G}}^*/i\tilde{{\cal A}}_{\tilde{G},F}^{\vee}$ se prolonge en une fonction rationnelle sur  $\tilde{{\cal A}}_{\tilde{G},{\mathbb C}}^*/i\tilde{{\cal A}}_{\tilde{G},F}^{\vee}$. Dans ce cas, on notera encore $\tilde{\lambda}\mapsto I^{\tilde{G}}(\tilde{\pi},\tilde{\lambda},\boldsymbol{\varphi})$ ce prolongement. 
Soulignons que cela n'implique nullement que  le membre de droite de (4) soit convergent pour $\tilde{\lambda}\not\in i\tilde{{\cal A}}_{\tilde{G}}^*/i\tilde{{\cal A}}_{\tilde{G},F}^{\vee}$. M\^eme s'il l'est, la fonction prolong\'ee n'a pas de raison d'\^etre \'egale \`a ce membre de droite.

\bigskip

\subsection{Une propri\'et\'e d'annulation}

On fixe un espace de Levi $\tilde{M}$. Soit   ${\bf f}\in I(\tilde{G}(F),\omega)\otimes Mes(G(F))$.     On a

(1) la fonction $^c\theta_{\tilde{M}}^{\tilde{G}}({\bf f})$ est de Schwartz;

(2) soit  $\tilde{\pi}$ une $\omega$-repr\'esentation temp\'er\'ee de $\tilde{M}(F)$; la fonction $\tilde{\lambda}\mapsto I^{\tilde{M}}(\tilde{\pi},\tilde{\lambda},{^c\theta}^{\tilde{G}}_{\tilde{M}}({\bf f})) $ sur $i\tilde{{\cal A}}_{\tilde{M}}^*/i\tilde{{\cal A}}_{\tilde{M},F}^{\vee}$ est la restriction d'une fonction rationnelle sur $\tilde{{\cal A}}_{\tilde{M},{\mathbb C}}^*/i\tilde{{\cal A}}_{\tilde{M},F}^{\vee}$; celle-ci a un nombre fini d'hyperplans polaires d'\'equations $e^{<\lambda,\check{\alpha}>}=c$, pour $\alpha\in \Sigma(A_{\tilde{M}})$.

 Le signification de ces assertions a \'et\'e expliqu\'ee en 1.2.

  Preuve. Relevons ${\bf f}$ en un \'el\'ement de $C_{c}^{\infty}(\tilde{G}(F))\otimes Mes(G(F))$. En utilisant la d\'efinition 1.5(1), il suffit par r\'ecurrence de prouver que la fonction
 $\phi_{\tilde{M}}^{\tilde{G}}({\bf f}) $
 v\'erifie les m\^emes propri\'et\'es. Or, soit $\tilde{\pi}$ une $\omega$-repr\'esentation temp\'er\'ee de $\tilde{M}(F)$. Par d\'efinition, la fonction $X\mapsto I^{\tilde{M}}(\tilde{\pi},X,\phi_{\tilde{M}}^{\tilde{G}}({\bf f}))$ est la transform\'ee de Fourier de la fonction $\tilde{\lambda}\mapsto J_{\tilde{M}}^{\tilde{G}}(\tilde{\pi}_{\tilde{\lambda}},{\bf f})$ sur $i\tilde{{\cal A}}_{\tilde{M}}^*/i\tilde{{\cal A}}_{\tilde{M},F}^{\vee}$. Celle-ci est  $C^{\infty}$ et se prolonge en une fonction  rationnelle sur $\tilde{{\cal A}}_{\tilde{M},{\mathbb C}}^*/i\tilde{{\cal A}}_{\tilde{M},F}^{\vee}$ dont les p\^oles ont la propri\'et\'e indiqu\'ee en (2). Cela prouve les assertions (1) et (2). $\square$

 Comme on l'a dit en 1.7, on  note encore 
 $$\tilde{\lambda}\mapsto  I^{\tilde{M}}(\tilde{\pi},\tilde{\lambda},{^c\theta}^{\tilde{G}}_{\tilde{M}}({\bf f}))$$
 la fonction rationnelle sur $\tilde{{\cal A}}_{\tilde{M},{\mathbb C}}^*/i\tilde{{\cal A}}_{\tilde{M},F}^{\vee}$  qui prolonge la fonction not\'ee ainsi qui est d\'ej\`a d\'efinie sur $i\tilde{{\cal A}}_{\tilde{M}}^*/i\tilde{{\cal A}}_{\tilde{M},F}^{\vee}$. 
 
   Soit $\nu\in {\cal A}_{\tilde{M}}^*$ un \'el\'ement tel que cette fonction  n'ait pas de p\^ole sur l'ensemble des \'el\'ements  $\tilde{\lambda}\in \tilde{{\cal A}}_{\tilde{M},{\mathbb C}}^*$  tels que $Re(\lambda)=\nu$. Pour $X\in \tilde{{\cal A}}_{\tilde{M},F}$, on pose
 $$I^{\tilde{M}}(\tilde{\pi},\nu,X,{^c\theta}^{\tilde{G}}_{\tilde{M}}({\bf f}))=\int_{\nu+i{\cal A}_{\tilde{M},F}^*} I^{\tilde{M}}(\tilde{\pi},\tilde{\lambda},{^c\theta}^{\tilde{G}}_{\tilde{M}}({\bf f}))e^{-<\tilde{\lambda},X>}\,d\lambda.$$
 On a la formule d'inversion
 $$ I^{\tilde{M}}(\tilde{\pi},\tilde{\lambda},{^c\theta}^{\tilde{G}}_{\tilde{M}}({\bf f}))=\sum_{X\in \tilde{{\cal A}}_{\tilde{M},F}}I^{\tilde{M}}(\tilde{\pi},\nu,X,{^c\theta}^{\tilde{G}}_{\tilde{M}}({\bf f}))e^{<\tilde{\lambda},X>}$$
 pour tout $\tilde{\lambda}\in \tilde{{\cal A}}_{\tilde{M}}^*/i\tilde{{\cal A}}_{\tilde{M},F}^{\vee}$ tel que $Re(\lambda)=\nu$.

 \ass{Proposition}{Supposons $\tilde{M}\not=\tilde{G}$ et $\tilde{\pi}$ elliptique.  Pour tout $\tilde{S}\in {\cal P}(\tilde{M})$, fixons un point $\nu_{\tilde{S}}$ comme en 1.3. Supposons ce point  "assez positif" pour $\tilde{S}$, cette notion d\'ependant de la repr\'esentation $\tilde{\pi}$. Alors on a l'\'egalit\'e
 $$\sum_{\tilde{S}\in {\cal P}(\tilde{M})}\omega_{\tilde{S}}(X)I^{\tilde{M}}(\tilde{\pi},\nu_{\tilde{S}},X,{^c\theta}^{\tilde{G}}_{\tilde{M}}({\bf f}))=0$$
 pour tout ${\bf f}\in I(\tilde{G}(F),\omega)\otimes Mes(G(F))$ et tout $X\in \tilde{{\cal A}}_{\tilde{M},F}$.}
 
 Preuve. On rel\`eve ${\bf f}$ en un \'el\'ement de $C_{c}^{\infty}(\tilde{G}(F))\otimes Mes(G(F))$. On utilise la d\'efinition 1.5(1). Il r\'esulte de la preuve de (1) et (2)  ci-dessus que
$$ I^{\tilde{M}}(\tilde{\pi},\tilde{\lambda},{^c\theta}^{\tilde{G}}_{\tilde{M}}({\bf f}))= I^{\tilde{M}}(\tilde{\pi},\tilde{\lambda}, \phi^{\tilde{G}}_{\tilde{M}}({\bf f}))-\sum_{\tilde{L}\in {\cal L}(\tilde{M}),\tilde{L}\not=\tilde{G}} I^{\tilde{M}}(\tilde{\pi},\tilde{\lambda},{^c\theta}^{\tilde{L}}_{\tilde{M}}(^c\phi_{\tilde{L}}^{\tilde{G}}({\bf f}))).$$
   On peut imiter les constructions ci-dessus pour chacune de ces fonctions et on utilise des notations similaires. On obtient que la somme de l'\'enonc\'e est \'egale \`a la somme de 
$$(3)\qquad \sum_{\tilde{S}\in {\cal P}(\tilde{M})}\omega_{\tilde{S}}(X)I^{\tilde{M}}(\tilde{\pi},\nu_{\tilde{S}},X,\phi^{\tilde{G}}_{\tilde{M}}({\bf f}))$$
et, pour tout $\tilde{L}\in {\cal P}(\tilde{M})$ tel que $\tilde{L}\not=\tilde{G}$, de 
$$(4) \qquad -\sum_{\tilde{S}\in {\cal P}(\tilde{M})}\omega_{\tilde{S}}(X)I^{\tilde{M}}(\tilde{\pi},\nu_{\tilde{S}},X,{^c\theta}^{\tilde{L}}_{\tilde{M}}(^c\phi_{\tilde{L}}^{\tilde{G}}({\bf f}))).$$

Fixons $\tilde{L}\not=\tilde{G}$. Par un argument d\'ej\`a utilis\'e plusieurs fois, pour $\tilde{S}\in {\cal P}(\tilde{M})$, on peut remplacer le point $\nu_{\tilde{S}}$ par $\nu_{\tilde{S}'}$, o\`u $\tilde{S}'=\tilde{S}\cap \tilde{L}$.  On peut ensuite, pour $\tilde{S}'$ fix\'e, remplacer la somme des $\omega_{\tilde{S}}(X)$ pour $\tilde{S}\cap\tilde{L}=\tilde{S}'$ par $\omega_{\tilde{S}'}(X)$. On obtient que la somme (4) est \'egale \`a l'oppos\'ee de celle de l'\'enonc\'e o\`u l'on remplace $\tilde{G}$ par $\tilde{L}$ et ${\bf f}$ par $^c\phi_{\tilde{L}}^{\tilde{G}}({\bf f})$. Dans l'\'enonc\'e figure l'hypoth\`ese $\tilde{M}\not=\tilde{G}$. Si $\tilde{L}\not=\tilde{M}$, elle reste v\'erifi\'ee et, par r\'ecurrence, la somme (4) est nulle. 

Consid\'erons l'espace $\tilde{L}=\tilde{M}$. Dans ce cas $^c\theta_{\tilde{M}}^{\tilde{M}}$ est l'identit\'e.  On a donc
$$ I^{\tilde{M}}(\tilde{\pi},\tilde{\lambda},{^c\theta}^{\tilde{M}}_{\tilde{M}}(^c\phi_{\tilde{M}}^{\tilde{G}}({\bf f})))=I^{\tilde{M}}(\tilde{\pi},\tilde{\lambda}, {^c\phi}_{\tilde{M}}^{\tilde{G}}({\bf f})).$$
Or $^c\phi_{\tilde{M}}^{\tilde{G}}({\bf f})$ est \`a support compact. Il en r\'esulte que cette fonction n'a pas de p\^ole. Ses coefficients de Fourier $I^{\tilde{M}}(\tilde{\pi},\nu,X, ^c\phi_{\tilde{M}}^{\tilde{G}}({\bf f}))$ sont donc ind\'ependants du point $\nu$. Ils sont \'egaux \`a
$$I^{\tilde{M}}(\tilde{\pi},X,{^c \phi}^{\tilde{G}}_{\tilde{M}}({\bf f}))=I^{\tilde{M}}(\tilde{\pi},0,X,{^c \phi}^{\tilde{G}}_{\tilde{M}}({\bf f}))=\int_{i{\cal A}_{\tilde{M},F}^*} I^{\tilde{M}}(\tilde{\pi}_{\tilde{\lambda}},{^c\phi}^{\tilde{G}}_{\tilde{M}}({\bf f}))e^{-<\tilde{\lambda},X>}\,d\lambda.$$
Dans la somme (4), on peut remplacer tous les points $\nu_{\tilde{S}}$ par $0$. Cette somme devient
$$-I^{\tilde{M}}(\tilde{\pi},X,{^c \phi}^{\tilde{G}}_{\tilde{M}}({\bf f}))\sum_{\tilde{S}\in {\cal P}(\tilde{M})}\omega_{\tilde{S}}(X),$$
qui est simplement $-I^{\tilde{M}}(\tilde{\pi},X,{^c \phi}^{\tilde{G}}_{\tilde{M}}({\bf f}))$. 
Par d\'efinition de $^c\phi_{\tilde{M}}^{\tilde{G}}({\bf f})$, on a
$$ I^{\tilde{M}}(\tilde{\pi}_{\tilde{\lambda}},{^c\phi}^{\tilde{G}}_{\tilde{M}}({\bf f}))={^cJ}_{\tilde{M}}^{\tilde{G}}(\tilde{\pi}_{\tilde{\lambda}},{\bf f})$$
pour $\tilde{\lambda}\in i\tilde{{\cal A}}_{\tilde{M}}^*/i\tilde{{\cal A}}_{\tilde{M},F}^{\vee}$. 
Parce que $\tilde{\pi}$ est elliptique, ceci n'est autre que
$$\sum_{Y\in \tilde{{\cal A}}_{\tilde{M},F}}\sum_{\tilde{S}\in {\cal P}(\tilde{M})}\omega_{\tilde{S}}(Y)\int_{\nu_{\tilde{S}}+i{\cal A}_{\tilde{M},F}^*}J_{\tilde{M}}^{\tilde{G}}(\tilde{\pi}_{\tilde{\lambda}+\tilde{\nu}},{\bf f})e^{-<\tilde{\nu},Y>}\,d\nu.$$
Par inversion de Fourier, on en d\'eduit
$$(5) \qquad I^{\tilde{M}}(\tilde{\pi},X,{^c \phi}^{\tilde{G}}_{\tilde{M}}({\bf f}))=\sum_{\tilde{S}\in {\cal P}(\tilde{M})}\omega_{\tilde{S}}(X)\int_{\nu_{\tilde{S}}+i{\cal A}_{\tilde{M},F}^*}J_{\tilde{M}}^{\tilde{G}}(\tilde{\pi}_{\tilde{\nu}},{\bf f})e^{-<\tilde{\nu},X>}\,d\nu.$$
En r\'esum\'e, la somme sur les $\tilde{L}\not=\tilde{G}$ des termes (4) est \'egale \`a l'oppos\'e du  membre de droite de (5).

Ainsi qu'il r\'esulte de la preuve de (1) et (2),  on a l'\'egalit\'e
$$ I^{\tilde{M}}(\tilde{\pi},\tilde{\lambda},\phi^{\tilde{G}}_{\tilde{M}}({\bf f}))=J_{\tilde{M}}^{\tilde{G}}(\tilde{\pi}_{\tilde{\lambda}},{\bf f})$$
pour tout $\tilde{\lambda}\in {\cal A}_{\tilde{M},{\mathbb C}}^*/i\tilde{{\cal A}}_{\tilde{M},F}^{\vee}$. Il r\'esulte alors des d\'efinitions que
$$I^{\tilde{M}}(\tilde{\pi},\nu_{\tilde{S}},X,\phi^{\tilde{G}}_{\tilde{M}}({\bf f}))=\int_{\nu_{\tilde{S}}+i{\cal A}_{\tilde{M},F}^*}J_{\tilde{M}}^{\tilde{G}}(\tilde{\pi}_{\tilde{\nu}},{\bf f})e^{-<\tilde{\nu},X>}\,d\nu$$
pour tout $\tilde{S}$. Donc la somme (3) est \'egale au membre de droite de (5). La somme de (3) et (4) est donc nulle. Cela ach\`eve la d\'emonstration. $\square$

\bigskip

\subsection{Une variante des int\'egrales orbitales pond\'er\'ees $\omega$-\'equivariantes}
Fixons un espace de Levi $\tilde{M}$. L'application $^c\phi_{\tilde{M}}^{\tilde{G}}$ d\'efinie en 1.3 v\'erifie les m\^emes propri\'et\'es formelles que $\phi_{\tilde{M}}^{\tilde{G}}$.  Dans beaucoup de constructions que l'on a faites, on peut remplacer la deuxi\`eme application par la premi\`ere. Ainsi, fixons pour un moment un sous-groupe compact sp\'ecial $K$ de $G(F)$ en bonne position relativement \`a $M$. Soit $\boldsymbol{\gamma}\in D_{g\acute{e}om}(\tilde{M}(F),\omega)\otimes Mes(M(F))^*$. On a d\'efini en [II] 1.5 la distribution
$${\bf f}\mapsto J_{\tilde{M}}^{\tilde{G}}(\boldsymbol{\gamma},{\bf f})$$
sur $C_{c}^{\infty}(\tilde{G}(F))\otimes Mes(G(F))$. On d\'efinit une nouvelle distribution par la formule
$$(1) \qquad ^cI_{\tilde{M}}^{\tilde{G}}(\boldsymbol{\gamma},{\bf f})=J_{\tilde{M}}^{\tilde{G}}(\boldsymbol{\gamma},{\bf f})-\sum_{\tilde{L}\in {\cal L}(\tilde{M}),\tilde{L}\not=\tilde{G}}{^cI}_{\tilde{M}}^{\tilde{L}}(\boldsymbol{\gamma},{^c\phi}_{\tilde{L}}^{\tilde{G}}({\bf f})).$$
Elle v\'erifie les m\^emes propri\'et\'es formelles que la distribution $I_{\tilde{M}}^{\tilde{G}}(\boldsymbol{\gamma},{\bf f})$. En particulier, elle est $\omega$-\'equivariante, c'est-\`a-dire se descend en une distribution sur $I(\tilde{G}(F),\omega)\otimes Mes(G(F))$. Elle est aussi compatible \`a l'induction, c'est-\`a-dire v\'erifie une propri\'et\'e analogue au lemme [II] 1.7. Mais elle v\'erifie une propri\'et\'e utile que la distribution $I_{\tilde{M}}^{\tilde{G}}(\boldsymbol{\gamma},{\bf f})$ ne v\'erifie pas. Pour tout sous-ensemble $\Omega$ de $\tilde{M}(F)$, posons $\Omega^{M}=\{m^{-1}\gamma m; m\in M(F), \gamma\in \Omega\}$. On a

(2) pour tout ${\bf f}\in I(\tilde{G}(F),\omega)\otimes Mes(G(F))$, il existe un sous-ensemble compact $\Omega\subset \tilde{M}(F)$ tel que $^cI_{\tilde{M}}^{\tilde{G}}(\boldsymbol{\gamma},{\bf f})=0$ si le support de $\boldsymbol{\gamma}$ ne coupe pas $\Omega^M$.

Preuve. On rel\`eve ${\bf f}$ en un \'el\'ement de $C_{c}^{\infty}(\tilde{G}(F))\otimes Mes(G(F))$. La formule (1) ram\`ene par r\'ecurrence \`a prouver que le terme $J_{\tilde{M}}^{\tilde{G}}(\boldsymbol{\gamma},{\bf f})$ v\'erifie la m\^eme propri\'et\'e. Notons $\Omega'$ le support de ${\bf f}$. Il existe un sous-ensemble compact $\Omega$ de $\tilde{M}(F)$ tel que $\Omega^M=(\Omega')^G\cap \tilde{M}(F)$. Il est clair que $J_{\tilde{M}}^{\tilde{G}}(\boldsymbol{\gamma},{\bf f})=0$ si le support de $\boldsymbol{\gamma}$ ne coupe pas $\Omega^M$. $\square$

Le lien entre nos deux distributions $\omega$-\'equivariantes est \'etabli par la formule suivante
$$(3)\qquad ^cI_{\tilde{M}}^{\tilde{G}}(\boldsymbol{\gamma},{\bf f})=\sum_{\tilde{L}\in {\cal L}(\tilde{M})}I_{\tilde{M}}^{\tilde{L}}(\boldsymbol{\gamma},{^c\theta}_{\tilde{L}}^{\tilde{G}}({\bf f})).$$

Preuve. Relevons ${\bf f}$ en un \'el\'ement de $C_{c}^{\infty}(\tilde{G}(F))\otimes Mes(G(F))$ et appliquons cette formule par r\'ecurrence au termes de la somme du membre de droite de la formule (1). On obtient 
$$^cI_{\tilde{M}}^{\tilde{G}}(\boldsymbol{\gamma},{\bf f})=J_{\tilde{M}}^{\tilde{G}}(\boldsymbol{\gamma},{\bf f})-\sum_{\tilde{L}_{1}\in {\cal L}(\tilde{M}),\tilde{L}_{1}\not=\tilde{G}}\sum_{\tilde{L}_{2}\in {\cal L}(\tilde{M}), \tilde{L}_{2}\subset \tilde{L}_{1}}I_{\tilde{M}}^{\tilde{L}_{2}}(\boldsymbol{\gamma},{^c\theta}_{\tilde{L}_{2}}^{\tilde{L}_{1}}(^c\phi_{\tilde{L}_{1}}^{\tilde{G}}({\bf f}))).$$
Renversons l'ordre de sommation. On obtient une somme sur $\tilde{L}_{2}\in {\cal L}(\tilde{M})$  de la distribution $I_{\tilde{M}}^{\tilde{L}_{2}}(\boldsymbol{\gamma},.)$ appliqu\'ee \`a la fonction 
$$\sum_{\tilde{L}_{1}\in {\cal L}(\tilde{L}_{2}), \tilde{L}_{1}\not=\tilde{G}}{^c\theta}_{\tilde{L}_{2}}^{\tilde{L}_{1}}(^c\phi_{\tilde{L}_{1}}^{\tilde{G}}({\bf f})).$$
Par d\'efinition de $^c\theta_{\tilde{L}_{2}}^{\tilde{G}}$, cette fonction est \'egale \`a $\phi_{\tilde{L}_{2}}^{\tilde{G}}({\bf f})-{^c\theta}_{\tilde{L}_{2}}^{\tilde{G}}({\bf f})$. En rempla\c{c}ant $\tilde{L}_{2}$ simplement par $\tilde{L}$, on obtient alors que 
$^cI_{\tilde{M}}^{\tilde{G}}(\boldsymbol{\gamma},{\bf f})$ est la somme de
$$J_{\tilde{M}}^{\tilde{G}}(\boldsymbol{\gamma},{\bf f}),$$
de 
$$-\sum_{\tilde{L}\in {\cal L}(\tilde{M})}I_{\tilde{M}}^{\tilde{L}}(\boldsymbol{\gamma},\phi_{\tilde{L}}^{\tilde{G}}({\bf f}))$$
et de 
$$\sum_{\tilde{L}\in {\cal L}(\tilde{M})}I_{\tilde{M}}^{\tilde{L}}(\boldsymbol{\gamma},{^c\theta}_{\tilde{L}}^{\tilde{G}}({\bf f})).$$
La somme des deux premi\`eres  expressions  est nulle par d\'efinition de la distribution $I_{\tilde{M}}^{\tilde{G}}(\boldsymbol{\gamma},.)$. Donc  $^cI_{\tilde{M}}^{\tilde{G}}(\boldsymbol{\gamma},{\bf f})$ est \'egal \`a la derni\`ere somme, ce qui est l'\'egalit\'e (3). $\square$

\bigskip

\section{Stabilisation de l'application $^c\theta_{\tilde{M}}$}

\bigskip

\subsection{Fonctions $\omega_{\tilde{S}}$ et endoscopie}
Soit ${\bf G}'=(G',{\cal G}',\tilde{s})$ une donn\'ee endoscopique elliptique et relevante de $(G,\tilde{G},{\bf a})$. Introduisons $\underline{les}$ paires de Borel $(B^*,T^*)$ et $(B^{_{'}*},T^{_{'}*})$ de $G$ et $G'$. Modulo certains choix dans les groupes duaux,  cf. [I] 1.5, on a un homomorphisme $\xi:T^*\to T^{_{'}*}$. Il se restreint en un homomorphisme $\xi:A_{\tilde{G}}\to A_{\tilde{G}'}$ qui est \'equivariant pour les actions galoisiennes. L'hypoth\`ese d'ellipticit\'e signifie que cet homomorphisme est un rev\^etement. Il s'en d\'eduit un isomorphisme ${\cal A}_{\tilde{G}}\simeq {\cal A}_{\tilde{G}'}$. On a

(1) il existe un unique isomorphisme $\tilde{\xi}:\tilde{{\cal A}}_{\tilde{G}}\simeq \tilde{{\cal A}}_{\tilde{G}'}$ tel que

(i) $\tilde{\xi}(X+H)=\xi(X)+\tilde{\xi}(H)$ pour tout $X\in {\cal A}_{\tilde{G}}$ et $H\in \tilde{{\cal A}}_{\tilde{G}}$;

(ii) pour tout couple $(\delta,\gamma)$ form\'e d'un \'el\'ement semi-simple $\delta\in \tilde{G}'(F)$ et d'un \'el\'ement semi-simple $\gamma\in \tilde{G}(F)$ qui se correspondent (cf. [I] 1.10), on ait l'\'egalit\'e
$$\tilde{\xi}(\tilde{H}_{\tilde{G}}(\gamma))=\tilde{H}_{\tilde{G}'}(\delta).$$

La preuve est la m\^eme qu'en [IV] 2.1(1). Signalons qu'en g\'en\'eral, les ensembles $\tilde{\xi}(\tilde{{\cal A}}_{\tilde{G},F})$ et $\tilde{{\cal A}}_{\tilde{G}',F}$ ne sont pas inclus l'un dans l'autre. Leur intersection est toutefois non vide d'apr\`es (ii) ci-dessus. En prenant pour point-base un \'el\'ement de cette intersection, ces ensembles s'identifient respectivement aux r\'eseaux $\xi({\cal A}_{\tilde{G},F})$ et ${\cal A}_{\tilde{G}',F}$. Ceux-ci ne sont pas non plus inclus l'un dans l'autre, mais sont commensurables (leur intersection est d'indice finie dans chacun d'eux).

Dans la suite, on abandonne les notations $\xi$ et $\tilde{\xi}$ et on identifie simplement ${\cal A}_{\tilde{G}}$ \`a ${\cal A}_{\tilde{G}'}$ par $\xi$, ainsi que $\tilde{{\cal A}}_{\tilde{G}}$ \`a $\tilde{{\cal A}}_{\tilde{G}'}$ par $\tilde{\xi}$.

Soit $M'$ un Levi de $G'$. Supposons $M'$ relevant. Il lui correspond donc un espace de Levi $\tilde{M}$ de $\tilde{G}$ et $M'$ se compl\`ete en une donn\'ee endoscopique elliptique  ${\bf M}'=(M',{\cal M}',\tilde{s})$ de $(M,\tilde{M},{\bf a})$. On a comme ci-dessus des identifications compatibles ${\cal A}_{\tilde{M}}\simeq {\cal A}_{\tilde{M}'}$ et $\tilde{{\cal A}}_{\tilde{M}}\simeq \tilde{{\cal A}}_{\tilde{M}'}$. Soit $\tilde{P}'\in {\cal P}(\tilde{M}')$. Des descriptions des ensembles de racines de nos diff\'erents groupes r\'esulte que  la cl\^oture de la chambre positive ${\cal A}_{\tilde{P}'}$ est r\'eunion de cl\^otures de chambres positives ${\cal A}_{\tilde{P}}$ pour certains $\tilde{P}\in {\cal P}(\tilde{M})$. Notons simplement $\tilde{P}\mapsto \tilde{P}'$ pour signifier que $\tilde{P}$ appartient \`a cet ensemble. On a fix\'e des fonctions $\omega_{\tilde{P}}$. On d\'efinit la fonction $\omega_{\tilde{P}'}$ sur $\tilde{{\cal A}}_{\tilde{M}'}$ par 
$$(2) \qquad \omega_{\tilde{P}'}=\sum_{\tilde{P}\in {\cal P}(\tilde{M});\tilde{P}\mapsto \tilde{P}'}\omega_{\tilde{P}}.$$
L'espace de Levi $\tilde{M}$ n'est pas unique, il n'est bien d\'efini qu'\`a conjugaison pr\`es, mais les propri\'et\'es d'invariance impos\'ees aux fonctions $\omega_{\tilde{P}}$ entra\^{\i}nent que la d\'efinition ci-dessus ne d\'epend pas du choix de $\tilde{M}$. 
Il est imm\'ediat que l'ensemble de fonctions ainsi d\'efinies v\'erifie les m\^emes conditions qu'en 1.1. 

Dans le cas o\`u $M'$ n'est pas relevant, on fixe sans r\'ef\'erence \`a $\tilde{G}$  des fonctions $\omega_{\tilde{P}'}$ pour $\tilde{P}'\in {\cal P}(\tilde{M}')$ qui v\'erifient les conditions de 1.1.

\bigskip

\subsection{Les applications $^cS\theta_{\tilde{M}}^{\tilde{G}}$}
Pour la suite de la section, on suppose $(G,\tilde{G},{\bf a})$ quasi-d\'eploy\'e et \`a torsion int\'erieure et on fixe  un espace de Levi $\tilde{M}$ de $\tilde{G}$. Nous allons d\'efinir une application lin\'eaire
$$^cS\theta_{\tilde{M}}^{\tilde{G}}:I(\tilde{G}(F))\otimes Mes(G(F))\to SI_{ac}(\tilde{M}(F))\otimes Mes(M(F)).$$
Comme toujours, on doit admettre par r\'ecurrence certaines de ses propri\'et\'es. Les premi\`eres sont formelles et permettent de d\'efinir des applications analogues pour des donn\'ees endoscopiques. On y revient ci-dessous. La seconde est que cette application est stable, c'est-\`a-dire qu'elle se factorise en une application d\'efinie sur $SI(\tilde{G}(F))\otimes Mes(G(F))$. Notons $p_{\tilde{M}}^{st}:I_{ac}(\tilde{M}(F))\otimes Mes(M(F))\to SI_{ac}(\tilde{M}(F))\otimes Mes(M(F))$ la projection naturelle. Pour ${\bf f}\in I(\tilde{G}(F))\otimes Mes(G(F))$, on pose
$$(1)\qquad  {^cS}\theta_{\tilde{M}}^{\tilde{G}}({\bf f})=p_{\tilde{M}}^{st}\circ{^c\theta}_{\tilde{M}}^{\tilde{G}}({\bf f})-\sum_{s\in Z(\hat{M})^{\Gamma_{F}}/Z(\hat{G})^{\Gamma_{F}},s\not=1}i_{\tilde{M}}(\tilde{G},\tilde{G}'(s)){^cS}\theta_{{\bf M}}^{{\bf G}'(s)}({\bf f}^{{\bf G}'(s)}).$$
Enon\c{c}ons la propri\'et\'e \'evoqu\'ee ci-dessus.

\ass{Proposition }{L'application $^cS\theta_{\tilde{M}}^{\tilde{G}}$ se quotiente en une application lin\'eaire
$$SI(\tilde{G}(F))\otimes Mes(G(F))\to SI_{ac}(\tilde{M}(F))\otimes Mes(M(F)).$$}

Cette proposition sera prouv\'ee en 4.2.

Revenons sur les questions formelles. Consid\'erons des extensions compatibles
$$1\to C_{\natural}\to G_{\natural}\to G\to 1\text{ et } \tilde{G}_{\natural}\to \tilde{G}$$
o\`u $C_{\natural}$ est un tore induit et $\tilde{G}_{\natural}$ est encore \`a torsion int\'erieure. Soit $\lambda_{\natural}$ un caract\`ere de $C_{\natural}(F)$. On a une projection naturelle ${\cal A}_{\tilde{M}_{\natural}}\to {\cal A}_{\tilde{M}}$ et une identification ${\cal P}(\tilde{M}_{\natural})={\cal P}(\tilde{M})$, notons-la $\tilde{P}_{\natural}\leftrightarrow \tilde{P}$. Pour $\tilde{P}\in {\cal P}_{\tilde{M}}$, on a une fonction $\omega_{\tilde{P}}$ sur ${\cal A}_{\tilde{M}}$. Pour $\tilde{P}_{\natural}\leftrightarrow \tilde{P}$, on choisit pour fonction $\omega_{\tilde{P}_{\natural}}$ la compos\'ee de cette fonction avec la projection pr\'ec\'edente. 
Fixons des mesures de Haar $dc$ sur $C_{\natural}(F)$, $dg$   sur $G(F)$ et $dm$   sur $M(F)$. Soient ${\bf f}\in C_{c,\lambda_{\natural}}^{\infty}(\tilde{G}_{\natural}(F))\otimes Mes(G(F))$. Ecrivons-la ${\bf f}=f\otimes dg$.  Choisissons une fonction $\dot{f}\in C_{c}^{\infty}(\tilde{G}_{\natural} (F))$ telle que
$$f(\gamma_{\natural})=\int_{C_{\natural}(F)}\dot{f}(c\gamma_{\natural})\lambda_{\natural}(c)\,dc$$
pour tout $\gamma_{\natural}\in \tilde{G}_{\natural}(F)$. De $dg$ et $dc$ se d\'eduit une mesure de Haar $dg_{\natural}$ sur $G_{\natural}(F)$. De m\^eme, de $dm$ et $dc$ se d\'eduit une mesure de Haar $dm_{\natural}$ sur $M_{\natural}(F)$. On d\'efinit alors $^cS\theta_{\tilde{M}_{\natural}}^{\tilde{G}_{\natural}}(\dot{f}\otimes dg_{\natural})$. Ecrivons-la $\dot{\varphi}\otimes dm_{\natural}$ avec $\dot{\varphi}\in SI_{ac}(\tilde{M}_{\natural}(F))$. Cette fonction n'est pas \`a support compact. Toutefois, on v\'erifie ais\'ement par r\'ecurrence sur la d\'efinition (1) que l'application $^cS\theta_{\tilde{M}_{\natural}}^{\tilde{G}_{\natural}}$ v\'erifie la relation 1.6(3) (la v\'erification pr\'ecise n\'ecessite encore quelques consid\'erations formelles que l'on passe). Il en r\'esulte que, pour tout $\gamma_{\natural}\in \tilde{M}_{\natural}(F)$, la fonction $c\mapsto \dot{\varphi}(c\gamma_{\natural})$ est \`a support compact. On peut alors d\'efinir $\varphi\in SI_{ac,\lambda_{\natural}}(M_{\natural}(F))$ par
$$\varphi(\gamma_{\natural})=\int_{C_{\natural}(F)}\dot{\varphi}(c\gamma_{\natural})\lambda_{\natural}(c)\,dc.$$

{\bf Remarque.} L'espace $SI_{ac,\lambda_{\natural}}(M_{\natural}(F))$ se d\'eduit de $C_{ac,\lambda_{\natural}}^{\infty}(M_{\natural}(F))$. Ce dernier est celui des fonctions $\psi:M_{\natural}(F)\to {\mathbb C}$ qui se transforment par $C_{\natural}(F)$ selon  le caract\`ere $\lambda_{\natural}$, qui sont biinvariantes par un sous-groupe ouvert compact de $M_{\natural}(F)$ et qui v\'erifient  la condition suivante. Soit $b$ une fonction \`a support compact sur $\tilde{{\cal A}}_{\tilde{M},F}$. Notons $b_{\natural}$ sa compos\'ee avec la projection $\tilde{{\cal A}}_{\tilde{M}_{\natural},F}\to \tilde{{\cal A}}_{\tilde{M},F}$. Alors la fonction $(b_{\natural}\circ \tilde{H}_{\tilde{M}_{\natural}})\psi$ appartient \`a $C_{c,\lambda_{\natural}}^{\infty}(\tilde{M}_{\natural}(F))$.

 \bigskip

Le terme $\varphi\otimes dm$ d\'epend a priori des  choix de mesures et du choix de la fonction $\dot{f}$. On v\'erifie facilement que changer de mesures ne modifie pas $\varphi\otimes dm$. Changer de $\dot{f}$ revient \`a ajouter \`a cette fonction une combinaison lin\'eaire de fonctions $\psi^c-\lambda_{\natural}(c)\psi$, o\`u $\psi\in C_{c}^{\infty}(\tilde{G}_{\natural}(F))$ et $c\in C_{\natural}(F)$. Notre syst\`eme de fonctions $\omega_{\tilde{P}}$ sur ${\cal A}_{\tilde{M}_{\natural}}$ v\'erifie par construction la condition 1.4(7) pour $Z=C_{\natural}$. On voit par r\'ecurrence que l'application $^cS\theta_{\tilde{M}_{\natural}}^{\tilde{G}_{\natural}}$ v\'erifie la propri\'et\'e 1.6(4). Il en r\'esulte qu'ajouter $\psi^c-\lambda_{\natural}(c)\psi$ \`a $\dot{f}$ revient \`a ajouter \`a $\dot{\varphi}$ une fonction de la forme $(\psi')^c-\lambda_{\natural}(c)\psi'$. Un tel terme dispara\^{\i}t par int\'egration et ne change pas $\varphi$. Ainsi le terme $\varphi\otimes dm$ est bien d\'efini et on peut poser
$$^cS\theta_{\tilde{M}_{\natural},\lambda_{\natural}}^{\tilde{G}_{\natural}}({\bf f})=\varphi\otimes dm.$$
Consid\'erons d'autres donn\'ees
$$1\to C_{\flat}\to G_{\flat}\to G\to 1 , \,\,\tilde{G}_{\flat}\to \tilde{G},\,\, \lambda_{\flat}$$
v\'erifiant les m\^emes hypoth\`eses. Introduisons le produit fibr\'e $G_{\natural,\flat}$ de $G_{\natural}$ et $G_{\flat}$ au-dessus de $G$ et le produit fibr\'e $\tilde{G}_{\natural,\flat}$ de $\tilde{G}_{\natural}$ et $\tilde{G}_{\flat}$ au-dessus de $\tilde{G}$. Supposons donn\'es un caract\`ere $\lambda_{\natural,\flat}$ de $G_{\natural,\flat}(F)$ dont la restriction \`a $C_{\natural}(F)\times C_{\flat}(F)$ est $\lambda_{\natural}\times \lambda_{\flat}^{-1}$ et une fonction non nulle $\tilde{\lambda}_{\natural,\flat}$ sur $\tilde{G}_{\natural,\flat}(F)$ qui se transforme selon le caract\`ere $\lambda_{\natural,\flat}$. On en d\'eduit un isomorphisme
$$C_{c,\lambda_{\natural}}^{\infty}(\tilde{G}_{\natural}(F))\otimes Mes(G(F))\simeq C_{c, \flat}^{\infty}(G_{\flat}(F))\otimes Mes(G(F)),$$
cf. [II] 1.10. On en d\'eduit de m\^eme un isomorphisme
$$SI_{ac,\lambda_{\natural}}(M_{\natural}(F))\otimes Mes(M(F))\simeq SI_{ac,\lambda_{\flat}}(M_{\flat}(F))\otimes Mes(M(F)).$$
Soient ${\bf f}_{\natural}\in C_{c,\lambda_{\natural}}^{\infty}(\tilde{G}_{\natural}(F))\otimes Mes(G(F))$ et ${\bf f}_{\flat}\in C_{c,\lambda_{\flat}}^{\infty}(\tilde{G}_{\flat}(F))\otimes Mes(G(F))$ qui se correspondent par le premier isomorphisme. Alors $^cS\theta_{\tilde{M}_{\natural},\lambda_{\natural}}^{\tilde{G}_{\natural}}({\bf f}_{\natural})$ et $^cS\theta_{\tilde{M}_{\flat},\lambda_{\flat}}^{\tilde{G}_{\flat}}({\bf f}_{\flat})$ se correspondent par le second. La preuve est similaire \`a celle de [II] 1.10(7). Rappelons-en seulement la structure. On commence par prouver que les deux termes sont \'egaux \`a des termes analogues relatifs aux extensions communes $G_{\natural,\flat}$ et $\tilde{G}_{\natural,\flat}$, o\`u le caract\`ere de $C_{\natural}(F)\times C_{\flat}(F)$ est $\lambda_{\natural}\otimes 1$ pour le premier terme et $1\otimes \lambda_{\flat}$ pour le second. On passe de l'un \`a l'autre par tensorisation par le caract\`ere affine $\tilde{\lambda}_{\natural,\flat}$. Cette tensorisation est compatible \`a nos constructions parce que l'on v\'erifie par r\'ecurrence que nos distributions v\'erifient l'analogue de 1.6(5).

 Cela \'etant, pour un \'el\'ement $s\in Z(\hat{M})^{\Gamma_{F}}/Z(\hat{G})^{\Gamma_{F}}$ tel que ${\bf G}'(s)$ est elliptique, on choisit des donn\'ees auxiliaires $G'_{1}(s),...,\Delta_{1}(s)$. On d\'efinit comme ci-dessus l'application lin\'eaire
 $$^cS\theta_{\tilde{M}_{1}(s),\lambda_{1}(s)}^{\tilde{G}'_{1}(s)}:C_{c,\lambda_{1}(s)}^{\infty}(\tilde{G}'_{1}(s;F))\otimes Mes(G(F))\to SI_{ac,\lambda_{1}(s)}(\tilde{M}_{1}(s;F))\otimes Mes(M(F)).$$
 La propri\'et\'e pr\'ec\'edente assure que, quand on change de donn\'ees auxiliaires, ces applications se recollent en une application lin\'eaire
 $$^cS\theta_{{\bf M}}^{{\bf G}'(s)}:C_{c}^{\infty}({\bf G}'(s))\otimes Mes(G(F))\to SI_{ac}({\bf M})\otimes Mes(M(F)).$$
 Elle se quotiente en tout cas en une application d\'efinie sur $I({\bf G}'(s))\otimes Mes(G(F))$ et, pourvu que la proposition soit d\'emontr\'ee pour $\tilde{G}'_{1}(s) $ et pour un choix de donn\'ees auxiliaires, en une application d\'efinie sur $SI({\bf G}'(s))\otimes Mes(G(F))$. 
 
 \bigskip
 
 \subsection{Commutation \`a l'induction}
 \ass{Lemme}{Soit $\tilde{R}$ un espace de Levi de $\tilde{M}$. Pour tout ${\bf f}\in I(\tilde{G}(F))\otimes Mes(G(F))$, on a l'\'egalit\'e 
$$(^cS\theta_{\tilde{M}}^{\tilde{G}}({\bf f}))_{\tilde{R}}= \sum_{\tilde{L}\in {\cal L}(\tilde{R})}e_{\tilde{R}}^{\tilde{G}}(\tilde{M},\tilde{L}){^cS}\theta_{\tilde{R}}^{\tilde{L}}({\bf f}_{\tilde{L}}).$$}

En utilisant le lemme 1.6, la preuve est similaire \`a celle de la proposition [II] 1.14(ii). $\square$

\bigskip

\subsection{Une propri\'et\'e d'annulation}
  Rappelons que l'on a d\'efini l'espace $D_{temp}^{st}(\tilde{M}(F))$  qui est le sous-espace des \'el\'ements de $D_{temp}(\tilde{M}(F))$ qui sont des distributions stables. On d\'efinit de m\^eme l'espace $D_{ell}^{st}(\tilde{M}(F))$. Fixons pour simplifier un espace de Levi minimal $\tilde{M}_{0}\subset \tilde{M}$. Moeglin a prouv\'e en [M] corollaire 2.1  que l'induction fournissait un isomorphisme
 $$D_{temp}^{st}(\tilde{M}(F))\otimes Mes(M(F))^*=$$
$$\oplus_{\tilde{R}\in {\cal L}^{\tilde{M}}(\tilde{M}_{0})/W^M(\tilde{M}_{0})}Ind_{\tilde{R}}^{\tilde{M}}((D_{ell}^{st}(\tilde{R}(F),\omega)\otimes Mes(R(F))^*)^{W^M(\tilde{R})}.$$

 Les d\'efinitions de 1.7 s'adaptent aux fonctions appartenant \`a $SI(\tilde{M}(F))\otimes Mes(M(F))$. Il suffit de se limiter dans les d\'efinitions de ce paragraphe  aux repr\'esentations appartenant \`a $D_{temp}^{st}(\tilde{M}(F))$.    Soit ${\bf f}\in I(\tilde{G}(F))\otimes Mes(G(F))$. On voit par r\'ecurrence que la fonction ${^cS\theta}_{\tilde{M}}^{\tilde{G}}({\bf f})$ est de Schwartz. 
Soit $\tilde{\pi}\in D_{temp}^{st}(\tilde{M}(F),\omega)\otimes Mes(M(F))^*$. On voit de m\^eme que la fonction 
 $\tilde{\lambda}\mapsto S^{\tilde{M}}(\tilde{\pi},\tilde{\lambda},{^cS\theta}_{\tilde{M}}^{\tilde{G}}({\bf f}))$ sur $i\tilde{{\cal A}}_{\tilde{M}}^*/i\tilde{{\cal A}}_{\tilde{M},F}^{\vee}$ v\'erifie la propri\'et\'e 1.8(2). On peut reprendre pour cette fonction les constructions et notations effectu\'ees dans ce paragraphe 1.8. 

\ass{Proposition}{Supposons $\tilde{M}\not=\tilde{G}$ et $\tilde{\pi}$ elliptique. Pour tout $\tilde{S}\in {\cal P}(\tilde{M})$, fixons un point $\nu_{\tilde{S}}$ comme en 1.3. Supposons ce point  "assez positif" pour $\tilde{S}$, cette notion d\'ependant de la repr\'esentation $\tilde{\pi}$. Alors on a l'\'egalit\'e
$$\sum_{\tilde{S}\in {\cal P}(\tilde{M})}\omega_{\tilde{S}}(X)S^{\tilde{M}}(\tilde{\pi},\nu_{\tilde{S}},X,{^cS\theta}_{\tilde{M}}^{\tilde{G}}({\bf f}))=0$$
pour tout $X\in \tilde{{\cal A}}_{\tilde{M},F}$.}

Preuve. Pour tous points $\nu$ et $X$, on a par d\'efinition l'\'egalit\'e 
$$S^{\tilde{M}}(\tilde{\pi},\nu,X,{^cS\theta}_{\tilde{M}}^{\tilde{G}}({\bf f}))=I^{\tilde{M}}(\tilde{\pi},\nu,X,{^c\theta}_{\tilde{M}}^{\tilde{G}}({\bf f}))-\sum_{s\in Z(\hat{M})^{\Gamma_{F}}/Z(\hat{G})^{\Gamma_{F}},s\not=1}S^{{\bf M}}(\tilde{\pi},\nu,X,{^cS\theta}_{{\bf M}}^{{\bf G}'(s)}({\bf f}^{{\bf G}'(s)})).$$
Le premier terme v\'erifie la relation de l'\'enonc\'e d'apr\`es la proposition 1.8. Fixons $s$ apparaissant ci-dessus. On doit prouver la relation de l'\'enonc\'e pour le terme index\'e par $s$ dans la somme ci-dessus. On a dit en 2.1 qu'il y avait une application naturelle $\tilde{S}\mapsto \tilde{S}'$ de ${\cal P}^{\tilde{G}}(\tilde{M})$ dans ${\cal P}^{\tilde{G}'(s)}(\tilde{M})$. Pour de tels $\tilde{S}$ et $\tilde{S}'$, la fonction $\tilde{\lambda}\mapsto S^{{\bf M}}(\tilde{\pi},\tilde{\lambda},{^cS}\theta_{{\bf M}}^{{\bf G}'(s)}({\bf f}^{{\bf G}'(s)}))$ n'a pas de p\^ole sur le segment joignant $\nu_{\tilde{S}}$ \`a $\nu_{\tilde{S}'}$. Cela nous permet de remplacer dans la relation de l'\'enonc\'e le point $\nu_{\tilde{S}}$ par $\nu_{\tilde{S}'}$. On vertu de la d\'efinition 2.1(2) de la fonction $\omega_{\tilde{S}'}$, la relation \`a d\'emontrer devient celle de l'\'enonc\'e pour $\tilde{G}'(s)$, ou plut\^ot ${\bf G}'(s)$, et la fonction ${\bf f}^{{\bf G}'(s)}$. Celle-ci est v\'erifi\'ee par r\'ecurrence puisque $s\not=1$. $\square$

\bigskip

\subsection{Une variante des int\'egrales orbitales pond\'er\'ees stables}
Soit $\boldsymbol{\delta}\in D_{g\acute{e}om}^{st}(\tilde{M}(F))\otimes Mes(M(F))^*$. On d\'efinit une forme lin\'eaire $^cS_{\tilde{M}}^{\tilde{G}}(\boldsymbol{\delta},.)$ par la formule habituelle
$$^cS_{\tilde{M}}^{\tilde{G}}(\boldsymbol{\delta},{\bf f})={^cI}_{\tilde{M}}^{\tilde{G}}(\boldsymbol{\delta},{\bf f})-\sum_{s\in Z(\hat{M})^{\Gamma_{F}}/Z(\hat{G})^{\Gamma_{F}},s\not=1}i_{\tilde{M}}(\tilde{G},\tilde{G}'(s)){^cS}_{{\bf M}}^{{\bf G}'(s)}(\boldsymbol{\delta},{\bf f}^{{\bf G}'(s)}).$$
Outre quelques formalit\'es que l'on passe, cette d\'efinition utilise par r\'ecurrence la propri\'et\'e suivante.

\ass{Proposition }{Pour tout $\boldsymbol{\delta}\in D_{g\acute{e}om}^{st}(\tilde{M}(F))\otimes Mes(M(F))^*$, la distribution
$${\bf f}\mapsto {^cS}_{\tilde{M}}^{\tilde{G}}(\boldsymbol{\delta},{\bf f})$$
est stable.}

Cette proposition sera prouv\'ee en 4.2.

Par r\'ecurrence, on voit que notre distribution v\'erifie la propri\'et\'e de compacite 1.9(2). 

 \bigskip

\section{L'application endoscopique $^c\theta_{\tilde{M}}^{\tilde{G},{\cal E}}$}

\bigskip

\subsection{D\'efinition d'une premi\`ere application endoscopique}
Le triplet $(G,\tilde{G},{\bf a})$ est quelconque et on fixe un espace de Levi $\tilde{M}$ de $\tilde{G}$. Soit ${\bf M}'=(M',{\cal M}',\tilde{\zeta})$ une donn\'ee endoscopique de $(M,\tilde{M},{\bf a})$ qui est elliptique et relevante. On d\'efinit une application lin\'eaire
$$^c\theta_{\tilde{M}}^{\tilde{G},{\cal E}}({\bf M}'):I(\tilde{G}(F),\omega)\otimes Mes(G(F))\to SI_{ac}({\bf M}')\otimes Mes(M'(F))$$
par l'\'egalit\'e
$$(1) \qquad ^c\theta_{\tilde{M}}^{\tilde{G},{\cal E}}({\bf M}',{\bf f})=\sum_{\tilde{s}\in \tilde{\zeta}Z(\hat{M})^{\Gamma_{F},\hat{\theta}}/Z(\hat{G})^{\Gamma_{F},\hat{\theta}}}i_{\tilde{M}'}(\tilde{G},\tilde{G}'(\tilde{s})){^cS\theta}_{{\bf M}'}^{{\bf G}'(\tilde{s})}({\bf f}^{{\bf G}'(\tilde{s})}).$$
Les m\^emes consid\'erations formelles qu'en 2.2 permettent de d\'efinir les termes du membre de droite. Ces termes sont bien d\'efinis car les distributions qui apparaissent v\'erifient  par r\'ecurrence la proposition 2.2,  c'est-\`a-dire sont stables, sauf dans un cas. C'est celui o\`u $(G,\tilde{G},{\bf a})$ est quasi-d\'eploy\'e et \`a torsion int\'erieure et o\`u ${\bf M}'={\bf M}$. On ne conna\^{\i}t pas alors les propri\'et\'es du  terme index\'e par $s=1$. On le remplace alors par ${^cS}\theta_{\tilde{M}}^{\tilde{G}}({\bf f})$. Par d\'efinition de ce dernier terme, on a alors l'\'egalit\'e tautologique
$$(2) \qquad ^c\theta_{\tilde{M}}^{\tilde{G},{\cal E}}({\bf M},{\bf f})=p_{\tilde{M}}^{st}\circ{^c\theta}_{\tilde{M}}^{\tilde{G}}({\bf f}).$$

Soit $b$ une fonction sur $\tilde{{\cal A}}_{\tilde{G}}$ \`a valeurs complexes. Pour tout $\tilde{s}$ intervenant dans la somme ci-dessus, on a une identification $\tilde{{\cal A}}_{\tilde{G}'(\tilde{s})}=\tilde{{\cal A}}_{\tilde{G}}$ et on peut consid\'erer $b$ comme une fonction sur le premier espace. D'autre part, la compos\'ee de la projection $\tilde{{\cal A}}_{\tilde{M}'}\to \tilde{{\cal A}}_{\tilde{G}'(\tilde{s})}$ et de l'identification ci-dessus est ind\'ependante de $\tilde{s}$. Elle co\"{\i}ncide en effet avec la compos\'ee de l'identification $\tilde{{\cal A}}_{\tilde{M}'}=\tilde{{\cal A}}_{\tilde{M}}$ et de la projection de cet espace sur $\tilde{{\cal A}}_{\tilde{G}}$. Notons $p_{\tilde{M}'}^{\tilde{G}}$ cette compos\'ee. D'apr\`es 2.1(1), on a  $({\bf f}(b\circ \tilde{H}_{\tilde{G}}))^{{\bf G}'(\tilde{s})}={\bf f}^{{\bf G}'(\tilde{s})}(b\circ\tilde{H}_{\tilde{G}'(\tilde{s})})$. On a d\'ej\`a dit que la propri\'et\'e 1.6(3) se propageait aux applications $^cS\theta_{\tilde{M}}^{\tilde{G}}$. Il r\'esulte de la d\'efinition (1) ci-dessus
qu'on a l'\'egalit\'e
$$(3) \qquad {^c\theta}_{\tilde{M}}^{\tilde{G},{\cal E}}({\bf M}',{\bf f}(b\circ \tilde{H}_{\tilde{G}}))={^c\theta}_{\tilde{M}}^{\tilde{G},{\cal E}}({\bf M}',{\bf f})(b\circ p_{\tilde{M}'}^{\tilde{G}}\circ\tilde{H}_{\tilde{M}'}).$$

{\bf Remarque.} L'identification $\tilde{{\cal A}}_{\tilde{M}'}=\tilde{{\cal A}}_{\tilde{M}}$ n'envoie pas n\'ecessairement $\tilde{{\cal A}}_{\tilde{M}',F}$ dans $\tilde{{\cal A}}_{\tilde{M},F}$. Donc $p_{\tilde{M}'}^{\tilde{G}}$ n'envoie pas n\'ecessairement $\tilde{{\cal A}}_{\tilde{M}',F}$ dans $\tilde{{\cal A}}_{\tilde{G},F}$. Le membre de gauche ci-dessus ne  d\'epend que de la restriction de $b$ \`a $\tilde{{\cal A}}_{\tilde{G},F}$. Il n'est pas  \'evident que ce soit le cas du membre de droite. Cela r\'esulte bien s\^ur de l'\'egalit\'e des deux membres. L'explication est que, par d\'efinition du transfert, les supports des transferts ${\bf f}^{{\bf G}'(\tilde{s})}$ sont  form\'es d'\'el\'ements $\delta\in \tilde{G}'(\tilde{s};F)$ tels que $\tilde{H}_{\tilde{G}'(\tilde{s})}(\delta)\in \tilde{{\cal A}}_{\tilde{G},F}$.

\bigskip

  Il r\'esulte de (1) et de ce que l'on a dit en 2.4 que

(4) pour tout ${\bf f}\in I(\tilde{G}(F),\omega)\otimes Mes(G(F))$, la fonction $^c\theta_{\tilde{M}}^{\tilde{G},{\cal E}}({\bf M}',{\bf f})$ est de Schwartz; pour tout $\tilde{\pi}\in D^{st}_{temp}({\bf M}')\otimes Mes(M'(F))^*$, la fonction $\tilde{\lambda}\mapsto SI^{{\bf M}'}(\tilde{\pi},\tilde{\lambda},{^c\theta}_{\tilde{M}}^{\tilde{G},{\cal E}}({\bf M}',{\bf f}))$ d\'efinie sur $i\tilde{{\cal A}}_{\tilde{M}'}^*/i\tilde{{\cal A}}_{\tilde{M}',F}^{\vee}$ se prolonge en une fonction rationnelle sur $\tilde{{\cal A}}_{\tilde{M}',{\mathbb C}}^*/i\tilde{{\cal A}}_{\tilde{M}',F}^{\vee}$ dont les p\^oles  sont de la forme  d\'ecrite en 1.8(2). 

\bigskip

\subsection{Action d'un groupe d'automorphismes}
On conserve les m\^emes donn\'ees. On a introduit en [I] 3.2 le groupe d'automorphismes $Aut(\tilde{M},{\bf M}')$. Il agit sur $SI_{ac}({\bf M}')\otimes Mes(M'(F))$.

\ass{Lemme}{L'application $^c\theta_{\tilde{M}}^{\tilde{G},{\cal E}}({\bf M}')$ prend ses valeurs dans le sous-espace d'invariants $(SI_{ac}({\bf M}')\otimes Mes(M'(F)))^{Aut(\tilde{M},{\bf M}')}$.}

La preuve est similaire \`a celle du lemme [II] 1.13. $\square$

\bigskip

\subsection{Commutation \`a l'induction}
On conserve les m\^emes donn\'ees. On consid\`ere les situations suivantes.

(a) Soit $R'$ un groupe de Levi de $M'$ qui est relevant. Modulo certains choix, cf. [I] 3.4, on construit un espace de Levi $\tilde{R}$ de $\tilde{M}$ et une donn\'ee endoscopique  ${\bf R}'$ de $\tilde{R}$ qui est elliptique et relevante. On a  l'homomorphisme "terme constant"
$$\begin{array}{ccc}SI_{ac}({\bf M}')\otimes Mes(M'(F))&\to&SI_{ac}({\bf R}')\otimes Mes(R'(F))\\ \boldsymbol{\varphi}&\mapsto&\boldsymbol{\varphi}_{{\bf R}'}.\\ \end{array}$$

(b) Soit $R'$ un groupe de Levi de $M'$ qui n'est pas relevant. La donn\'ee ${\bf R}'$ n'est plus d\'efinie et l'homomorphisme du (i) n'a plus de sens. Par contre, pour $\boldsymbol{\varphi}\in SI_{ac}({\bf M}')\otimes Mes(M'(F))$, la relation $\boldsymbol{\varphi}_{\tilde{R}'}=0$ conserve un sens. On fixe des donn\'ees auxilaires $M'_{1},...,\Delta_{1}$ pour ${\bf M}'$. On identifie $\boldsymbol{\varphi}$ \`a un \'el\'ement $\boldsymbol{\varphi}_{1}\in SI_{ac,\lambda_{1}}(\tilde{M}'_{1}(F))\otimes Mes(M'(F))$. La relation signifie que $(\boldsymbol{\varphi}_{1})_{\tilde{R}'_{1}}=0$, ce qui ne d\'epend pas du choix des donn\'ees auxiliaires.

\ass{Proposition}{Soit ${\bf f}\in I(\tilde{G}(F),\omega)\otimes Mes(G(F))$.

(i) Dans la situation (a), on a l'\'egalit\'e
$$(^c\theta_{\tilde{M}}^{\tilde{G},{\cal E}}({\bf M}',{\bf f}))_{{\bf R}'}=\sum_{\tilde{L}\in {\cal L}(\tilde{R})}d_{\tilde{R}}^{\tilde{G}}(\tilde{M},\tilde{L}){^c\theta}_{\tilde{R}}^{\tilde{L},{\cal E}}({\bf R}',{\bf f}_{\tilde{L},\omega})$$

(ii) Dans la situation (b), on a l'\'egalit\'e
$$(^c\theta_{\tilde{M}}^{\tilde{G},{\cal E}}({\bf M}',{\bf f}))_{\tilde{R}'}=0.$$}

Preuve. Dans la situation (a), soit $\boldsymbol{\delta}\in D_{g\acute{e}om}^{st}({\bf R}')\otimes Mes(R'(F))$. On doit prouver l'\'egalit\'e
$$SI^{{\bf M}'}(\boldsymbol{\delta}^{{\bf M}'},^c\theta_{\tilde{M}}^{\tilde{G},{\cal E}}({\bf M}',{\bf f}))=\sum_{\tilde{L}\in {\cal L}(\tilde{R})}d_{\tilde{R}}^{\tilde{G}}(\tilde{M},\tilde{L})SI^{{\bf R}'}(\boldsymbol{\delta},{^c\theta}_{\tilde{R}}^{\tilde{L},{\cal E}}({\bf R}',{\bf f}_{\tilde{L},\omega})).$$
Dans la situation (b), fixons des donn\'ees auxiliaires $M'_{1},...,\Delta_{1}$ pour ${\bf M}'$. Soit $\boldsymbol{\delta}\in D_{g\acute{e}om,\lambda_{1}}^{st}(\tilde{R}'_{1}(F))\otimes Mes(R'(F))$. L'induite $\boldsymbol{\delta}^{\tilde{M}'_{1}}$ appartient \`a $D_{g\acute{e}om,\lambda_{1}}^{st}(\tilde{M}'_{1}(F))\otimes Mes(M'(F))$. On l'identifie \`a un \'el\'ement de $D_{g\acute{e}om}^{st}({\bf M}')\otimes Mes(M'(F))$ que l'on note $\boldsymbol{\delta}^{{\bf M}'}$. 
On doit prouver l'\'egalit\'e
 $$SI^{{\bf M}'}(\boldsymbol{\delta}^{{\bf M}'},^c\theta_{\tilde{M}}^{\tilde{G},{\cal E}}({\bf M}',{\bf f}))=0.$$
 La preuve de ces assertions est la m\^eme que celle de (i) et (iii) de la proposition [II] 1.14, en utilisant la relation initiale fournie par le lemme 1.6. $\square$
 
\bigskip
\subsection{D\'efinition de $^c\theta_{\tilde{M}}^{\tilde{G},{\cal E}}$}
\ass{Proposition}{Il existe une unique application lin\'eaire
$$^c\theta_{\tilde{M}}^{\tilde{G},{\cal E}}:I(\tilde{G}(F),\omega)\otimes Mes(G(F))\to I_{ac}(\tilde{M}(F),\omega)\otimes Mes(M(F))$$
telle que, pour toute donn\'ee endoscopique ${\bf M}'$ de $(M,\tilde{M},{\bf a})$ qui est elliptique et relevante et pour tout ${\bf f}\in I(\tilde{G}(F),\omega)\otimes Mes(G(F))$, on ait l'\'egalit\'e
$$(^c\theta_{\tilde{M}}^{\tilde{G},{\cal E}}({\bf f}))^{{\bf M}'}={^c\theta}_{\tilde{M}}^{\tilde{G},{\cal E}}({\bf M}',{\bf f}).$$}

Preuve. Fixons un ensemble de repr\'esentants ${\cal E}(\tilde{M},{\bf a})$ des classes d'\'equivalence de donn\'ees endoscopiques elliptiques et relevantes de $(M,\tilde{M},{\bf a})$. Soit ${\bf f}\in I(\tilde{G}(F),\omega)\otimes Mes(G(F))$. Pour ${\bf M}'\in {\cal E}(\tilde{M},{\bf a})$, posons $\boldsymbol{\varphi}_{{\bf M}'}={^c\theta}_{\tilde{M}}^{\tilde{G},{\cal E}}({\bf M}',{\bf f})$. 
On a $\boldsymbol{\varphi}_{{\bf M}'}\in SI_{ac}({\bf M}')\otimes Mes(M'(F))$. Oublions pour un instant l'indice $ac$, c'est-\`a-dire supposons $\boldsymbol{\varphi}_{{\bf M}'}\in SI({\bf M}')\otimes Mes(M'(F))$. On a d\'ecrit dans la proposition [I] 4.11 l'image $I^{{\cal E}}(\tilde{M}(F),\omega)\otimes Mes(M(F))$ de l'application de transfert
$$I(\tilde{M}(F),\omega)\otimes Mes(M(F))\to \oplus_{{\bf M}'\in {\cal E}(\tilde{M},{\bf a})}SI({\bf M}')\otimes Mes(M'(F)).$$
Il s'agit de montrer que la collection $(\boldsymbol{\varphi}_{{\bf M}'})_{{\bf M}'\in {\cal E}(\tilde{M},{\bf a})}$ appartient \`a cette image. Elle doit pour cela v\'erifier les conditions (1), (2) et (3) de [I] 4.11. La condition d'invariance par automorphismes  est le lemme  3.2 Les deux autres conditions r\'esultent de la proposition 3.3. Cela conclut sous l'hypoth\`ese faite ci-dessus. Levons cette hypoth\`ese. Soit $b$ une fonction \`a support compact sur $\tilde{{\cal A}}_{\tilde{M}}$. Pour ${\bf M}'\in {\cal E}(\tilde{M},{\bf a})$, on a $\boldsymbol{\varphi}_{{\bf M}'}(b\circ\tilde{H}_{\tilde{M}'})\in SI({\bf M}')\otimes Mes(M'(F))$. On peut appliquer le m\^eme raisonnement \`a la collection $(\boldsymbol{\varphi}_{{\bf M}'}(b\circ\tilde{H}_{\tilde{M}'}))_{{\bf M}'\in {\cal E}(\tilde{M},{\bf a})}$. Celle-ci est donc le transfert d'un unique \'el\'ement $\boldsymbol{\varphi}[b]\in I(\tilde{M}(F),\omega)\otimes Mes(M(F))$. En choisissant une suite de fonctions $(b_{i})_{i\in {\mathbb N}}$ localement finie et  telle que $\sum_{i\in {\mathbb N}}b_{i}$ soit constante de valeur $1$, on v\'erifie ais\'ement que la suite de fonctions $(\boldsymbol{\varphi}[b_{i}])_{i\in {\mathbb N}}$ est localement finie. La  s\'erie $\sum_{i\in {\mathbb N}}\boldsymbol{\varphi}[b_{i}]$ converge donc et est son transfert est bien la collection $(\boldsymbol{\varphi}_{{\bf M}'})_{{\bf M}'\in {\cal E}(\tilde{M},{\bf a})}$. Il faut toutefois v\'erifier que la somme de la s\'erie appartient \`a $I_{ac}(\tilde{M}(F),\omega)\otimes Mes(M(F))$, c'est-\`a-dire qu'elle provient d'une fonction sur $\tilde{M}(F)$ qui est biinvariante par un sous-groupe ouvert compact. Pour chaque ${\bf M}'\in {\cal E}(\tilde{M},{\bf a})$, fixons des donn\'ees auxiliaires $M'_{1},...,\Delta_{1}$.  Identifions $\boldsymbol{\varphi}_{{\bf M}'}$ \`a un \'el\'ement de $SI_{ac,\lambda_{1}}(\tilde{M}'_{1}(F))\otimes Mes(M'(F))$. Par d\'efinition de cet espace, on peut fixer un sous-groupe ouvert compact $K_{M'_{1}}$ de $M'_{1}(F)$ de sorte que $\boldsymbol{\varphi}_{{\bf M}'}$ provienne d'une fonction sur $\tilde{M}'_{1}(F)$ biinvariante par $K_{M'_{1}}$. Il en r\'esulte que, pour tout $i\in {\mathbb N}$, $\boldsymbol{\varphi}_{{\bf M}'}(b_{i}\circ H_{\tilde{M}'})$ provient d'une fonction sur $\tilde{M}'_{1}(F)$ biinvariante par $K_{M'_{1}}$. D'apr\`es [M], il existe un sous-groupe ouvert compact $K_{M}$ de $M(F)$, qui ne d\'epend que des $K_{M'_{1}}$ (donc qui ne d\'epend pas de $i$), de sorte que $\boldsymbol{\varphi}[b_{i}]$ provienne d'une fonction sur $\tilde{M}(F)$ qui est biinvariante par $K_{M}$. Alors la somme de la s\'erie $\sum_{i\in {\mathbb N}}\boldsymbol{\varphi}[b_{i}]$ provient elle-aussi d'une fonction biinvariante par $K_{M}$. Cela ach\`eve la preuve. $\square$

\bigskip
\subsection{Commutation \`a l'induction}
\ass{Lemme}{Soient $\tilde{R}$ un espace de Levi de $\tilde{M}$ et ${\bf f}\in I(\tilde{G}(F),\omega)\otimes Mes(G(F))$. On a l'\'egalit\'e
$$^c\theta_{\tilde{M}}^{\tilde{G},{\cal E}}({\bf f})_{\tilde{R},\omega}=\sum_{\tilde{L}\in {\cal L}(\tilde{M})}d_{\tilde{R}}^{\tilde{G}}(\tilde{M},\tilde{L}){^c\theta}_{\tilde{R}}^{\tilde{L},{\cal E}}({\bf f}_{\tilde{L},\omega}).$$}

Preuve. Il suffit de voir que, pour toute donn\'ee endoscopique ${\bf R}'$ de $(R,\tilde{R},{\bf a})$ qui est elliptique et relevante, les deux membres ont m\^eme transfert \`a ${\bf R}'$. Pour une telle donn\'ee ${\bf R}'$, on peut  trouver une donn\'ee endoscopique ${\bf M}'$ de $(M,\tilde{M},{\bf a})$ qui est elliptique et relevante, de sorte que ${\bf R}'$ soit une "donn\'ee de Levi" de ${\bf M}'$. L'\'egalit\'e cherch\'ee r\'esulte alors de la proposition  3.3 et de la relation
$$(\boldsymbol{\varphi}_{\tilde{R},\omega})^{{\bf R}'}=(\boldsymbol{\varphi}^{{\bf M}'})_{{\bf R}'}$$
pour tout $\boldsymbol{\varphi}\in I(\tilde{M}(F),\omega)\otimes Mes(M(F))$. $\square$

\bigskip

\subsection{$^c\theta_{\tilde{M}}^{\tilde{G},{\cal E}}({\bf f})$ est de Schwartz}

\ass{Lemme}{Soit ${\bf f}\in I(\tilde{G}(F),\omega)\otimes Mes(G(F))$. Alors $^c\theta_{\tilde{M}}^{\tilde{G},{\cal E}}({\bf f})$ est une fonction de Schwartz. Pour tout $\tilde{\pi}\in D_{temp}(\tilde{M}(F),\omega)\otimes Mes(M(F))^*$, la fonction $\tilde{\lambda}\mapsto I^{\tilde{G}}(\tilde{\pi},\tilde{\lambda},{^c\theta}_{\tilde{M}}^{\tilde{G},{\cal E}}({\bf f}))$ d\'efinie pour $\tilde{\lambda}\in i\tilde{{\cal A}}_{\tilde{M}}^*/i\tilde{{\cal A}}_{\tilde{M},F}^{\vee}$ se prolonge en une fonction rationnelle sur $\tilde{{\cal A}}_{\tilde{M},{\mathbb C}}^*/i\tilde{{\cal A}}_{\tilde{M},F}^{\vee}$ dont les p\^oles sont de la forme d\'ecrite en 1.8(2). }
 
 Preuve.  Soit $\tilde{\pi}\in D_{temp}(\tilde{M}(F),\omega)\otimes Mes(M(F))^*$. On doit montrer que la fonction $X\mapsto I^{\tilde{M}}(\tilde{\pi},X,{^c\theta}_{\tilde{M}}^{\tilde{G},{\cal E}}({\bf f}))$ d\'efinie sur $\tilde{{\cal A}}_{\tilde{M},F}$ est \`a d\'ecroissance rapide. On doit ensuite \'etudier la fonction $\tilde{\lambda}\mapsto I^{\tilde{G}}(\tilde{\pi},\tilde{\lambda},{^c\theta}_{\tilde{M}}^{\tilde{G},{\cal E}}({\bf f}))$. Par lin\'earit\'e, on peut supposer soit que $\tilde{\pi}$ est elliptique, soit qu'il existe un espace de Levi propre $\tilde{R}$ de $\tilde{M}$ et une $\omega$-repr\'esentation elliptique $\tilde{\sigma}$ de $\tilde{R}(F)$ de sorte que $\tilde{\pi}=Ind_{\tilde{R}}^{\tilde{M}}(\tilde{\sigma})$.

 Supposons d'abord $\tilde{\pi}$ elliptique. 
Fixons un ensemble de repr\'esentants ${\cal E}(\tilde{M},{\bf a})$ des classes d'\'equivalence de donn\'ees endoscopiques elliptiques et relevantes de $(M,\tilde{M},{\bf a})$. 
Rappelons que le transfert d\'efinit un isomorphisme
$$\oplus_{{\bf M}'\in {\cal E}(\tilde{M},{\bf a})}(D_{ell}^{st}({\bf M}')\otimes Mes(M'(F))^*)^{Aut({\bf M}')}\to D_{ell}(\tilde{M}(F),\omega)\otimes Mes(M(F))^*.$$
 On \'ecrit conform\'ement 
$$(1) \qquad \tilde{\pi}=\sum_{{\bf M}'\in {\cal E}(\tilde{M},{\bf a})}transfert(\tilde{\pi}_{{\bf M}'}),$$
o\`u $\tilde{\pi}_{{\bf M}'}\in (D_{ell}^{st}({\bf M}')\otimes Mes(M'(F))^*)^{Aut({\bf M}')}$. Pour $\tilde{\lambda}\in \tilde{{\cal A}}_{\tilde{M},{\mathbb C}}^*$, il est clair que
$$\tilde{\pi}_{\tilde{\lambda}}=\sum_{{\bf M}'\in {\cal E}(\tilde{M},{\bf a})}transfert(\tilde{\pi}_{{\bf M}',\tilde{\lambda}}).$$
Le membre de gauche ne d\'epend que de $\tilde{\lambda}$ modulo $i\tilde{{\cal A}}_{\tilde{M},F}^{\vee}$. Il n'est pas clair a priori que chaque  terme $\tilde{\pi}_{{\bf M}',\tilde{\lambda}}$  v\'erifie la m\^eme propri\'et\'e, car les ensembles $\tilde{{\cal A}}_{\tilde{M},F}^{\vee}$ et $\tilde{{\cal A}}_{\tilde{M}',F}^{\vee}$ peuvent \^etre distincts. Mais l'unicit\'e de la d\'ecomposition ci-dessus assure qu'il en est bien ainsi: $\tilde{\pi}_{{\bf M}',\tilde{\lambda}}$  ne d\'epend que de $\tilde{\lambda}$ modulo $i\tilde{{\cal A}}_{\tilde{M},F}^{\vee}$. 

  Pour $X\in \tilde{{\cal A}}_{\tilde{M},F}$, on a l'\'egalit\'e 
$$I^{\tilde{M}}(\tilde{\pi},X,{^c\theta}_{\tilde{M}}^{\tilde{G},{\cal E}}({\bf f}))=I^{\tilde{M}}(\tilde{\pi},({\bf 1}_{X}\circ \tilde{H}_{\tilde{M}}){^c\theta}_{\tilde{M}}^{\tilde{G},{\cal E}}({\bf f})),$$
o\`u ${\bf 1}_{X}$ est la fonction caract\'eristique de $X$ dans $\tilde{{\cal A}}_{\tilde{M}}$. Appliquons (1). Pour tout ${\bf M}'\in {\cal E}(\tilde{M},{\bf a})$, le transfert \`a ${\bf M}'$ de la fonction ${^c\theta}_{\tilde{M}}^{\tilde{G},{\cal E}}({\bf f})$ est ${^c\theta}_{\tilde{M}}^{\tilde{G},{\cal E}}({\bf M}',{\bf f})$. Le transfert de la fonction $({\bf 1}_{X}\circ \tilde{H}_{\tilde{M}}){^c\theta}_{\tilde{M}}^{\tilde{G},{\cal E}}({\bf f})$ est $({\bf 1}_{X}\circ \tilde{H}_{\tilde{M}'}){^c\theta}_{\tilde{M}}^{\tilde{G},{\cal E}}({\bf M}',{\bf f})$. D'o\`u l'\'egalit\'e
$$I^{\tilde{M}}(\tilde{\pi},X,{^c\theta}_{\tilde{M}}^{\tilde{G},{\cal E}}({\bf f}))=\sum_{{\bf M}'\in {\cal E}(\tilde{M},{\bf a})}SI^{{\bf M}'}(\tilde{\pi}_{{\bf M}'},({\bf 1}_{X}\circ \tilde{H}_{\tilde{M}'}){^c\theta}_{\tilde{M}}^{\tilde{G},{\cal E}}({\bf M}',{\bf f}))$$
$$=\sum_{{\bf M}'\in {\cal E}(\tilde{M},{\bf a})}SI^{{\bf M}'}(\tilde{\pi}_{{\bf M}'},X, {^c\theta}_{\tilde{M}}^{\tilde{G},{\cal E}}({\bf M}',{\bf f})).$$
D'apr\`es 3.1(4), les fonctions ${^c\theta}_{\tilde{M}}^{\tilde{G},{\cal E}}({\bf M}',{\bf f})$ sont de Schwartz. L'expression ci-dessus est donc \`a d\'ecroissance rapide en $X$. On peut donc d\'efinir la fonction $\tilde{\lambda}\mapsto I^{\tilde{G}}(\tilde{\pi},\tilde{\lambda},{^c\theta}_{\tilde{M}}^{\tilde{G},{\cal E}}({\bf f}))$ sur $i\tilde{{\cal A}}_{\tilde{M}}^*/i\tilde{{\cal A}}_{\tilde{M},F}^{\vee}$ par
$$I^{\tilde{G}}(\tilde{\pi},\tilde{\lambda},{^c\theta}_{\tilde{M}}^{\tilde{G},{\cal E}}({\bf f}))=\sum_{X\in \tilde{{\cal A}}_{\tilde{M},F}}e^{<\tilde{\lambda},X>}I^{\tilde{M}}(\tilde{\pi},X,{^c\theta}_{\tilde{M}}^{\tilde{G},{\cal E}}({\bf f})).$$
Pour ${\bf M}'\in {\cal E}(\tilde{M},{\bf a})$ on voit appara\^{\i}tre la somme
$$\sum_{X\in \tilde{{\cal A}}_{\tilde{M},F}}e^{<\tilde{\lambda},X>}SI^{{\bf M}'}(\tilde{\pi}_{{\bf M}'},X, {^c\theta}_{\tilde{M}}^{\tilde{G},{\cal E}}({\bf M}',{\bf f})).$$
 Pour $X\not\in \tilde{{\cal A}}_{\tilde{M}',F}$, les termes $ SI^{{\bf M}'}(\tilde{\pi}_{{\bf M}'}, X, {^c\theta}_{\tilde{M}}^{\tilde{G},{\cal E}}({\bf M}',{\bf f}))$ sont nuls par d\'efinition. Parce que $\tilde{\pi}_{{\bf M}',\tilde{\lambda}}$  ne d\'epend que de $\tilde{\lambda}$ modulo $i\tilde{{\cal A}}_{\tilde{M},F}^{\vee}$, les m\^emes termes sont nuls pour $X\in \tilde{{\cal A}}_{\tilde{M}',F}- (\tilde{{\cal A}}_{\tilde{M},F}\cap \tilde{{\cal A}}_{\tilde{M}',F})$. 
On peut donc remplacer l'ensemble de sommation ci-dessus par $  \tilde{{\cal A}}_{\tilde{M}',F}$. La somme devient $SI^{{\bf M}'}(\tilde{\pi}_{{\bf M}'}, \tilde{\lambda}, {^c\theta}_{\tilde{M}}^{\tilde{G},{\cal E}}({\bf f}) )$.  On obtient
$$(1)\qquad I^{\tilde{G}}(\tilde{\pi},\tilde{\lambda},{^c\theta}_{\tilde{M}}^{\tilde{G},{\cal E}}({\bf f}))=\sum_{{\bf M}'\in {\cal E}(\tilde{M},{\bf a})} SI^{{\bf M}'}(\tilde{\pi}_{{\bf M}'}, \tilde{\lambda}, {^c\theta}_{\tilde{M}}^{\tilde{G},{\cal E}}({\bf f}) ).$$
 Gr\^ace \`a 3.1(4), cette fonction se prolonge en une fonction rationnelle sur $\tilde{{\cal A}}_{\tilde{M},{\mathbb C}}/i\tilde{{\cal A}}_{\tilde{M},F}^*$ avec des hyperplans polaires de la forme habituelle. Cela d\'emontre les propri\'et\'es voulues pour une $\omega$-repr\'esentation elliptique. 
 
 Fixons maintenant un espace de Levi propre $\tilde{R}$ de $\tilde{M}$ et une $\omega$-repr\'esentation elliptique $\tilde{\sigma}$ de $\tilde{R}(F)$. Posons $\tilde{\pi}=Ind_{\tilde{R}}^{\tilde{M}}(\tilde{\sigma})$.  Soit $X\in \tilde{{\cal A}}_{\tilde{M},F}$. D'apr\`es le lemme 3.5, on a l'\'egalit\'e
 $$(2) \qquad I^{\tilde{G}}(\tilde{\pi},({\bf 1}_{X}\circ\tilde{H}_{\tilde{M}}){^c\theta}_{\tilde{M}}^{\tilde{G},{\cal E}}({\bf f}))=\sum_{\tilde{L}\in {\cal L}(\tilde{R})}d_{\tilde{R}}^{\tilde{G}}(\tilde{M},\tilde{L}) I^{\tilde{L}}(\tilde{\sigma},({\bf 1}_{X}\circ \tilde{H}_{\tilde{M}}){^c\theta}_{\tilde{R}}^{\tilde{L},{\cal E}}({\bf f}_{\tilde{L},\omega})).$$
Fixons  
 un espace de Levi $\tilde{L}\in {\cal L}(\tilde{R})$ tel que $d_{\tilde{R}}^{\tilde{G}}(\tilde{M},\tilde{L})\not=0$.  
 On a une projection $Y\mapsto Y_{\tilde{M}}$ de $\tilde{{\cal A}}_{\tilde{R},F}$ sur $\tilde{{\cal A}}_{\tilde{M},F}$. On a l'\'egalit\'e
 $$({\bf 1}_{X}\circ \tilde{H}_{\tilde{M}}){^c\theta}_{\tilde{R}}^{\tilde{L},{\cal E}}({\bf f}_{\tilde{L},\omega})=\sum_{Y\in \tilde{{\cal A}}_{\tilde{R},F}; Y_{\tilde{M}}=X}({\bf 1}_{Y}\circ \tilde{H}_{\tilde{R}}){^c\theta}_{\tilde{R}}^{\tilde{L},{\cal E}}({\bf f}_{\tilde{L},\omega}).$$
 La projection dans $\tilde{{\cal A}}_{\tilde{L}}$ du support de ${^c\theta}_{\tilde{R}}^{\tilde{L},{\cal E}}({\bf f}_{\tilde{L},\omega})$ est compacte. L'hypoth\`ese $d_{\tilde{R}}^{\tilde{G}}(\tilde{M},\tilde{L})\not=0$ entra\^{\i}ne que l'ensemble des $Y\in \tilde{{\cal A}}_{\tilde{R},F}$ dont la projection dans $\tilde{{\cal A}}_{\tilde{L}}$ appartient \`a ce compact et dont la projection dans $\tilde{{\cal A}}_{\tilde{M}}$ est $X$ est fini. La somme ci-dessus n'a donc qu'un nombre fini de termes non nuls. Il en r\'esulte l'\'egalit\'e
  $$(3) \qquad I^{\tilde{L}}(\tilde{\sigma},({\bf 1}_{X}\circ \tilde{H}_{\tilde{M}}){^c\theta}_{\tilde{R}}^{\tilde{L},{\cal E}}({\bf f}_{\tilde{L},\omega}))=\sum_{Y\in \tilde{{\cal A}}_{\tilde{R},F}; Y_{\tilde{M}}=X}I^{\tilde{L}}(\tilde{\sigma}, Y,{^c\theta}_{\tilde{R}}^{\tilde{L},{\cal E}}({\bf f}_{\tilde{L},\omega})).$$
  Puisque $\tilde{R}\not=\tilde{M}$, on sait par r\'ecurrence que la fonction $Y\mapsto I^{\tilde{L}}(\tilde{\sigma}, Y,{^c\theta}_{\tilde{R}}^{\tilde{L},{\cal E}}({\bf f}_{\tilde{L},\omega}))$ est \`a d\'ecroissance rapide. Donc l'expression ci-dessus est \`a d\'ecroissance rapide en $X$. D'apr\`es (2), il en est de m\^eme de $ I^{\tilde{G}}(\tilde{\pi},({\bf 1}_{X}\circ\tilde{H}_{\tilde{M}}){^c\theta}_{\tilde{M}}^{\tilde{G},{\cal E}}({\bf f}))$. 
  En multipliant  (3) par $e^{<\tilde{\lambda},X>}$ et en sommant en $X$, on obtient
 $$\sum_{X\in \tilde{{\cal A}}_{\tilde{M},F}}e^{<\tilde{\lambda},X>} I^{\tilde{L}}(\tilde{\sigma},({\bf 1}_{X}\circ \tilde{H}_{\tilde{M}}){^c\theta}_{\tilde{R}}^{\tilde{L},{\cal E}}({\bf f}_{\tilde{L},\omega}))=I^{\tilde{L}}(\tilde{\sigma}, \tilde{\lambda},{^c\theta}_{\tilde{R}}^{\tilde{L},{\cal E}}({\bf f}_{\tilde{L},\omega})).$$
 D'o\`u, gr\^ace \`a (2),
$$ I^{\tilde{G}}(\tilde{\pi},\tilde{\lambda},{^c\theta}_{\tilde{M}}^{\tilde{G},{\cal E}}({\bf f}))=\sum_{\tilde{L}\in {\cal L}(\tilde{R})}d_{\tilde{R}}^{\tilde{G}}(\tilde{M},\tilde{L}) I^{\tilde{L}}(\tilde{\sigma}, \tilde{\lambda},{^c\theta}_{\tilde{R}}^{\tilde{L},{\cal E}}({\bf f}_{\tilde{L},\omega})).$$
Les propri\'et\'es voulues du terme de gauche r\'esultent des m\^emes propri\'et\'es du membre de droite, qui sont connues par r\'ecurrence. $\square$
 
\bigskip

\subsection{Une propri\'et\'e d'annulation}
Soient ${\bf f}\in I(\tilde{G}(F),\omega)\otimes Mes(G(F))$ et $\tilde{\pi}\in D_{temp}(\tilde{M}(F),\omega)\otimes Mes(M(F))^*$. La fonction $\tilde{\lambda}\mapsto I^{\tilde{M}}(\tilde{\pi},\tilde{\lambda},{^c\theta}_{\tilde{M}}^{\tilde{G},{\cal E}}({\bf f}))$ a les m\^emes propri\'et\'es qu'en 1.8 et on reprend pour cette fonction les constructions et notations de ce paragraphe.

\ass{Proposition}{Supposons $\tilde{M}\not=\tilde{G}$ et supposons que $\tilde{\pi}$ est elliptique. Pour tout $\tilde{S}\in {\cal P}(\tilde{M})$, fixons un point $\nu_{\tilde{S}}$ comme en 1.3. Supposons ce point  "assez positif" pour $\tilde{S}$, cette notion d\'ependant de la repr\'esentation $\tilde{\pi}$. Alors on a l'\'egalit\'e
$$\sum_{\tilde{S}\in {\cal P}(\tilde{M})}\omega_{\tilde{S}}(X)I^{\tilde{M}}(\tilde{\pi},\nu_{\tilde{S}},X,{^c\theta}_{\tilde{M}}^{\tilde{G},{\cal E}}({\bf f}))=0$$
pour tout $X\in \tilde{{\cal A}}_{\tilde{M},F}$.}

Preuve. On reprend les notations de la premi\`ere partie de la preuve du lemme 3.6. En vertu des \'egalit\'es 3.6(1) et 3.1(1), on peut fixer ${\bf M}'=(M',{\cal M}',\tilde{\zeta})$ et $\tilde{s}\in \tilde{\zeta}Z(\hat{M})^{\Gamma_{F},\hat{\theta}}/Z(\hat{G})^{\Gamma_{F},\hat{\theta}}$ tel que ${\bf G}'(\tilde{s})$ soit elliptique et prouver l'\'egalit\'e
$$(1) \qquad \sum_{\tilde{S}\in {\cal P}(\tilde{M})}\omega_{\tilde{S}}(X)S^{{\bf M}'}(\tilde{\pi}_{{\bf M}'},\nu_{\tilde{S}},X,{^cS\theta}_{{\bf M}'}^{{\bf G}'(\tilde{s})}({\bf f}^{{\bf G}'(\tilde{s})}))=0.$$
 Comme on l'a dit dans la preuve du lemme 3.6,  les termes intervenant  sont nuls si $X\not\in \tilde{{\cal A}}_{\tilde{M},F}\cap \tilde{{\cal A}}_{\tilde{M}',F}$. On peut donc supposer que $X$ appartient \`a cette intersection. La proposition 2.4 nous dit alors que l'on a une \'egalit\'e
$$\sum_{\tilde{S}'\in {\cal P}^{\tilde{G}'(\tilde{s})}(\tilde{M}')}\omega_{\tilde{S}'}(X)S^{{\bf M}'}(\tilde{\pi}_{{\bf M}'},\nu_{\tilde{S}'},X,{^cS\theta}_{{\bf M}'}^{{\bf G}'(\tilde{s})}({\bf f}^{{\bf G}'(\tilde{s})}))=0.$$
 Le m\^eme calcul que dans la preuve de la proposition 2.4 montre que cette \'egalit\'e est \'equivalente \`a (1). $\square$
 
 \bigskip
 
 \subsection{Egalit\'e de deux applications lin\'eaires}
 
 \ass{Proposition (\`a prouver)}{Soient ${\bf f}\in I(\tilde{G}(F),\omega)\otimes Mes(G(F))$. On a l'\'egalit\'e
 $$^c\theta_{\tilde{M}}^{\tilde{G},{\cal E}}({\bf f})={^c\theta}_{\tilde{M}}^{\tilde{G}}({\bf f}).$$}
 
 Cette proposition sera prouv\'ee en 4.3 dans le cas o\`u $(G,\tilde{G},{\bf a})$ est quasi-d\'eploy\'e et \`a torsion int\'erieure. Elle sera prouv\'ee en g\'en\'eral dans le m\^eme paragraphe, sous une hypoth\`ese qui ne sera v\'erifi\'ee que plus tard.
 
 \bigskip
 
 \subsection{Variante des int\'egrales orbitales pond\'er\'ees elliptiques}
  Soit ${\bf M}'=(M',{\cal M}',\tilde{\zeta})$ une donn\'ee endoscopique de $(M,\tilde{M},{\bf a})$ qui est elliptique et relevante. Soit $\boldsymbol{\delta}\in D_{g\acute{e}om,\tilde{G}-reg}^{st}({\bf M}')\otimes Mes(M'(F))^*$ et soit ${\bf f}\in I(\tilde{G}(F),\omega)\otimes Mes(G(F))$. On pose
 $$^cI_{\tilde{M}}^{\tilde{G},{\cal E}}({\bf M}',\boldsymbol{\delta},{\bf f})=\sum_{\tilde{s}\in \tilde{\zeta}Z(\hat{M})^{\Gamma_{F},\hat{\theta}}/Z(\hat{G})^{\Gamma_{F},\hat{\theta}}}i_{\tilde{M}'}(\tilde{G},\tilde{G}'(\tilde{s})){^cS}_{{\bf M}'}^{{\bf G}'(\tilde{s})}(\boldsymbol{\delta},{\bf f}^{{\bf G}'(\tilde{s})}).$$ 
Comme toujours, il y a un cas particulier. Si $(G,\tilde{G},{\bf a})$ est quasi-d\'eploy\'e et \`a torsion int\'erieure et si ${\bf M}'={\bf M}$, on doit remplacer le terme $^cS_{{\bf M}}^{{\bf G}'(1)}(\boldsymbol{\delta},{\bf f}^{{\bf G}'(1)})$ par $^cS_{\tilde{M}}^{\tilde{G}}(\boldsymbol{\delta},{\bf f})$. Dans ce cas, on a tautologiquement $^cI_{\tilde{M}}^{\tilde{G},{\cal E}}({\bf M},\boldsymbol{\delta},{\bf f})={^cI}_{\tilde{M}}^{\tilde{G}}(\boldsymbol{\delta},{\bf f})$.

Le terme ainsi d\'efini jouit des m\^emes propri\'et\'es que le terme $I_{\tilde{M}}^{\tilde{G},{\cal E}}({\bf M}',\boldsymbol{\delta},{\bf f})$ d\'efini en [II] 1.12. En particulier les propri\'et\'es d\'emontr\'ees en [II] 1.13,  1.14 et 1.15 valent aussi pour nos nouveaux objets. On peut alors poser la m\^eme d\'efinition qu'en [II] 1.15. C'est-\`a-dire, fixons un ensemble ${\cal E}(\tilde{M},{\bf a})$ de repr\'esentants des classes d'\'equivalence de donn\'ees endoscopiques de $(M,\tilde{M},{\bf a})$ qui sont elliptiques et relevantes. Soit $\boldsymbol{\gamma}\in D_{g\acute{e}om,\tilde{G}-reg}(\tilde{M}(F),\omega)\otimes Mes(M(F))^*$. On peut l'\'ecrire sous la forme
$$(1) \qquad \boldsymbol{\gamma}=\sum_{{\bf M}'\in {\cal E}(\tilde{M},{\bf a})}transfert(\boldsymbol{\delta}_{{\bf M}'}),$$
avec des $\boldsymbol{\delta}_{{\bf M}'}\in  D_{g\acute{e}om,\tilde{G}-reg}^{st}({\bf M}')\otimes Mes(M'(F))^*$. On pose
$$^cI_{\tilde{M}}^{\tilde{G},{\cal E}}(\boldsymbol{\gamma},{\bf f})=\sum_{{\bf M}'\in {\cal E}(\tilde{M},{\bf a})}{^cI}_{\tilde{M}}^{\tilde{G},{\cal E}}({\bf M}',\boldsymbol{\delta}_{{\bf M}'},{\bf f}).$$
La d\'ecomposition (1) n'est pas unique mais le terme ci-dessus ne d\'epend pas de la d\'ecomposition choisie. 

La propri\'et\'e de compacit\'e 1.9(2) se propage \`a notre distribution $^cI_{\tilde{M}}^{\tilde{G},{\cal E}}(\boldsymbol{\gamma},.)$. 

\ass{Proposition (\`a prouver)}{Pour tout $\boldsymbol{\gamma}\in D_{g\acute{e}om,\tilde{G}-reg}(\tilde{M}(F),\omega)\otimes Mes(M(F))^*$ et tout ${\bf f}\in I(\tilde{G}(F),\omega)\otimes Mes(G(F))$, on a l'\'egalit\'e
$$^cI_{\tilde{M}}^{\tilde{G},{\cal E}}(\boldsymbol{\gamma},{\bf f})={^cI}_{\tilde{M}}^{\tilde{G}}(\boldsymbol{\gamma},{\bf f}).$$}

 Cette proposition sera prouv\'ee en 4.3 dans le cas o\`u $(G,\tilde{G},{\bf a})$ est quasi-d\'eploy\'e et \`a torsion int\'erieure. Elle sera prouv\'ee en g\'en\'eral dans le m\^eme paragraphe, sous une hypoth\`ese qui ne sera v\'erifi\'ee que plus tard.

{\bf Remarque.} On s'est limit\'e ici \`a des distributions $\boldsymbol{\gamma}$ \`a support fortement r\'egulier dans $\tilde{G}$. Cela parce qu'une d\'efinition correcte pour les \'el\'ements ne v\'erifiant pas cette propri\'et\'e n\'ecessite l'introduction du syst\`eme de fonctions $B^{\tilde{G}}$ de [II] 1.11. Cela compliquerait grandement les choses alors que les distributions \`a support fortement r\'egulier suffiront aux applications. Mais, dans le cas o\`u $(G,\tilde{G},{\bf a})$ est quasi-d\'eploy\'e et \`a torsion int\'erieure, le recours au syst\`eme de fonctions $B^{\tilde{G}}$ se trivialise. Dans ce cas, on peut poser les m\^emes d\'efinitions que ci-dessus et \'enoncer la m\^eme proposition pour une distribution $\boldsymbol{\gamma}$ \`a support quelconque.

\bigskip

\section{Les preuves  et l'application $\epsilon_{\tilde{M}}$}

\bigskip

\subsection{Lien entre les int\'egrales orbitales pond\'er\'ees  stables ou endoscopiques et leurs variantes}
\ass{Proposition}{(i) Soient $\boldsymbol{\gamma}\in D_{g\acute{e}om,\tilde{G}-reg}(\tilde{M}(F),\omega)\otimes Mes(M(F))^*$ et ${\bf f}\in I(\tilde{G}(F),\omega)\otimes Mes(G(F))$. On a l'\'egalit\'e
$$^cI_{\tilde{M}}^{\tilde{G},{\cal E}}(\boldsymbol{\gamma},{\bf f})=\sum_{\tilde{L}\in {\cal L}(\tilde{M})}I_{\tilde{M}}^{\tilde{L},{\cal E}}(\boldsymbol{\gamma},{^c\theta}_{\tilde{L}}^{\tilde{G},{\cal E}}({\bf f})).$$ 

 (ii) Supposons $(G,\tilde{G},{\bf a})$ quasi-d\'eploy\'e et \`a torsion int\'erieure. Soient $\boldsymbol{\delta}\in D_{g\acute{e}om}^{st}(\tilde{M}(F))\otimes Mes(M(F))^*$ et ${\bf f}\in I(\tilde{G}(F))\otimes Mes(G(F))$. On a l'\'egalit\'e
 $$^cS_{\tilde{M}}^{\tilde{G}}(\boldsymbol{\delta},{\bf f})=\sum_{\tilde{L}\in {\cal L}(\tilde{M})}S_{\tilde{M}}^{\tilde{L}}(\boldsymbol{\delta},{^cS\theta}_{\tilde{L}}^{\tilde{G}}({\bf f})).$$}
 
 {\bf Remarque.} L'expression (ii) a un sens puisqu'on sait que les distributions $S_{\tilde{M}}^{\tilde{L}}(\boldsymbol{\delta},.)$ sont stables, cf. [II] th\'eor\`eme 1.10. 
 
 \bigskip

Preuve de (i). Par lin\'earit\'e, on peut fixer une donn\'ee endoscopique ${\bf M}'=(M',{\cal M}',\tilde{\zeta})$ de $(M,\tilde{M},{\bf a})$, elliptique et relevante, et un \'el\'ement $\boldsymbol{\delta}\in D_{g\acute{e}om,\tilde{G}-reg}^{st}({\bf M}')\otimes Mes(M'(F))$ et supposer  que $\boldsymbol{\gamma}=transfert(\boldsymbol{\delta})$. D'apr\`es les d\'efinitions, l'\'egalit\'e \`a prouver devient
$$(1) \qquad ^cI_{\tilde{M}}^{\tilde{G},{\cal E}}({\bf M}',\boldsymbol{\delta},{\bf f})=\sum_{\tilde{L}\in {\cal L}(\tilde{M})}I_{\tilde{M}}^{\tilde{L},{\cal E}}({\bf M}',\boldsymbol{\delta},{^c\theta}_{\tilde{L}}^{\tilde{G},{\cal E}}({\bf f})).$$ 

Consid\'erons d'abord le cas o\`u $(G,\tilde{G},{\bf a})$ est quasi-d\'eploy\'e et \`a torsion int\'erieure et o\`u ${\bf M}'={\bf M}$. Comme on l'a dit, on a tautologiquement l'\'egalit\'e
$$^cI_{\tilde{M}}^{\tilde{G},{\cal E}}({\bf M},\boldsymbol{\delta},{\bf f})={^cI}_{\tilde{M}}^{\tilde{G}}(\boldsymbol{\delta},{\bf f}).$$
De fa\c{c}on \'egalement tautologique, les distributions
$I_{\tilde{M}}^{\tilde{L},{\cal E}}({\bf M},\boldsymbol{\delta},.)$ et $I_{\tilde{M}}^{\tilde{L}}(\boldsymbol{\delta},.)$ sont \'egales. L'\'egalit\'e \`a prouver prend la forme
$$(2) \qquad ^cI_{\tilde{M}}^{\tilde{G}}(\boldsymbol{\delta},{\bf f})=\sum_{\tilde{L}\in {\cal L}(\tilde{M})}I_{\tilde{M}}^{\tilde{L}}(\boldsymbol{\delta},{^c\theta}_{\tilde{L}}^{\tilde{G},{\cal E}}({\bf f})).$$
Si $\tilde{L}\not=\tilde{M}$, nos habituelles hypoth\`eses de r\'ecurrence nous autorisent \`a utiliser la proposition 3.8: ${^c\theta}_{\tilde{L}}^{\tilde{G},{\cal E}}({\bf f})={^c\theta}_{\tilde{L}}^{\tilde{G}}({\bf f})$. Consid\'erons le cas $\tilde{L}=\tilde{M}$. Puisque $\boldsymbol{\delta}$ est stable, la distribution $I^{\tilde{M}}(\boldsymbol{\delta},.)$ est stable et on peut remplacer le terme ${^c\theta}_{\tilde{M}}^{\tilde{G},{\cal E}}({\bf f})$ par sa projection dans $SI(\tilde{M}(F))\otimes Mes(M(F))$. L'\'egalit\'e  3.1(2) n'est qu'une fa\c{c}on de dire que les projections dans    
 $SI(\tilde{M}(F))\otimes Mes(M(F))$ de ${^c\theta}_{\tilde{M}}^{\tilde{G},{\cal E}}({\bf f})$ et de ${^c\theta}_{\tilde{M}}^{\tilde{G}}({\bf f})$ sont \'egales.  On a ainsi obtenu le m\^eme r\'esultat que dans le cas $\tilde{L}\not=\tilde{M}$, \`a savoir que l'on peut supprimer les exposants ${\cal E}$ figurant dans la formule (2). Alors cette formule devient simplement 1.9(3).
 
 Excluons  maintenant le cas o\`u $(G,\tilde{G},{\bf a})$ est quasi-d\'eploy\'e et \`a torsion int\'erieure et o\`u ${\bf M}'={\bf M}$. On utilise la d\'efinition
 $$^cI_{\tilde{M}}^{\tilde{G},{\cal E}}({\bf M}',\boldsymbol{\delta},{\bf f})=\sum_{\tilde{s}\in \tilde{\zeta}Z(\hat{M})^{\Gamma_{F},\hat{\theta}}/Z(\hat{G})^{\Gamma_{F},\hat{\theta}}}i_{\tilde{M}'}(\tilde{G},\tilde{G}'(\tilde{s})){^cS}_{{\bf M}'}^{{\bf G}'(\tilde{s})}(\boldsymbol{\delta},{\bf f}^{{\bf G}'(\tilde{s})}).$$
 On a $dim(G'(\tilde{s})_{SC})<dim(G_{SC})$ pour tous les $\tilde{s}$ intervenant et on peut utiliser par r\'ecurrence le (ii) de la pr\'esente proposition. En utilisant les m\^emes notations que dans la preuve de la proposition [II] 2.6 (dont la pr\'esente preuve est une variante), on obtient
  $$^cI_{\tilde{M}}^{\tilde{G},{\cal E}}({\bf M}',\boldsymbol{\delta},{\bf f})=\sum_{\tilde{s}\in \tilde{\zeta}Z(\hat{M})^{\Gamma_{F},\hat{\theta}}/Z(\hat{G})^{\Gamma_{F},\hat{\theta}}}i_{\tilde{M}'}(\tilde{G},\tilde{G}'(\tilde{s}))\sum_{\tilde{L}'_{\tilde{s}}\in {\cal L}^{\tilde{G}'(\tilde{s})}(\tilde{M}')}S_{{\bf M}'}^{{\bf L}'_{\tilde{s}}}(\boldsymbol{\delta},{^cS}\theta_{{\bf L}'_{\tilde{s}}}^{{\bf G}'(\tilde{s})}({\bf f}^{{\bf G}'(\tilde{s})})).$$
  On regroupe les paires $(\tilde{s},\tilde{L}'_{\tilde{s}})$ selon l'espace de Levi $\tilde{L}$ de $\tilde{G}$ d\'etermin\'e par l'\'egalit\'e ${\cal A}_{\tilde{L}}={\cal A}_{\tilde{L}'_{\tilde{s}}}$. On obtient
 $$^cI_{\tilde{M}}^{\tilde{G},{\cal E}}({\bf M}',\boldsymbol{\delta},{\bf f})=\sum_{\tilde{L}\in {\cal L}(\tilde{M})} \sum_{\tilde{s}\in \tilde{\zeta}Z(\hat{M})^{\Gamma_{F},\hat{\theta}}/Z(\hat{L})^{\Gamma_{F},\hat{\theta}},{\bf L}'(\tilde{s})\text{ elliptique}}$$
 $$\sum_{\tilde{t}\in \tilde{s}Z(\hat{L})^{\Gamma_{F},\hat{\theta}}/Z(\hat{G})^{\Gamma_{F},\hat{\theta}}}i_{\tilde{M}'}(\tilde{G},\tilde{G}'(\tilde{t})) S_{{\bf M}'}^{{\bf L}'(\tilde{s})}(\boldsymbol{\delta},{^cS}\theta_{{\bf L}'(\tilde{s})}^{{\bf G}'(\tilde{t})}({\bf f}^{{\bf G}'(\tilde{t})})).$$
 On a l'\'egalit\'e
 $$i_{\tilde{M}'}(\tilde{G},\tilde{G}'(\tilde{t}))=i_{\tilde{M}'}(\tilde{L},\tilde{L}'(\tilde{s}))i_{\tilde{L}'(\tilde{s})}(\tilde{G},\tilde{G}'(\tilde{t}))$$
 et l'expression ci-dessus devient
 $$^cI_{\tilde{M}}^{\tilde{G},{\cal E}}({\bf M}',\boldsymbol{\delta},{\bf f})=\sum_{\tilde{L}\in {\cal L}(\tilde{M})} \sum_{\tilde{s}\in \tilde{\zeta}Z(\hat{M})^{\Gamma_{F},\hat{\theta}}/Z(\hat{L})^{\Gamma_{F},\hat{\theta}} }i_{\tilde{M}'}(\tilde{L},\tilde{L}'(\tilde{s}))$$
 $$\sum_{\tilde{t}\in \tilde{s}Z(\hat{L})^{\Gamma_{F},\hat{\theta}}/Z(\hat{G})^{\Gamma_{F},\hat{\theta}}}i_{\tilde{L}'(\tilde{s})}(\tilde{G},\tilde{G}'(\tilde{t})) S_{{\bf M}'}^{{\bf L}'(\tilde{s})}(\boldsymbol{\delta},{^cS}\theta_{{\bf L}'(\tilde{s})}^{{\bf G}'(\tilde{t})}({\bf f}^{{\bf G}'(\tilde{t})})).$$
 Pour tous $\tilde{L}$ et  $\tilde{s}$, on reconna\^{\i}t
 $$\sum_{\tilde{t}\in \tilde{s}Z(\hat{L})^{\Gamma_{F},\hat{\theta}}/Z(\hat{G})^{\Gamma_{F},\hat{\theta}}}i_{\tilde{L}'(\tilde{s})}(\tilde{G},\tilde{G}'(\tilde{t})){^cS}\theta_{{\bf L}'(\tilde{s})}^{{\bf G}'(\tilde{t})}({\bf f}^{{\bf G}'(\tilde{t})})={^c\theta}_{\tilde{L}}^{\tilde{G},{\cal E}}({\bf L}'(\tilde{s}),{\bf f})=(^c\theta_{\tilde{L}}^{\tilde{G},{\cal E}}({\bf f}))^{{\bf L}'(\tilde{s})}.$$
 L'expression plus haut devient
 $$^cI_{\tilde{M}}^{\tilde{G},{\cal E}}({\bf M}',\boldsymbol{\delta},{\bf f})=\sum_{\tilde{L}\in {\cal L}(\tilde{M})} \sum_{\tilde{s}\in \tilde{\zeta}Z(\hat{M})^{\Gamma_{F},\hat{\theta}}/Z(\hat{L})^{\Gamma_{F},\hat{\theta}} }i_{\tilde{M}'}(\tilde{L},\tilde{L}'(\tilde{s}))S_{{\bf M}'}^{{\bf L}'(\tilde{s})}(\boldsymbol{\delta},(^c\theta_{\tilde{L}}^{\tilde{G},{\cal E}}({\bf f}))^{{\bf L}'(\tilde{s})}).$$
 Par d\'efinition,  pour tout $\tilde{L}$, la somme en $\tilde{s}$ n'est autre que
 $$I_{\tilde{M}}^{\tilde{L},{\cal E}}({\bf M}',\boldsymbol{\delta},{^c\theta}_{\tilde{L}}^{\tilde{G},{\cal E}}({\bf f})).$$
 Alors l'expression ci-dessus devient simplement la formule (1) qu'il fallait prouver. Cela prouve le (i) de l'\'enonc\'e.
 
 Preuve du (ii).  Maintenant $(G,\tilde{G},{\bf a})$ est quasi-d\'eploy\'e et \`a torsion int\'erieure. On identifie $\boldsymbol{\delta}$ \`a un \'el\'ement de $D_{g\acute{e}om}^{st}({\bf M})\otimes Mes(M(F))^*$. On a prouv\'e l'identit\'e (1) (pour ${\bf M}'={\bf M}$) par une m\'ethode sp\'ecifique. Mais on peut aussi reprendre le calcul du cas g\'en\'eral. Comme toujours, il y a un terme auquel on ne peut plus appliquer par r\'ecurrence l'\'egalit\'e du (ii) de l'\'enonc\'e. C'est celui index\'e par $s=1$. On le remplace par le membre de droite du (ii) de l'\'enonc\'e plus la diff\'erence $X$ entre le membre de gauche et ce membre de droite. Le calcul se poursuit, la seule diff\'erence est l'addition de ce terme $X$. On obtient
$$ ^cI_{\tilde{M}}^{\tilde{G},{\cal E}}({\bf M},\boldsymbol{\delta},{\bf f})=X+\sum_{\tilde{L}\in {\cal L}(\tilde{M})}I_{\tilde{M}}^{\tilde{L},{\cal E}}({\bf M},\boldsymbol{\delta},{^c\theta}_{\tilde{L}}^{\tilde{G},{\cal E}}({\bf f})).$$
Puisqu'on a d\'ej\`a prouv\'e l'\'egalit\'e analogue sans le terme $X$, cela implique $X=0$, ce que l'on devait prouver. $\square$ 

\bigskip

\subsection{Preuves des propositions 2.2 et 2.5}
On suppose $(G,\tilde{G},{\bf a})$ quasi-d\'eploy\'e et \`a torsion int\'erieure. Soit ${\bf f}\in I(\tilde{G}(F))\otimes Mes(G(F))$ dont l'image dans $SI(\tilde{G}(F))\otimes Mes(G(F))$ est nulle. On veut prouver que $^cS\theta_{\tilde{M}}^{\tilde{G}}({\bf f})=0$ et que, pour tout $\boldsymbol{\delta}\in D_{g\acute{e}om}^{st}(\tilde{M}(F))\otimes Mes(M(F))^*$, on a $^cS_{\tilde{M}}^{\tilde{G}}(\boldsymbol{\delta},{\bf f})=0$. Le lemme 2.3 et nos hypoth\`eses de r\'ecurrence assurent que $(^cS\theta_{\tilde{M}}^{\tilde{G}}({\bf f}))_{\tilde{R}}=0$ pour tout espace de Levi propre $\tilde{R}$ de $\tilde{M}$, autrement dit $^cS\theta_{\tilde{M}}^{\tilde{G}}({\bf f})$ est cuspidale.
Pour $\tilde{L}\in {\cal L}(\tilde{M})$, avec $\tilde{L}\not=\tilde{M}$, nos hypoth\`eses de r\'ecurrence assurent que $^cS\theta_{\tilde{L}}^{\tilde{G}}({\bf f})=0$. Pour $\boldsymbol{\delta}$ comme ci-dessus, la formule de la proposition 4.1(ii) se simplifie en
$$(1) \qquad ^cS_{\tilde{M}}^{\tilde{G}}(\boldsymbol{\delta},{\bf f})=S^{\tilde{M}}(\boldsymbol{\delta},{^cS}\theta_{\tilde{M}}^{\tilde{G}}({\bf f})).$$
On utilise la propri\'et\'e de compacit\'e 1.9(2) dont on a dit qu'elle se propageait \`a la distribution de gauche ci-dessus: le terme est nul quand le support de $\boldsymbol{\delta}$ ne coupe pas   $\Omega^{G}$ pour un certain sous-ensemble compact $\Omega$ de $\tilde{G}(F)$. Il en r\'esulte que ${^cS}\theta_{\tilde{M}}^{\tilde{G}}({\bf f})$, qui est a priori un \'el\'ement de $SI_{ac}(\tilde{M}(F))\otimes Mes(M(F))$, est "\`a support compact", c'est-\`a-dire appartient \`a $SI(\tilde{M}(F))\otimes Mes(M(F))$. Il en r\'esulte que, pour $\tilde{\pi}\in D_{ell}^{st}(\tilde{M}(F))\otimes Mes(M(F))^*$, la fonction $\tilde{\lambda}\mapsto S^{\tilde{M}}(\tilde{\pi},\tilde{\lambda},{^cS}\theta_{\tilde{M}}^{\tilde{G}}({\bf f}))$ n'a pas de p\^ole. Ses coefficients de Fourier $S^{\tilde{M}}(\tilde{\pi},\nu,X,{^cS}\theta_{\tilde{M}}^{\tilde{G}}({\bf f}))$ sont donc ind\'ependants du point $\nu\in {\cal A}_{\tilde{M}}^*$. Dans la formule de la proposition 2.4, on peut remplacer chaque $\nu_{\tilde{S}}$ par $0$. Puisque $\sum_{\tilde{S}\in {\cal P}(\tilde{M})}\omega_{\tilde{S}}(X)=1$ pour tout $X$, on obtient que $S^{\tilde{M}}(\tilde{\pi},0,X,{^cS}\theta_{\tilde{M}}^{\tilde{G}}({\bf f}))=0$ pour tout $X$. Par inversion de Fourier, 
$S^{\tilde{M}}(\tilde{\pi},{^cS}\theta_{\tilde{M}}^{\tilde{G}}({\bf f}))=0$. Un \'el\'ement cuspidal de $SI(\tilde{M}(F))\otimes Mes(M(F))$ annul\'e par les distributions $S^{\tilde{M}}(\tilde{\pi},.)$ pour tout $\tilde{\pi}\in D_{ell}^{st}(\tilde{M}(F))\otimes Mes(M(F))^*$ est nul ([M] corollaire 4.2(i)). Donc ${^cS}\theta_{\tilde{M}}^{\tilde{G}}({\bf f})=0$, ce qui prouve la proposition 2.2. Le membre de droite de (1) est maintenant nul, donc aussi celui de gauche. Cela  prouve la proposition 2.5. 

\bigskip

\subsection{Preuve conditionnelle des propositions 3.8 et 3.9}
Ici $(G,\tilde{G},{\bf a})$ est quelconque mais on pose l'hypoth\`ese suivante:

(1) pour tout $\boldsymbol{\gamma}\in D_{g\acute{e}om,\tilde{G}-reg}(\tilde{M}(F),\omega)\otimes Mes(M(F))^*$ et tout ${\bf f}\in I(\tilde{G}(F),\omega)\otimes Mes(G(F))$, on a l'\'egalit\'e
 $$I^{\tilde{G},{\cal E}}_{\tilde{M}}(\boldsymbol{\gamma},{\bf f})=I^{\tilde{G}}_{\tilde{M}}(\boldsymbol{\gamma},{\bf f}).$$
 
 Soient $\boldsymbol{\gamma}$ et ${\bf f}$ comme ci-dessus. Consid\'erons  le membre de droite de la formule du (i) de la proposition 4.1. Pour $\tilde{L}\in {\cal L}(\tilde{M})$, la distribution $I_{\tilde{M}}^{\tilde{L},{\cal E}}(\boldsymbol{\gamma},.)$ est \'egale \`a $I_{\tilde{M}}^{\tilde{L}}(\boldsymbol{\gamma},.)$ soit par hypoth\`ese de r\'ecurrence si $\tilde{L}\not=\tilde{G}$, soit par l'hypoth\`ese (1) si $\tilde{L}=\tilde{G}$. Si $\tilde{L}\not=\tilde{M}$, on a aussi ${^c\theta}_{\tilde{L}}^{\tilde{G},{\cal E}}({\bf f})={^c\theta}_{\tilde{L}}^{\tilde{G}}({\bf f})$ par nos hypoth\`eses de r\'ecurrence. La formule devient
 $$^cI_{\tilde{M}}^{\tilde{G},{\cal E}}(\boldsymbol{\gamma},{\bf f})=I^{\tilde{M}}(\boldsymbol{\gamma},{^c\theta}_{\tilde{M}}^{\tilde{G},{\cal E}}({\bf f}))+\sum_{\tilde{L}\in {\cal L}(\tilde{M}),\tilde{L}\not=\tilde{M}}I_{\tilde{M}}^{\tilde{L}}(\boldsymbol{\gamma},{^c\theta}_{\tilde{L}}^{\tilde{G}}({\bf f})).$$
 En comparant avec 1.9(3), on obtient
 $$(2) \qquad ^cI_{\tilde{M}}^{\tilde{G}}(\boldsymbol{\gamma},{\bf f})-{^cI}_{\tilde{M}}^{\tilde{G},{\cal E}}(\boldsymbol{\gamma},{\bf f})=I^{\tilde{M}}(\boldsymbol{\gamma},\boldsymbol{\varphi}),$$
 o\`u
 $$\boldsymbol{\varphi}={^c\theta}_{\tilde{M}}^{\tilde{G}}({\bf f})-{^c\theta}_{\tilde{M}}^{\tilde{G},{\cal E}}({\bf f}).$$
 La propri\'et\'e de compacit\'e 1.9(3) se propage comme on l'a dit au membre de gauche de (2). Il en r\'esulte que $\boldsymbol{\varphi}$, qui est a priori un \'el\'ement de $I_{ac}(\tilde{M}(F),\omega)\otimes Mes(M(F))$, est "\`a support compact", c'est-\`a-dire est un \'el\'ement de $I(\tilde{M}(F),\omega)\otimes Mes(M(F))$. Pour une $\omega$-repr\'esentation elliptique $\tilde{\pi}$ de $\tilde{M}(F)$, la fonction $\tilde{\lambda}\mapsto I^{\tilde{M}}(\tilde{\pi},\tilde{\lambda},\boldsymbol{\varphi})$ n'a donc pas de p\^ole. Ses coefficients de Fourier $I^{\tilde{M}}(\tilde{\pi},\nu,X,\boldsymbol{\varphi})$ ne d\'ependent donc pas du point $\nu\in {\cal A}_{\tilde{M}}^*$. Par construction, on a
 $$I^{\tilde{M}}(\tilde{\pi},\nu,X,\boldsymbol{\varphi})=I^{\tilde{M}}(\tilde{\pi},\nu,X,{^c\theta}_{\tilde{M}}^{\tilde{G}}({\bf f}))-I^{\tilde{M}}(\tilde{\pi},\nu,X,{^c\theta}_{\tilde{M}}^{\tilde{G},{\cal E}}({\bf f})).$$
 On utilise les propositions 1.8 et 3.7 qui nous disent que
 $$\sum_{\tilde{S}\in {\cal P}(\tilde{M})}\omega_{\tilde{S}}(X)I^{\tilde{M}}(\tilde{\pi},\nu_{\tilde{S}},X,\boldsymbol{\varphi})=0$$
 pour tout $X$. Comme dans le paragraphe pr\'ec\'edent, on peut remplacer les points $\nu_{\tilde{S}}$ par $0$ et conclure que $I^{\tilde{M}}(\tilde{\pi},0,X,\boldsymbol{\varphi})=0$ pour tout $X$. Par inversion de Fourier, on obtient $I^{\tilde{M}}(\tilde{\pi},\boldsymbol{\varphi})=0$. Par ailleurs, les lemmes 1.6 et  3.5 et nos hypoth\`eses de r\'ecurrence impliquent que $\boldsymbol{\varphi}_{\tilde{R}}=0$ pour tout espace de Levi propre $\tilde{R}$ de $\tilde{M}$, c'est-\`a-dire que $\boldsymbol{\varphi}$ est cuspidale. Etant annul\'ee par la distribution $I^{\tilde{M}}(\tilde{\pi},.)$ pour toute $\omega$-repr\'esentation elliptique $\tilde{\pi}$ de $\tilde{M}(F)$, cette fonction est nulle. Donc 
 $${^c\theta}_{\tilde{M}}^{\tilde{G}}({\bf f})={^c\theta}_{\tilde{M}}^{\tilde{G},{\cal E}}({\bf f}),$$
 ce qui prouve la proposition 3.8.
 
 La fonction $\boldsymbol{\varphi}$ \'etant nulle, l'\'egalit\'e (2) entra\^{\i}ne
 $$^cI_{\tilde{M}}^{\tilde{G}}(\boldsymbol{\gamma},{\bf f})={^cI}_{\tilde{M}}^{\tilde{G},{\cal E}}(\boldsymbol{\gamma},{\bf f}),$$
 ce qui prouve la proposition 3.9.
 
 Dans le cas o\`u $(G,\tilde{G},{\bf a})$ est quasi-d\'eploy\'e et \`a torsion int\'erieure, l'hypoth\`ese (1) est le th\'eor\`eme [II] 1.16(ii), que l'on a d\'ej\`a prouv\'e. Donc les propositions 3.8 et 3.9 sont prouv\'ees inconditionnellement dans ce cas.  Comme on l'a dit, la restriction que $\boldsymbol{\delta}$ est \`a support fortement r\'egulier dans $\tilde{G}$ ne sert qu'\`a \'eviter dans le  cas g\'en\'eral l'introduction du d\'elicat syst\`eme de fonctions $B^{\tilde{G}}$. Dans le cas quasi-d\'eploy\'e et \`a torsion int\'erieure, ce syst\`eme de fonctions est trivial, la preuve vaut aussi bien pour $\boldsymbol{\delta}$ quelconque. 
 
 \bigskip
 
 \subsection{L'application $\epsilon_{\tilde{M}}$}
 \bigskip
On a d\'efini en [III] 6.3 des triplets particuliers $(G,\tilde{G},{\bf a})$.  Consid\'erons un tel triplet. Le groupe $G$ est simplement connexe et quasi-d\'eploy\'e. On note $\Theta_{F}$ l'ensemble des \'el\'ements  semi-simples $\eta\in \tilde{G}(F)$  tels que $ad_{\eta}$ conserve une paire de Borel \'epingl\'ee d\'efinie sur $F$. On a vu en [III] 6.3 que ces \'el\'ements sont contenus dans un nombre fini de classes de conjugaison stable. 
 On d\'efinit l'espace $I(\tilde{G}(F),\omega)^{00}$ comme le sous-espace des \'el\'ements ${\bf f}\in I(\tilde{G}(F),\omega)$ dont les int\'egrales orbitales sont nulles en tout point $\gamma\in \tilde{G}(F)$ dont la partie semi-simple appartient \`a une classe de conjugaison stable coupant $ \Theta_{F}$.

Si $(G,\tilde{G},{\bf a})$ n'est pas l'un des triplets particuliers d\'efinis en [III] 6.3, on pose simplement  $I(\tilde{G}(F),\omega)^{00}=I(\tilde{G}(F),\omega)$.

On note naturellement $I_{ac,cusp}(\tilde{M}(F),\omega)$ le sous-espace des \'el\'ements cuspidaux de $I_{ac}(\tilde{M}(F),\omega)$. 

 \ass{Proposition}{Il existe une unique application lin\'eaire
 $$\epsilon_{\tilde{M}}:I(\tilde{G}(F),\omega)^{00}\otimes Mes(G(F))\to I_{ac,cusp}(\tilde{M}(F),\omega)\otimes Mes(M(F))$$
 telle que
 $$I_{\tilde{M}}^{\tilde{G},{\cal E}}(\boldsymbol{\gamma},{\bf f})-I_{\tilde{M}}^{\tilde{G}}(\boldsymbol{\gamma},{\bf f})=I^{\tilde{M}}(\boldsymbol{\gamma},\epsilon_{\tilde{M}}({\bf f}))$$
 pour tout $\boldsymbol{\gamma}\in D_{g\acute{e}om}(\tilde{M}(F),\omega)\otimes Mes(M(F))^*$ et tout ${\bf f}\in I(\tilde{G}(F),\omega)^{00}\otimes Mes(G(F))$.}
 
 Preuve. Il est clair que l'application $\epsilon_{\tilde{M}}$ est unique si elle existe. Fixons ${\bf f}\in 
 I(\tilde{G}(F),\omega)^{00}\otimes Mes(G(F))$. Consid\'erons l'\'egalit\'e
 
 (1) $I_{\tilde{M}}^{\tilde{G},{\cal E}}(\boldsymbol{\gamma},{\bf f})-I_{\tilde{M}}^{\tilde{G}}(\boldsymbol{\gamma},{\bf f})=I^{\tilde{M}}(\boldsymbol{\gamma},\boldsymbol{\varphi}).$
 
 Commen\c{c}ons par prouver
 
 (2) pour toute classe de conjugaison stable semi-simple ${\cal O}$ dans $\tilde{M}(F)$, il existe $\boldsymbol{\varphi}\in I(\tilde{M}(F),\omega)\otimes Mes(M(F))$ de sorte que
  (1) soit v\'erifi\'ee 
 pour tout $\boldsymbol{\gamma}\in D_{g\acute{e}om,\tilde{G}-\acute{e}qui}(\tilde{M}(F),\omega)\otimes Mes(M(F))^*$ assez proche de ${\cal O}$.
 
 On rappelle qu'on dit que $\boldsymbol{\gamma}$ est assez proche de ${\cal O}$ si les parties semi-simples des \'el\'ements du support de $\boldsymbol{\gamma}$ sont assez proches de ${\cal O}$.  On utilise les d\'eveloppements en germes des deux termes du membre de gauche fournis par les propositions [II] 2.3 et [II] 2.6. On a ainsi
  $$I_{\tilde{M}}^{\tilde{G},{\cal E}}(\boldsymbol{\gamma},{\bf f})-I_{\tilde{M}}^{\tilde{G}}(\boldsymbol{\gamma},{\bf f})=\sum_{\tilde{L}\in {\cal L}(\tilde{M})}I_{\tilde{L}}^{\tilde{G},{\cal E}}(g_{\tilde{M},{\cal O}}^{\tilde{L},{\cal E}}(\boldsymbol{\gamma}),{\bf f})-I_{\tilde{L}}^{\tilde{G}}(g_{\tilde{M},{\cal O}}^{\tilde{L}}(\boldsymbol{\gamma}),{\bf f}).$$
  Si $\tilde{L}\not=\tilde{M}$, on applique par r\'ecurrence le th\'eor\`eme [II] 1.16(i): $I_{\tilde{L}}^{\tilde{G},{\cal E}}=I_{\tilde{L}}^{\tilde{G}}$. Si $\tilde{L}\not=\tilde{G}$, on utilise par r\'ecurrence la proposition [II] 2.7: $g_{\tilde{M},{\cal O}}^{\tilde{L},{\cal E}}=g_{\tilde{M},{\cal O}}^{\tilde{L}}$. La formule se simplifie en
 $$I_{\tilde{M}}^{\tilde{G},{\cal E}}(\boldsymbol{\gamma},{\bf f})-I_{\tilde{M}}^{\tilde{G}}(\boldsymbol{\gamma},{\bf f})=I^{\tilde{G}}(g_{\tilde{M},{\cal O}}^{\tilde{G},{\cal E}}(\boldsymbol{\gamma})-g_{\tilde{M},{\cal O}}^{\tilde{G}}(\boldsymbol{\gamma}),{\bf f})$$
 $$+I_{\tilde{M}}^{\tilde{G},{\cal E}}(g_{\tilde{M},{\cal O}}^{\tilde{M}}(\boldsymbol{\gamma}),{\bf f})-  I_{\tilde{M}}^{\tilde{G}}(g_{\tilde{M},{\cal O}}^{\tilde{M}}(\boldsymbol{\gamma}),{\bf f}).$$
Les termes $g_{\tilde{M},{\cal O}}^{\tilde{G},{\cal E}}(\boldsymbol{\gamma})$ et $g_{\tilde{M},{\cal O}}^{\tilde{G}}(\boldsymbol{\gamma})$ sont en tout cas des \'el\'ements de $D_{g\acute{e}om}({\cal O}^{\tilde{G}},\omega)\otimes Mes(G(F))^*$ o\`u ${\cal O}^{\tilde{G}}$ est la classe de conjugaison stable dans $\tilde{G}(F)$ contenant ${\cal O}$.  Si $(G,\tilde{G},{\bf a})$ est l'un des triplets construits en [III] 6.3  et si ${\cal O}^{\tilde{G}}$ coupe $\Theta_{F}$,   l'hypoth\`ese que  ${\bf f}\in I(\tilde{G}(F),\omega)^{00}\otimes Mes(G(F))$ entra\^{\i}ne que le premier terme du membre de droite ci-dessus est nul. Si au contraire  $(G,\tilde{G},{\bf a})$ n'est pas l'un de ces triplets  ou si c'est l'un d'eux mais   ${\cal O}^{\tilde{G}}$ ne coupe pas  $\Theta_{F}$, la proposition [III] 8.5 affirme que $g_{\tilde{M},{\cal O}}^{\tilde{G},{\cal E}}(\boldsymbol{\gamma})=g_{\tilde{M},{\cal O}}^{\tilde{G}}(\boldsymbol{\gamma})$. Le premier terme ci-dessus  est donc encore nul. L'\'egalit\'e se simplifie en
 $$(3) \qquad I_{\tilde{M}}^{\tilde{G},{\cal E}}(\boldsymbol{\gamma},{\bf f})-I_{\tilde{M}}^{\tilde{G}}(\boldsymbol{\gamma},{\bf f})=I_{\tilde{M}}^{\tilde{G},{\cal E}}(g_{\tilde{M},{\cal O}}^{\tilde{M}}(\boldsymbol{\gamma}),{\bf f})-  I_{\tilde{M}}^{\tilde{G}}(g_{\tilde{M},{\cal O}}^{\tilde{M}}(\boldsymbol{\gamma}),{\bf f}).$$
 Choisissons une fonction $\boldsymbol{\varphi}\in I(\tilde{M}(F),\omega)\otimes Mes(M(F))$ telle que
 $$I^{\tilde{M}}(\boldsymbol{\tau},\boldsymbol{\varphi})=I_{\tilde{M}}^{\tilde{G},{\cal E}}(\boldsymbol{\tau},{\bf f})-I_{\tilde{M}}^{\tilde{G}}(\boldsymbol{\tau},{\bf f})$$
 pour tout $\boldsymbol{\tau}\in D_{g\acute{e}om}({\cal O},\omega)\otimes Mes(M(F))^*$. C'est possible  puisque $D_{g\acute{e}om}({\cal O},\omega)\otimes Mes(M(F))^*$ est un sous-espace de dimension finie du dual de $I(\tilde{M}(F),\omega)\otimes Mes(M(F))$. Cela entra\^{\i}ne en particulier que
 $$I^{\tilde{M}}(g_{\tilde{M},{\cal O}}^{\tilde{M}}(\boldsymbol{\gamma}),\boldsymbol{\varphi})=I_{\tilde{M}}^{\tilde{G},{\cal E}}(g_{\tilde{M},{\cal O}}^{\tilde{M}}(\boldsymbol{\gamma}),{\bf f})-  I_{\tilde{M}}^{\tilde{G}}(g_{\tilde{M},{\cal O}}^{\tilde{M}}(\boldsymbol{\gamma}),{\bf f}).$$
 Par le d\'eveloppement en germes de Shalika ordinaires, le membre de gauche ci-dessus est \'egal \`a $I^{\tilde{M}}(\boldsymbol{\gamma},\boldsymbol{\varphi})$ pourvu que $\boldsymbol{\gamma}$ soit assez proche de ${\cal O}$. Alors l'\'egalit\'e (3) devient (1).  Cela d\'emontre  l'assertion (2).
 
 Prouvons que
 
 (4) il existe un sous-ensemble compact $\Omega\subset \tilde{M}(F)$ et une fonction $\boldsymbol{\varphi}\in I_{ac}(\tilde{M}(F),\omega)\otimes Mes(M(F))$ de sorte que l'on ait l'\'egalit\'e
   (1)
 pour tout $\boldsymbol{\gamma}\in D_{g\acute{e}om}(\tilde{M}(F),\omega)\otimes Mes(M(F))^*$
 dont le support ne coupe pas $\Omega^{M}$.
 
 On utilise la relation 1.9(3) et la proposition 4.1(i). On obtient l'\'egalit\'e
 $$I_{\tilde{M}}^{\tilde{G},{\cal E}}(\boldsymbol{\gamma},{\bf f})-I_{\tilde{M}}^{\tilde{G}}(\boldsymbol{\gamma},{\bf f})={^cI}_{\tilde{M}}^{\tilde{G},{\cal E}}(\boldsymbol{\gamma},{\bf f})-{^cI}_{\tilde{M}}^{\tilde{G}}(\boldsymbol{\gamma},{\bf f})$$
 $$-\sum_{\tilde{L}\in {\cal L}(\tilde{M}), \tilde{L}\not=\tilde{G}}I_{\tilde{M}}^{\tilde{L},{\cal E}}(\boldsymbol{\gamma},{^c\theta}_{\tilde{L}}^{\tilde{G},{\cal E}}({\bf f}))-I_{\tilde{M}}^{\tilde{L}}(\boldsymbol{\gamma},{^c\theta}_{\tilde{L}}^{\tilde{G}}({\bf f})).$$
 On peut utiliser le th\'eor\`eme [II] 1.16(i) par r\'ecurrence puisque $\tilde{L}\not=\tilde{G}$:  on a $I_{\tilde{M}}^{\tilde{L},{\cal E}}=I_{\tilde{M}}^{\tilde{L}}$. Si $\tilde{L}\not=\tilde{M}$, on peut aussi utiliser la proposition 3.8: $^c\theta_{\tilde{L}}^{\tilde{G},{\cal E}}({\bf f})={^c\theta}_{\tilde{L}}^{\tilde{G}}({\bf f})$. L'\'egalit\'e se simplifie en
  $$(5) \qquad I_{\tilde{M}}^{\tilde{G},{\cal E}}(\boldsymbol{\gamma},{\bf f})-I_{\tilde{M}}^{\tilde{G}}(\boldsymbol{\gamma},{\bf f})={^cI}_{\tilde{M}}^{\tilde{G},{\cal E}}(\boldsymbol{\gamma},{\bf f})-{^cI}_{\tilde{M}}^{\tilde{G}}(\boldsymbol{\gamma},{\bf f})+I^{\tilde{M}}(\boldsymbol{\gamma},\boldsymbol{\varphi}),$$
  o\`u 
  $$\boldsymbol{\varphi}={^c\theta}_{\tilde{M}}^{\tilde{G},{\cal E}}({\bf f})-{^c\theta}_{\tilde{M}}^{\tilde{G}}({\bf f}).$$
  Cette fonction $\boldsymbol{\varphi}$ est bien un \'el\'ement de $I_{ac}(\tilde{M}(F),\omega)\otimes Mes(M(F))$. La propri\'et\'e de compacit\'e 1.9(2), qui se propage au premier terme du membre de droite de (5), assure l'existence d'un sous-ensemble compact $\Omega$ de $\tilde{M}(F)$ tel que les deux premiers termes de ce membre de droite soient nuls si le support de $\boldsymbol{\gamma}$ ne coupe pas $\Omega^{M}$. Cela prouve (4).
 
On peut \'evidemment supposer $\Omega^M$ ouvert, ferm\'e et invariant par conjugaison stable.  Par partition de l'unit\'e, l'assertion (2) permet de construire une fonction $\boldsymbol{\varphi}\in I(\tilde{M}(F),\omega)\otimes Mes(M(F))$ telle que (1) soit v\'erifi\'ee 
 pour tout $\boldsymbol{\gamma}\in D_{g\acute{e}om,\tilde{G}-\acute{e}qui}(\tilde{M}(F),\omega)\otimes Mes(M(F))^*$ \`a support contenu dans $\Omega^M$. L'assertion (4) construit une telle fonction telle que (1) soit v\'erifi\'ee pour tout  $\boldsymbol{\gamma}\in D_{g\acute{e}om}(\tilde{M}(F),\omega)\otimes Mes(M(F))^*$ \`a support disjoint  de $\Omega^M$. En recollant ces deux fonctions, on obtient une fonction $\boldsymbol{\varphi}$ telle que (1) soit v\'erifi\'ee pour tout $\boldsymbol{\gamma}\in D_{g\acute{e}om,\tilde{G}-\acute{e}qui}(\tilde{M}(F),\omega)\otimes Mes(M(F))^*$. Montrons que l'on peut supprimer la restriction sur le support de $\boldsymbol{\gamma}$. Pour cela, il suffit de fixer une classe de conjugaison stable semi-simple ${\cal O}$ dans $\tilde{M}(F)$ et de prouver que la relation (1) est encore v\'erifi\'ee pour $\boldsymbol{\gamma}\in D_{g\acute{e}om}({\cal O},\omega)\otimes Mes(M(F))^*$. Pour de tels ${\cal O}$ et $\boldsymbol{\gamma}$, on peut choisir $\boldsymbol{\gamma}'$ \`a support fortement r\'egulier dans $\tilde{G}$ et proche de ${\cal O}$ de sorte que $g_{\tilde{M},{\cal O}}^{\tilde{M}}(\boldsymbol{\gamma}')=\boldsymbol{\gamma}$, cf. [II] 2.1(1). La relation (3) appliqu\'ee \`a ce $\boldsymbol{\gamma}'$ nous dit que
 $$ I_{\tilde{M}}^{\tilde{G},{\cal E}}(\boldsymbol{\gamma}',{\bf f})-I_{\tilde{M}}^{\tilde{G}}(\boldsymbol{\gamma}',{\bf f})= I_{\tilde{M}}^{\tilde{G},{\cal E}}(\boldsymbol{\gamma},{\bf f})-I_{\tilde{M}}^{\tilde{G}}(\boldsymbol{\gamma},{\bf f}).$$
 De m\^eme, on a
 $$I^{\tilde{M}}(\boldsymbol{\gamma}',\boldsymbol{\varphi})=I^{\tilde{M}}(\boldsymbol{\gamma},\boldsymbol{\varphi}).$$
 Puisque les membres de gauche des deux \'egalit\'es pr\'ec\'edentes sont \'egaux, ceux de droite le sont aussi.
 
 On a ainsi construit une fonction $\boldsymbol{\varphi}\in I_{ac}(\tilde{M}(F))\otimes Mes(M(F))$ pour laquelle (1) est v\'erif\'ee pour tout $\boldsymbol{\gamma}\in D_{g\acute{e}om}(\tilde{M}(F),\omega)\otimes Mes(M(F))^*$. Soit $\tilde{R}$ un espace de Levi propre de $\tilde{M}$ et $\boldsymbol{\tau}\in D_{g\acute{e}om}(\tilde{R}(F),\omega)\otimes Mes(R(F))^*$. Les relations de descente du lemme [II] 1.7 et de [II] 1.15(1), jointes \`a nos hypoth\`eses de r\'ecurrence, assurent que le membre de gauche de (1) est nul pour $\boldsymbol{\gamma}=\boldsymbol{\tau}^{\tilde{M}}$. Donc $I^{\tilde{M}}(\boldsymbol{\tau}^{\tilde{M}},\boldsymbol{\varphi})=0$. Cela assure que $\boldsymbol{\varphi}$ est cuspidale. $\square$
 
 Il r\'esulte de la preuve que, pour tout ${\bf f}\in I(\tilde{G}(F),\omega)^{00}\otimes Mes(G(F))$, la fonction $\epsilon_{\tilde{M}}({\bf f})$ est la somme de $^c\theta_{\tilde{M}}^{\tilde{G},{\cal E}}({\bf f})-{^c\theta}_{\tilde{M}}^{\tilde{G}}({\bf f})$ et d'une fonction \`a support compact, c'est-\`a-dire d'un \'el\'ement de $I(\tilde{M}(F),\omega)\otimes Mes(M(F))$. Il r\'esulte alors de 1.8 et 3.6 que
 
 (6) la fonction $\epsilon_{\tilde{M}}({\bf f})$ est de Schwartz;
 
 (7) soit $\tilde{\pi}\in D_{temp}(\tilde{M}(F),\omega)\otimes Mes(M(F))^*$; alors la fonction $\tilde{\lambda}\mapsto I^{\tilde{M}}(\tilde{\pi},\tilde{\lambda},\epsilon_{\tilde{M}}({\bf f}))$ d\'efinie sur $i\tilde{{\cal A}}_{\tilde{M}}^*/i\tilde{{\cal A}}_{\tilde{M},F}^{\vee}$ se prolonge en une fonction rationnelle sur $\tilde{{\cal A}}_{\tilde{M},{\mathbb C}}^*/i\tilde{{\cal A}}_{\tilde{M},F}^{\vee}$; ses p\^oles sont de la forme d\'ecrite en 1.8(2). 
 \bigskip
 
 {\bf Bibliographie.}
 
 [A1] J. Arthur: {\it A stable trace formula III. Proof of the main theorems }, Annals of Math. 158 (2003), p. 769-873
 
 [A2] ------------: {\it The trace formula in invariant form}, Annals of Math. 114 (1981), p. 1-74
 
 [A3] -----------: {\it Intertwining operators and residues I. Weighted characters}, J. of Funct. Analysis 84 (1989), p. 19-84
 
 [HL] G. Henniart, B. Lemaire: {\it La transform\'ee de Fourier pour les espaces tordus sur un groupe r\'eductif $p$-adique I. Le th\'eor\`eme de Paley-Wiener }, pr\'epublication 2013
 
 [M] C. Moeglin: {\it  Repr\'esentations elliptiques; caract\'erisation et formule de transfert de caract\`eres}, pr\'epublication 2013
 
 [W] J.-L. Waldspurger: {\it La formule des traces locale tordue}, pr\'epublication 2012
 
 [I], [II], [III], [IV]  ------------------------: {\it  Stabilisation de la formule des traces tordue I: endoscopie tordue sur un corps local}, {\it II:   int\'egrales orbitales et endoscopie sur un corps local non archim\'edien; d\'efinitions et \'enonc\'es des r\'esultats},  {\it III: int\'egrales orbitales et endoscopie sur un corps local non archim\'edien; r\'eductions et preuves},  {\it IV: transfert spectral archim\'edien}, pr\'epublications 2014
  
  \bigskip
  
  e-mail: jean-loup.waldspurger@imj-prg.fr
  
\end{document}